\documentclass{amsart}
\usepackage{amsthm,amsmath,amsfonts}
\usepackage{mathrsfs}
\usepackage{tocvsec2}
\usepackage{framed}
 \usepackage[colorlinks=true]{hyperref}
\hypersetup{urlcolor=blue, citecolor=red}

\numberwithin{equation}{section}
\numberwithin{table}{section}
\font\tenscrpt=eusm10
\font\sevenscrpt=eusm10 scaled 700
\font\fivescrpt=eusm10 scaled 500
\newfam\eusmfam
\textfont\eusmfam=\tenscrpt
\scriptfont\eusmfam=\sevenscrpt
\scriptscriptfont\eusmfam=\fivescrpt

\newtheorem{thm}{Theorem}[section]
\newtheorem{cor}{Corollary}[section]
\newtheorem{lem}{Lemma}[section]
\newtheorem{prop}{Proposition}[section]

\theoremstyle{definition}
\newtheorem{defn}{Definition}[section]

\newtheorem{rem}{Remark}[section]
\newtheorem{notn}{Notation}[section]
\newcommand{\thmref}[1]{Theorem~\ref{#1}}
\newcommand{\secref}[1]{Section~\ref{#1}}

\newcommand{\lemref}[1]{Lemma~\ref{#1}}
\newcommand{\coref}[1]{Corollary~\ref{#1}}
\newcommand{\propref}[1]{Proposition~\ref{#1}}
\newcommand{\defnref}[1]{Definition~\ref{#1}}
\newcommand{\remref}[1]{Remark~\ref{#1}}
\newcommand{\eqnref}[1]{{\rm (\ref{#1})}}
\newcommand{\tablref}[1]{Table~\ref{#1}}
\def\ds{\displaystyle}
\def\ts{\textstyle}

\def\h{\delta_n}

\def\hd{\delta_n^d}
\def\sdh{\delta_n^{{d}/{2}}}

\def\d{\delta}

\def\dd{\delta^d}

\def\sdd{\delta^{{d}/{2}}}

\def\ptxy{{K}^{\text{\tiny{\sc{BM}}}^d}_{t;x,y}}
\def\psxy{{K}^{\text{\tiny{\sc{BM}}}^d}_{s;x,y}}

\def\psoxz{{K}^{\text{\tiny{\sc{BM}}}^d}_{s_{1};x}}
\def\pstxz{{K}^{\text{\tiny{\sc{BM}}}^d}_{s_{2};x}}

\def\ptzs{{K}^{\text{\tiny{\sc{BM}}}}_{t;0,s}}
\def\ptzz{{K}^{\text{\tiny{\sc{BM}}}}_{t;0,0}}

\def\pszri{{K}^{\text{\tiny{\sc{BM}}}}_{s;0,r_{i}}}
\def\pszro{{K}^{\text{\tiny{\sc{BM}}}}_{s;0,r_{1}}}
\def\pszrt{{K}^{\text{\tiny{\sc{BM}}}}_{s;0,r_{2}}}
\def\ptzso{{K}^{\text{\tiny{\sc{BM}}}}_{t;0,s_{1}}}
\def\ptzst{{K}^{\text{\tiny{\sc{BM}}}}_{t;0,s_{2}}}
\def\puzro{{K}^{\text{\tiny{\sc{BM}}}}_{u;0,r_{1}}}
\def\puzrt{{K}^{\text{\tiny{\sc{BM}}}}_{u;0,r_{2}}}
\def\pvzrt{{K}^{\text{\tiny{\sc{BM}}}}_{v;0,r_{2}}}
\def\proprtz{K^{{\text{\tiny{\sc{BM}}}^d}}_{r_{1}+r_{2};0}}
\def\prox{K^{{\text{\tiny{\sc{BM}}}^d}}_{r_{1};x}}
\def\prtx{K^{{\text{\tiny{\sc{BM}}}^d}}_{r_{2};x}}
\def\prix{K^{{\text{\tiny{\sc{BM}}}^d}}_{r_{i};x}}
\def\prixpz{K^{{\text{\tiny{\sc{BM}}}^d}}_{r_{i};x+z}}
\def\proprtz{K^{{\text{\tiny{\sc{BM}}}^d}}_{r_{1}+r_{2};0}}
\def\proprtzz{K^{{\text{\tiny{\sc{BM}}}^d}}_{r_{1}+r_{2};z}}
\def\Mpisxy {{K}^{\text{\tiny{\sc{CLKS}}}^d}_{is;x,y}}
\def\K{{\mathbb K}}
\def\BTBM{{\K}^{\text{\tiny{\sc{BTBM}}}^d}}
\def\KBtxy{{\K}^{\text{\tiny{\sc{BTBM}}}^d}_{t;x,y}}
\def\KBtx{{\K}^{\text{\tiny{\sc{BTBM}}}^d}_{t;x}}
\def\KKStxy{{\K}^{\text{\tiny{\sc{LKS}}}^d}_{t;x,y}}

\def\KBtsxy{{\K}^{\text{\tiny{\sc{BTBM}}}^d}_{t-s;x,y}}
\def\KBtsx{{\K}^{\text{\tiny{\sc{BTBM}}}^d}_{t-s;x}}
\def\KBrsx{{\K}^{\text{\tiny{\sc{BTBM}}}^d}_{r-s;x}}
\def\KBsptmrx{{\K}^{{\text{\tiny{\sc{BTBM}}}^d}}_{s+(t-r);x}}
\def\KBsx{{\K}^{{\text{\tiny{\sc{BTBM}}}^d}}_{s;x}}
\def\KBsxpz{{\K}^{{\text{\tiny{\sc{BTBM}}}^d}}_{s;x+z}}
\def\KBux{{\K}^{{\text{\tiny{\sc{BTBM}}}^d}}_{u;x}}
\def\KKStsxy{{\K}^{\text{\tiny{\sc{LKS}}}^d}_{t-s;x,y}}
\def\BTRW{{\K}^{{\text{\tiny{\sc{BTRW}}}^d_\delta}}}
\def\Ktx{{\K}^{{\text{\tiny{\sc{BTRW}}}^d_\delta}}_{t;x}}
\def\Ktxy{{\K}^{{\text{\tiny{\sc{BTRW}}}^d_{\delta}}}_{t;x,y}}
\def\Ktsx{{\K}^{{\text{\tiny{\sc{BTRW}}}^d_\delta}}_{t-s;x}}
\def\Krsx{{\K}^{{\text{\tiny{\sc{BTRW}}}^d_\delta}}_{r-s;x}}
\def\Ksx{{\K}^{{\text{\tiny{\sc{BTRW}}}^d_\delta}}_{s;x}}
\def\Kux{{\K}^{{\text{\tiny{\sc{BTRW}}}^d_\delta}}_{u;x}}
\def\Kvx{{\K}^{{\text{\tiny{\sc{BTRW}}}^d_\delta}}_{v;x}}
\def\Ksptmrx{{\K}^{{\text{\tiny{\sc{BTRW}}}^d_\delta}}_{s+(t-r);x}}

\def\Ktsxy{{\K}^{{\text{\tiny{\sc{BTRW}}}^d_\delta}}_{t-s;x,y}}
\def\Ktsxz{{\K}^{{\text{\tiny{\sc{BTRW}}}^d_\delta}}_{t-s;x,z}}
\def\Ktsyz{{\K}^{{\text{\tiny{\sc{BTRW}}}^d_\delta}}_{t-s;y,z}}
\def\Krsxz{{\K}^{{\text{\tiny{\sc{BTRW}}}^d_\delta}}_{r-s;x,z}}
\def\Kroxz{{\K}^{{\text{\tiny{\sc{BTRW}}}^d_\delta}}_{\rho;x,z}}
\def\Ksx{{\K}^{{\text{\tiny{\sc{BTRW}}}^d_\delta}}_{s;x}}
\def\Ksxpz{{\K}^{{\text{\tiny{\sc{BTRW}}}^d_\delta}}_{s;x+z}}
\def\Kzxn{{\K}^{{\text{\tiny{\sc{BTRW}}}^d_{\h}}}_{0;x}}
\def\Kzxyn{{\K}^{{\text{\tiny{\sc{BTRW}}}^d_{\h}}}_{0;x,y}}
\def\Ktxn{{\K}^{{\text{\tiny{\sc{BTRW}}}^d_{\h}}}_{t,x}}
\def\Ktxzn{{\K}^{{\text{\tiny{\sc{BTRW}}}^d_{\h}}}_{t;x,0}}
\def\Ktsxyn{{\K}^{{\text{\tiny{\sc{BTRW}}}^d_{\h}}}_{t-s;x,y}}
\def\Ktsxykn{{\K}^{{\text{\tiny{\sc{BTRW}}}^d_{\h}}}_{t-s;x,y^{(k)}}}
\def\Ktausxyn{{\K}^{{\text{\tiny{\sc{BTRW}}}^d_{\h}}}_{\tau-s;x,y}}

\def\Ktauitxiyn{{\K}^{{\text{\tiny{\sc{BTRW}}}^d_{\h}}}_{\tau_i-t;x_i,y}}

\def\Ktaujtxjyn{{\K}^{{\text{\tiny{\sc{BTRW}}}^d_{\h}}}_{\tau_j-t;x_j,y}}
\def\Ktxyn{{\K}^{{\text{\tiny{\sc{BTRW}}}^d_{\h}}}_{t;x,y}}
\def\Ktxynint{{\K}^{{\text{\tiny{\sc{BTRW}}}^d_{\h}}}_{t;[x]_{\h},[y]_{\h}}}

\def\KBBsptmrz{{\K}^{{\text{\tiny{\sc{BTRW}}}^d_\delta(2)}}_{s+(t-r),s+(t-r);0}}
\def\KBBBsptmrz{{\K}^{{\text{\tiny{\sc{BTBM}}}^d(2)}}_{s+(t-r),s+(t-r);0}}
\def\KBBsptmrsz{{\K}^{{\text{\tiny{\sc{BTRW}}}^d_\delta(2)}}_{s+(t-r),s;0}}
\def\KBBBsptmrsz{{\K}^{{\text{\tiny{\sc{BTBM}}}^d(2)}}_{s+(t-r),s;0}}
\def\KBBsz{{\K}^{{\text{\tiny{\sc{BTRW}}}^d_\delta(2)}}_{s,s;0}}
\def\KBBBsz{{\K}^{{\text{\tiny{\sc{BTBM}}}^d(2)}}_{s,s;0}}

\def\KBBBuz{{\K}^{{\text{\tiny{\sc{BTBM}}}^d(2)}}_{u,u;0}}
\def\KBBuvz{{\K}^{{\text{\tiny{\sc{BTRW}}}^d_\delta(2)}}_{u,v;0}}
\def\KBBBuvz{{\K}^{{\text{\tiny{\sc{BTBM}}}^d(2)}}_{u,v;0}}

\def\qsxyn{K^{{\text{\tiny{\sc{RW}}}}_{\h}^d}_{s;x,y}}

\def\qzxn{K^{{\text{\tiny{\sc{RW}}}}_{\h}^d}_{0;x}}
\def\qtxn{K^{{\text{\tiny{\sc{RW}}}}_{\h}^d}_{t;x}}
\def\qsxn{K^{{\text{\tiny{\sc{RW}}}}_{\h}^d}_{s;x}}
\def\qtxyn{K^{{\text{\tiny{\sc{RW}}}}_{\h}^d}_{t;x,y}}
\def\qtxynint{K^{{\text{\tiny{\sc{RW}}}}_{\h}^d}_{t;[x]_{\h},[y]_{\h}}}
\def\qsxynint{K^{{\text{\tiny{\sc{RW}}}}_{\h}^d}_{s;[x]_{\h},[y]_{\h}}}

\def\qrox{K^{{\text{\tiny{\sc{RW}}}^d_{\delta}}}_{r_{1};x}}
\def\qrtx{K^{{\text{\tiny{\sc{RW}}}^d_{\delta}}}_{r_{2};x}}

\def\qroprtz{K^{{\text{\tiny{\sc{RW}}}^d_{\delta}}}_{r_{1}+r_{2};0}}

\def\qrox{K^{{\text{\tiny{\sc{RW}}}^d_{\delta}}}_{r_{1};x}}
\def\qrtx{K^{{\text{\tiny{\sc{RW}}}^d_{\delta}}}_{r_{2};x}}

\def\utn{\tilde{U}_n}
\def\utnlf{\tilde{U}_{n,l}}
\def\utnk{\tilde{U}_{n_k}}
\def\utsyn{\tilde{U}_n^{y}(s)}

\def\uttxn{\tilde{U}_n^{x}(t)}
\def\uttxnk{\tilde{U}_{n_k}^{x}(t)}
\def\uttxnlf{\tilde{U}_{n,l}^{x}(t)}
\def\uttauxnlf{\tilde{U}_{n,l}^{x}(\tau)}
\def\xtaunlf{X_{n,l}^{\tau}}

\def\xttauoxnlf{X_{n,l}^{\tau_1,x_{1}}(t)}
\def\xttauoynlf{X_{n,l}^{\tau_1,y}(t)}

\def\xttautxnlf{X_{n,l}^{\tau_2,x_{2}}(t)}
\def\xttautynlf{X_{n,l}^{\tau_2,y}(t)}

\def\xttauxnlf{X_{n,l}^{\tau,x}(t)}
\def\xttauxonlf{X_{n,l}^{\tau,x_1}(t)}

\def\xttauxinlf{X_{n,l}^{\tau_{i},x_i}(t)}

\def\mtauxo{M^{\tau_1,x_{1}}(\cdot)}
\def\mtauxt{M^{\tau_2,x_{2}}(\cdot)}
\def\mtautx{M^{\tau,x}(t)}
\def\mtautxi{M^{\tau_i,x_{i}}(t)}
\def\mtautxj{M^{\tau_j,x_{j}}(t)}

\def\mtauuxj{M^{\tau_j,x_{j}}(u)}

\def\ntautxot{N^{\tau_{1,2},x_{1,2}}(t)}
\def\xztauxnlf{{X_{n,l}^{\tau,x}(0)}}
\def\xtautauxnlf{X_{n,l}^{\tau,x}(\tau)}
\def\xttxnlf{X_{n,l}^{t,x}(t)}
\def\utsynlf{\tilde{U}_{n,l}^{y}(s)}
\def\utrynlf{\tilde{U}_{n,l}^{y}(r)}
\def\xstauynlf{X_{n,l}^{\tau,y}(s)}
\def\xttauxrnlf{X_{n,l}^{\tau,x_r}(t)}
\def\xstauiynlf{X_{n,l}^{\tau_i,y}(s)}

\def\xstauxinlf{X_{n,l}^{\tau_i,x_{i}}(s)}
\def\uttx{\tilde{U}^{x}(t)}
\def\uttxD{\tilde{U}^x_D(t)}
\def\uttxR{\tilde{U}^x_R(t)}
\def\uttxnD{\tilde{U}^x_{n,D}(t)}
\def\uttxinD{\tilde{U}^{x_i}_{n,D}(t)}

\def\utsxjnD{\tilde{U}^{x_j}_{n,D}(s)}

\def\utuxinD{\tilde{U}^{x_i}_{n,D}(u)}
\def\uttxnR{\tilde{U}^x_{n,R}(t)}
\def\uttynR{\tilde{U}^y_{n,R}(t)}
\def\utrxnR{\tilde{U}^x_{n,R}(r)}
\def\uttxnlfR{\tilde{U}^x_{n,l,R}(t)}
\def\utrxnlfR{\tilde{U}^x_{n,l,R}(r)}
\def\uttynlfR{\tilde{U}^y_{n,l,R}(t)}
\def\uttxD{\tilde{U}^x_{D}(t)}
\def\uttxR{\tilde{U}^x_{R}(t)}
\def\utsxn{\tilde{U}_n^{x}(s)}
\def\utzeroxn{\tilde{U}_n^{x}(0)}
\def\uotn{{\tilde U}_n^{(1)}}
\def\utttn{{\tilde U}_n^{(2)}}
\def\yk{Y_k}
\def\ytxk{Y_k(t,x)}
\def\ytyk{Y_k(t,y)}
\def\yrxk{Y_k(r,x)}
\def\utx{U(t,x)}
\def\uty{U(t,y)}
\def\urx {U(r,x)}

\def\usy{U(s,y)}
\def\ytx{Y(t,x)}
\def\yty{Y(t,y)}
\def\yrx{Y(r,x)}
\def\uity{u_i^y(t)}
\def\ujty{u_j^y(t)}

\def\unx{u_0(x)}
\def\unxo{u_0(x_1)}
\def\unxt{u_0(x_2)}
\def\unxi{u_0(x_i)}

\def\uny{u_0(y)}
\def\un{u_0}
\def\zivtx{U^{(0)}(t,x)}
\def\zivsy{U^{(0)}(s,y)}
\def\nivtx{U^{(n)}(t,x)}
\def\nivsy{U^{(n)}(s,y)}
\def\npoivtx{U^{(n+1)}(t,x)}

\def\nmoivsy{U^{(n-1)}(s,y)}
\def\dnptx{D_{n,p}(t,x)}

\def\dnmopsy{D_{n-1,p}(s,y)}
\def\mnpt{D^*_{n,p}(t)}

\def\mnmops{D^*_{n-1,p}(s)}
\def\sW{\mathscr W}

\def\sTns{{\mathscr T}_{n,s}}
\def\sAn{{\mathscr A}_{n}}
\def\sAns{{\mathscr A}_{n}^{2}}
\def\wryn{W_n^{y}(r)}
\def\wtxn{W_n^{x}(t)}
\def\wtyn{W_n^{y}(t)}
\def\wtxm{W_m^{x}(t)}
\def\wsxn{W_n^{x}(s)}
\def\wsyn{W_n^{y}(s)}
\def\DxBt{{\mathbb D}^x_B(t)}
\def\DxkBet{{\mathbb D}^{x,k}_{B,e}(t)}
\def\DxiBet{{\mathbb D}^{x,\infty}_{B,e}(t)}
\def\DxBsnkchnt{{\mathbb D}^{x}_{B,SC}(t)}
\def\P{{\mathbb P}}
\def\Pt{\tilde{\mathbb P}}
\def\EP{{\mathbb E}_{\P}}
\def\EPt{{\mathbb E}_{\Pt}}
\def\E{{\mathbb E}}
\def\N{{\mathbb N}}
\def\Ns{{\mathbb N}^{*}}
\def\Rd{{\mathbb R}^{d}}
\def\Q{\mathbb Q}
\def\TQ{{\T}_{\Q}}
\def\Zd{{\mathbb Z}^{d}}
\def\R{\mathbb R}
\def\Rs{{\mathbb R}^2}
\def\B{{\mathbb B}}
\def\S{\mathbb S}
\def\T{\mathbb T}
\def\X{\mathbb X}
\def\Rp{{\R}_+}

\def\Xnd{{\mathbb X}_n^{d}}
\def\Xmd{{\mathbb X}_m^{d}}
\def\Xd{{\mathbb X}^{d}}
\def\Xnds{{\mathbb X}_n^{2d}}
\def\Xnldf{{\mathbb X}_{n,l}^{d}}
\def\Xnlxdf{{\mathbb X}_{n,l,x}^{d}}
\def\sF{{\mathscr F}}

\def\sFt{{\mathscr F}_t}
\def\sFs{{\mathscr F}_s}
\def\sFr{{\mathscr F}_r}

\def\OFFtP{(\Omega,\sF,\{\sFt\},\P)}
\def\OFFtPt{(\tilde{\Omega},\tilde{\sF},\{\tilde{\sFt}\},\tilde{\P})}
\def\OFPt{(\tilde{\Omega},\tilde{\sF},\tilde{\P})}

\def\OtauotFtauotPtauott{(\tilde{\Omega}_{\tau_{1,2}},\tilde{\sF}_{\tau_{1,2}},\tilde{\P}_{\tau_{1,2}})}
\def\Skspace{(\Omega^S,\sF^S,\{\sFt^S\},\P^S)}
\def\locm{{\mathscr{M}}_2^{c,loc}}
\def\C{\mathrm C}
\def\H{\mathrm H}
\def\eqdef{:=}
\def\eqd{\overset{{\mathscr L}}{=}}

\def\sm{\setminus}
\def\D{\Delta}
\def\lap{\Delta}
\def\Ds{\Delta^2}
\def\df#1#2{\ds{\frac{#1}{#2}}}
\def\tf#1#2{\ts{\frac{#1}{#2}}}
\def\lbl#1{\label{#1}}
\def\intrd{\int_{\Rd}}
\def\intrdzt{\int_{\Rd}\int_0^t}
\def\intzt{\int_0^t}

\def\pa{\partial}
\def\pas{\partial^2}
\def\lqv{\left<}
\def\rqv{\right>}

\def\lab{\left|}
\def\rab{\right|}
\def\lpa{\left(}
\def\rpa{\right)}

\def\lbk{\left[}
\def\rbk{\right]}
\def\lbr{\left\{}
\def\rbr{\right\}}
\def\bdf{\begin{defn}}
\def\edf{\end{defn}}
\def\bcr{\begin{cor}}
\def\ecr{\end{cor}}
\def\bnt{\begin{notn}}
\def\ent{\end{notn}}
\def\brm{\begin{rem}}
\def\erm{\end{rem}}
\def\blm{\begin{lem}}
\def\elm{\end{lem}}
\def\bpf{\begin{proof}}
\def\epf{\end{proof}}
\def\bpfs{\begin{pf*}}
\def\epfs{\end{pf*}}
\def\eqrf{\eqref}
\def\beq{\begin{equation}}
\def\beqs{\begin{equation*}}
\def\eeq{\end{equation}}
\def\eeqs{\end{equation*}}
\def\bsp{\begin{split}}
\def\esp{\end{split}}
\def\bc{\begin{cases}}
\def\ec{\end{cases}}
\def\bt{\begin{tabular}}
\def\et{\end{tabular}}

\def\bthm{\begin{thm}}
\def\ethm{\end{thm}}
\def\bpr{\begin{prop}}
\def\epr{\end{prop}}
\def\bfr{\begin{framed}}
\def\efr{\end{framed}}
\def\bsh{\begin{shaded}}
\def\esh{\end{shaded}}
\def\bcm{\iffalse}
\def\babs{\begin{abstract}}
\def\eabs{\end{abstract}}
\def\bit{\begin{itemize}}
\def\eit{\end{itemize}}
\def\ben{\begin{enumerate}}

\def\een{\end{enumerate}}
\def\babs{\begin{abstract}}
\def\eabs{\end{abstract}}
\def\Dn{\Delta_n}
\def\Dns{\Delta_n^2}

\def\btsie{e_{_{\mbox{\tiny BTBM}}}^{\mbox{\tiny SIE}}(a,u_0)}
\def\btsieb{e_{_{\mbox{\tiny BTBM}}}^{\mbox{\tiny SIE}}(a,b,u_0)}

\def\btrwsien{e_{_{\mbox{\tiny BTRW}}}^{\mbox{\tiny SIE}}(a,u_0,n)}
\def\btrwsienl{e_{_{\mbox{\tiny BTRW}}}^{\mbox{\tiny t-SIE}}(a,u_0,n,l)}
\def\btrwsienlaux{e_{_{\mbox{\tiny BTRW}}}^{\mbox{\tiny aux-SIE}}(a,u_0,n,l,\tau)}

\def\btsdden{e_{_{\mbox{\tiny BTRW}}}^{\mbox{\tiny SDDE}}(a,u_0,n)}
\def\pbtsdden{e_{_{\mbox{\tiny BTRW}}}^{\mbox{\tiny PSDDE}}(a,u_0,n)}
\def\btspde{e_{_{\mbox{\tiny BTP}}}^{\mbox{\tiny SPDE}}(a,u_0)}

\def\sAKtauot{{\mathscr{A}}_{\Upsilon}^{\tau_{1,2}}}
\def\sAKtauotim{{\mathscr{A}}_{\Upsilon,i_m}^{\tau_{1,2}}}
\def\uttxin{\tilde{U}^\xi_n (t)}

\def\utsyn{{\tilde U}_n^{y}(s)}

\def\utstxdn{{\tilde U}_n^{x,\cdot}(s,t)}

\def\utttxdn{{\tilde U}_n^{x,\cdot}(t,t)}
\def\utstxyn{{\tilde U}_n^{x,y}(s,t)}

\def\utrtxyn{{\tilde U}_n^{x,y}(r,t)}

\def\utryn{{\tilde U}_n^{y}(r)}

\def\utttxxn{\tilde{U}_n^{x,x}(t,t)}
\def\utztxyn{\tilde{U}_n^{x,y}(0,t)}
\def\ig{\iffalse}

\title[Ultra Regular BTBM SIE{\scriptsize s} on $\Rp\times\Rd$]
{Brownian-time Brownian motion SIE{\scriptsize s} on $\Rp\times\Rd$: Ultra Regular direct and lattice-limits solutions 
 and fourth order SPDE{\scriptsize s} links}
\author[Hassan Allouba]{}
\email{allouba@math.kent.edu}


\begin{document}
\maketitle

\vspace{-1mm}
\centerline{\scshape Hassan Allouba}
\medskip
{\footnotesize
 \centerline{Department of Mathematical Sciences, Kent State University, Kent, Ohio 44242, USA}
} 
\begin{abstract}
We delve deeper into the compelling regularizing effect of the Brownian-time Brownian motion density, $\KBtxy$, on the space-time-white-noise-driven stochastic integral equation we call BTBM SIE: \beq\lbl{z} U(t,x)=\intrd\KBtxy \uny dy+ \intrdzt\KBtsxy a(\usy)\sW(ds\times dy),\eeq which we recently introduced in \cite{Abtpspde}.  In sharp contrast to traditional second order heat-operator-based SPDEs---whose real-valued mild solutions are confined to $d=1$---we prove the existence of solutions to \eqref{z} in $d=1,2,3$  with dimension-dependent and striking H\"older regularity,  under both less than Lipschitz and Lipschitz conditions on $a$.   In space, we show an unprecedented nearly local Lipschitz regularity for $d=1,2$---roughly, $U$ is spatially twice as regular as the Brownian sheet in these dimensions---and we prove nearly local H\"older $1/2$ regularity in $d=3$.  In time, our solutions are locally $\gamma$-H\"older continuous with exponent $\gamma\in\lpa0,\tf{4-d}{8}\rpa,\ 1\le d\le3$.  To investigate \eqref{z} under less than Lipschitz conditions on $a$, we (a) introduce the Brownian-time random walk---a special case of lattice processes we call Brownian-time chains---and we use it to formulate the spatial lattice version of \eqref{z}; and (b) develop a delicate variant of Stroock-Varadhan martingale approach, the K-martingale approach, tailor-made for a wide variety of kernel SIEs including \eqref{z} and the mild forms of many SPDEs of different orders on the lattice. Solutions to \eqref{z} are defined as limits of their lattice version.
Along the way, we prove interesting aspects of Brownian-time random walk, including a fourth order differential-difference equation connection.  We also prove existence, pathwise uniqueness, and the same H\"older regularity for \eqref{z}, without discretization, in the Lipschitz case.  The SIE \eqref{z} is intimately connected to intriguing fourth order SPDEs in two ways. First, we show that \eqref{z} is connected to the diagonals of a new unconventional fourth order SPDE we call parametrized BTBM SPDE. Second, replacing $\KBtxy$ by the intimately connected kernel of our recently-introduced imaginary-Brownian-time-Brownian-angle process (IBTBAP), \eqref{z} becomes the mild form of a Kuramoto-Sivashinsky (KS) SPDE with linearized PDE part.   Ideas and tools developed here are adapted in separate papers to give an entirely new approach, via our explicit IBTBAP representation, to many linear and nonlinear KS-type SPDEs in multi-spatial dimensions.
\end{abstract}

\tableofcontents
\section{Introduction and statement of results}\lbl{intro}
\subsection{Motivation and the first main theorem}
A fascinating aspect of the Brownian-time processes\footnote{A
BTP, in its simplest form,  is a process $X^x\lpa\lab B_t\rab\rpa$ in which $X^x$ is a Markov process starting at $x\in\Rd$ and
$B_t$ is an independent one dimensional BM starting at $0$.
A Brownian-time Brownian motion (BTBM) is a BTP in which
$X^x$ is also a Brownian motion.  BTPs include many new and quite interesting processes (see \cite{Abtp1,Abtp2,DeB,Nanesd}), which we are currently investigating in several directions (e.g., \cite{AN,AD,AL,AX}).
With the exception of the Markov snake of Le Gall (\cite{LeG}), BTPs fall outside the classical theory
of Markov, Gaussian, or semimartingale processes.  See also related multiparameter processes and their PDEs in \cite{AN,Asheet}.} (BTPs) we introduced in \cite{Abtp1,Abtp2} is the rich interplay between them and fourth order PDEs.  On
one hand, BTPs solve non-Markovian (memory preserving) fourth order PDEs involving a positive bi-Laplacian that is coupled with a time-scaled
positive Laplacian so as to produce smooth solutions, for all times and all spatial dimensions.  The canonical such deterministic PDE is
\beq\lbl{btppdedet}
\bc
\df{\pa u}{\pa t}= \df{\Delta\un }{\sqrt{8\pi t}}+\df18\Delta^2 u;&(t,x)\in(0,\infty)\times\Rd, \cr
u(0,x) = \unx;&x\in\Rd.
\ec
\eeq
From a PDE perspective, BTPs give rise to regular ($\mathrm{C}^{1,4}$) solutions to new different fourth
order PDEs---like \eqref{btppdedet}---that involve the positive bi-Laplacian, for all times and in all spatial dimensions.  They are also connected to equivalent time-fractional PDEs as shown in (\cite{AN,Abtp1,MNV09,Nanesd}).  The intrigue comes
not only from directly connecting these BTPs to PDEs despite the lack of classical properties for the underlying processes (non-semimartingales, non-Marokvian, and non-Gaussian), but also from the fact that typical positive bi-Laplacian PDEs are not well behaved.  Nevertheless, BTPs lead to equations in which
the positive bi-Laplacian is coupled in a very specific way---dictated by the BTP probability density function---with
a time-scaled Laplacian acting on the smooth initial data whose smoothing effect gets arbitrarily large as time $t\searrow0$ and fades away as $t\nearrow\infty$ at the rate of $1/\sqrt{8\pi t}$; and the BTP solutions to these BTP PDEs are eternally highly regular.

On the other hand, tweaking the BTPs a little by running
the Brownian-time on the imaginary axis and adding a Brownian ``angle'' we obtain the imaginary-Brownian-time-Brownian-angle process (IBTBAP\footnote{Since we introduce several new concepts, stochastic processes, and equations in this paper, as well as build on our recent work, we have included a glossary of frequently used acronyms and abbreviations at the end of the paper for the convenience of the reader (see Appendix \ref{glossary}).}
) \cite{Aks}.  The IBTBAP in turn gives a probabilistically inspired representation for the solution of a linearized version of the prominent fourth order Kuramoto-Sivashinsky (KS) PDE of modern applied mathematics
\begin{equation} \label{lksdet}
 \begin{cases} \displaystyle\frac{\partial u}{\partial t}=
-\frac18\lap^2u-\frac12\lap u-\frac12u, & (t,x)\in(0,\infty )\times\Rd;
\cr u(0,x)=\unx, & x\in\Rd,
\end{cases}
\end{equation}
as we showed in \cite{Aks}.  Because of the intimate relation between BTPs and the IBTBAP, their kernels have a similar regularizing effect on solutions to their respective PDEs.  This is  despite the fact that BTP PDEs involve the positive bi-Laplacian while the KS PDE contains the more traditional well behaved negative bi-Laplacian.
It is at least as intriguing to study BTP-connected stochastic equations driven by additive and multiplicative space-time white noise.
In this paper, we continue our study of the quite-notable regularizing effect the Brownian-time Brownian motion (BTBM) kernel has on the space-time white noise driven BTBM  stochastic integral equation (BTBM SIE), first introduced in \cite{Abtpspde}:
\beq\lbl{btpsol}
\bsp
U(t,x)=\intrd\KBtxy \uny dy+ \intrdzt\KBtsxy a(\usy)\sW(ds\times dy)
\end{split}
\eeq
where $\sW$ is the white noise on $\Rp\times\Rd$ and $\KBtxy$ is the density of a Brownian-time Brownian motion given by:
\beq\lbl{btpden}\KBtxy=2\int_0^\infty \psxy\ptzs ds
\eeq
with  $\psxy=\tf{e^{-|x-y|^2/2s}}{(2\pi s)^{{d}/{2}}}$ and $\ptzs=\tf{e^{-s^2/2t}}{\sqrt{2\pi t}}$.  We call \eqref{btpsol} BTBM SIE  (or BTP SIE) since it is expressed in terms of the density (or kernel) of some BTP, a Brownian-time Brownian motion in this case of \eqref{btpsol}.  We denote our BTBM SIE \eqrf{btpsol} by $\btsie$.

In \cite{Abtpspde} we considered the additive noise case $a\equiv1$ for $\btsie$, and we
proved the existence of a pathwise unique continuous BTBM SIE solution\footnote{Unlike the deterministic case $a\equiv0$, the BTBM SIE \eqref{btpsol} is not the kernel integral form of the stochastic PDE version of \eqref{btppdedet}, viz.~the BTBM SPDE
\beq\lbl{btpspde}
\bc
\df{\pa U}{\pa t}= \df{\Delta\un }{\sqrt{8\pi t}}+\df18\Delta^2 U+
a(U) \df{\pa^{d+1} W}{ \pa t\pa x},&(t,x)\in(0,\infty)\times\Rd; \cr
U(0,x) = \unx,&x\in\Rd.
\ec
\eeq
I.e., the BTBM SIE \eqref{btpsol} and the BTBM SPDE \eqref{btpspde} are two different stochastic versions of \eqref{btppdedet}.  It is $\btsie$ that has ultra regular solutions in $d=1,2,3$---capturing the smoothing effect of the BTBM density in the stochastic setting---and it is the equation intimately related to the KS and other important fourth order SPDEs of modern applied mathematics.  We therefore focus the bulk of our investigation and our main results in this article (\thmref{lip} and \thmref{mainthm2})  on \eqref{btpsol}.  In \secref{spdeslinks} we discuss further the links between $\btsie$ and fourth order SPDEs of KS type as well as the new parametrized BTBM SPDE, relative to which \eqref{btpspde} may be thought of as a rougher version.  The proof of Theorem 3.1 in \cite{Abtpspde} and the associated existence, uniqueness, and dimension-dependent $L^{p}$-regularity are all correct for the BTBM SIE.  Only the name should change from BTP SPDE solutions to BTBM SIE solutions in  Theorem 3.1 and Definition 3.1  of \cite{Abtpspde}.  We analyze BTBM SPDEs and other related stochastic fractional PDEs in upcoming articles.} $U(t,x)$
for $x\in\Rd$ and $1\le d\le 3$, such that
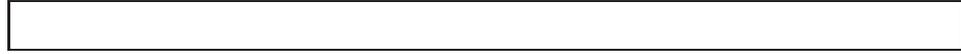
$$\sup_{x\in\Rd}\EP|U(t,x)|^{2p}\le C\lbk1+t^{\tf{(4-d)p}{4}}\rbk;\quad t>0,\ 1\le d\le 3,\  p\ge1.$$
In the second main result of this paper, we prove existence and finer dimension-dependent H\"older regularity results for \eqref{btpsol} under the following less than Lipschitz conditions\footnote{Here, $\mathrm{C}_{b}^{2,\gamma}(\Rd;\R)$ denotes the subspace of the standard space $\mathrm{C}_{b}^2(\Rd;\R)$ (see Appendix \ref{glossary}) in which all second derivatives are H\"older continuous, with some H\"older exponent $0<\gamma\le1$.
Also, the boundedness conditions on $\un$ and its derivatives may easily be relaxed as in \cite{AN}.}
on the Borel-measurable diffusion coefficient $a$:
\beq\lbl{cnd}\tag{\mbox{NLip}}
\bc
(a)&\hspace{-2.5mm} a(u) \mbox { is continuous in $u$;}\quad u\in\R,\\
(b)&\hspace{-2.5mm} a^2(u)\le C (1+u^2);\quad u\in\R,\\
(c)&\hspace{-2.5mm} u_0\in \mathrm{C}_{b}^{2,\gamma}(\Rd;\R)\mbox{ and nonrandom for some } 0<\gamma\le1,\ \forall\ 1\le d\le3.
\ec
\eeq
Our first main result gives stronger existence, as well as uniqueness, and the same H\"older regularity  for $\btsie$ under an added Lipschitz condition on $a$:
\beq\lbl{lcnd}
\bc
\hspace{-2.5mm}&(a) \lab a(u)-a(v)\rab\le C\lab u-v\rab \mbox \quad u,v\in\R;\\
\hspace{-2.5mm}&(b) \mbox{ and } (c)\mbox{ same as in \eqref{cnd}}.
\ec
\tag{\mbox{Lip}}
\eeq
More precisely, we denote by\footnote{Throughout the paper, $\T=[0,T]$ for some fixed but arbitrary $T>0$.} $\H^{\gamma_t^{-},\gamma_s^{-}}(\T\times\Rd;\R)$ the space of real-valued locally H\"older functions on $\T\times\Rd$
whose time and space H\"older exponents are in $(0,\gamma_t)$ and $(0,\gamma_s)$, respectively.  Our first main result is now stated directly for $\btsie$, without any lattice approximations.
\bfr
\bthm[Direct solutions to $\btsie$ for dimensions $1\le d\le3$]\lbl{lip}
Assume that \eqref{lcnd} holds.  Then there exists a pathwise-unique strong solution $(U,\sW)$ to $\btsie$
on $\Rp\times\Rd$, for $1\le d\le3$, which is $L^p(\Omega)$-bounded
on $\T\times\Rd$ for all $p\ge2$.  Furthermore, $U\in\H^{{\tf{4-d}{8}}^{-},{\lpa\tf{4-d}{2}\wedge 1\rpa}^{-}}(\T\times\Rd;\R)$ for every $1\le d\le3$.
\ethm
\efr
\thmref{lip} states that the stochastic kernel integral equation \eqref{btpsol} has ultra regular strong\footnote{Here strong is in the stochastic sense of the noise $\sW$ and its probability space $\OFFtP$ being fixed a priori.  Throughout this article, whenever needed, we will assume that our filtrations satisfy the usual conditions without explicitly stating so.} solutions on $\Rp\times\Rd$, namely $U\in\H^{{\tf{3}{8}}^{-},1^{-}}(\T\times\R;\R)$, $U\in\H^{{\tf{1}{4}}^{-},1^{-}}(\T\times\R^{2};\R)$, and $U\in\H^{{\tf{1}{8}}^{-},{\tf12}^{-}}(\T\times\R^{3};\R)$.   I.e., in space, we show a rather remarkable---and initially-surprising---nearly local Lipschitz regularity for $d=1,2$; and we prove nearly local H\"older $1/2$ regularity in $d=3$.  This is remarkable because the BTBM kernel is able, in $d=1,2$, to spatially regularize such solutions beyond the traditional H\"older-$1/2^{-}$ spatial regularity of the underlying Brownian sheet corresponding to the driving space-time white noise\footnote{It is important to note here that the common ``folklore wisdom'' of solutions of space-time-white-noise driven equations not being smoother than the associated Brownian sheet---in either space or time---originated from  the predominant case of SPDEs, in which either the underlying kernel is that of a Brownian motion or the spatial operator is a Laplacian.  The kernel $\KBtxy$, however, is much more regularizing to the space-time-white-noise driven $\btsie$ than the density of BM is to its corresponding equation. This becomes evidently clear in \lemref{3rdQinequality}, \lemref{4thQinequality}, and \lemref{2ndQinequality} (compare to the more traditional BM and random walk case in \cite{Asdde2}).}.  This degree of smoothness is unprecedented for space-time white noise driven kernel equations or their corresponding SPDEs; and the BTBM SIE is thus the first such example.  In time, our solutions are locally $\gamma$-H\"older continuous with dimension-dependent exponent $\gamma\in\lpa0,\tf{4-d}{8}\rpa$ for $1\le d\le3$.  This is in sharp contrast to traditional second order reaction-diffusion (RD) and other heat-operator-based SPDEs driven by space-time white noise, whose fundamental kernel is the Brownian motion density and whose real-valued mild solutions are confined to the case $d=1$.  This sharp contrast in regularity\footnote{We observe in passing that---roughly speaking---the paths of $\btsie$ in $d=1$ are effectively $3/2$ times as smooth as the RD SPDE paths in $d=1$, in $d=2$ the BTBM SIE is as smooth as an RD SPDE in $d=1$, and in $d=3$ our BTBM SIE is half as smooth as an RD SPDE in $d=1$.  Also, for $d=2,3$, the spatial regularity is roughly four times the temporal one, and in $d=1$ the spatial regularity is maximized at a near Lipschitz vs near H\"older $3/8$ in time (see also Table \ref{tabl}).} is summarized in \tablref{tabl}.  In this regard, the dichotomy between the rougher paths of BTBMs as compared to standard Brownian motions on the one hand (quartic vs.~quadratic variations) and the stronger regularizing properties of the BTBM density vs.~the BM one on the other hand is certainly another interesting point to make.  Here, it is indeed interesting to note that  random field  solutions for the BTBM SIE exist only for $d=1,2,3$; which are the same dimensions for which there exists a local time for the iterated Brownian motion---a special case of the BTBM class (see \cite{Abtp1}) with the same density \eqref{btpden}---as proved by Xiao in Theorem 1 of \cite{X}.  To the best of our knowledge, some of the earlier  computations linking local times and white noise driven SPDEs, through the underlying kernels, appeared in Lemma A.2.3 in the paper by Reimers \cite{Reim} for heat SPDEs.    Later, Foondun et~al. \cite{FKN} showed the equivalence of the local time existence and the existence of random fields solutions for a class of white noise driven SPDEs in which the spatial operator is the generator of a L\'evy process (and a weaker result if it is the generator of a Markov process)\footnote{Such SPDEs in \cite{FKN} do not possess random field solution in $\Rd$ for $d>1$ in the L\'evy case and do not have random field solutions in $\Rd$ for $d\ge2$ in all cases (see \cite{FKN} for the details on the existence of ``random field'' solutions for $d=2-\epsilon$ on  fractal subsets of $\Rd$).}.   We believe this equivalence extends to our BTBM SIE---which is outside the setting in \cite{FKN}---as the above discussion suggests.  We plan to address, in a future article, this issue in a setting that includes the BTBM SIE case and covers dimensions $d>1$.

Even under the Lipschitz condition on $a$, the proof of the fine dimension-dependent and uncoventional spatio-temporal H\"older regularity in \thmref{lip} requires delicate estimates on the spatial and temporal differences of the Brownian-time Brownian motion density $\KBtxy$ (\lemref{4thQinequality} and \lemref{3rdQinequality} below).  Some of the computations there are probabilistically-flavored with connections to yet other new processes that we introduce in the course of these computations  (e.g., $2$-Brownian-times random walk and $2$-Brownian-times Brownian motion).

Second, due to the intimate BTPs-IBTBAP connection; and just as we showed in the deterministic PDEs case \cite{Aks,Abtp2,Abtp1}---where methods from BTPs PDEs were adapted to prove results for the linearized KS PDE \eqref{lksdet} {\it in all dimensions}---our study here is very helpful in our related investigation of many prominent, as well as new,  KS-type SPDEs which we treat in separate papers \cite{Aksspde,AL,AD} and followup articles.     The methods presented here and in \cite{Aks} are adapted and generalized in \cite{Aksspde,AL,AD} to give an entirely new approach, in terms of the IBTBAP kernel and related probabilistically-motivated concepts, to the SPDEs version of these famous fourth order applied mathematics PDEs in $d=1,2,3$.  Here, it is noteworthy that even the existence/nonexistence of the KS semigroup is not known for $d>1$ using standard analytical methods\footnote{We also note here that our probabilistic density (or kernel) approach, allows us to obtain sharp dimension-dependent regularity estimates in both space and time, simultaneously.  It is instructive, for example, to compare our spatio-temporal regularity to ones obtained for other bi-Laplacian SPDEs by analytical semigroup methods (e.g., \cite{DaDe}).  Our effective regularity agrees with \cite{DaDe} in $d=3$---H\"older $\lpa\tf18\rpa^{-1}$---but is much sharper in $d=1,2$. (H\"older $\lpa\tf38\rpa^{-1}$ and $\lpa\tf14\rpa^{-1}$ vs.~H\"older $\lpa\tf18\rpa^{-1}$ in $d=1,2$, respectively).  We'll address this further in planned future articles.}.  The regularity of the SPDEs in \cite{Aksspde,AL,AD} is very similar to this paper's result.  We therefore regard $\btsie$ as a cousin of and a companion model for such important SPDEs.

Before we can precisely state our results under the conditions in \eqref{cnd} (\thmref{mainthm2}), we need to introduce the lattice version of $\btsie$ as well as introduce the new associated process we call Brownian-time random walk and define the lattice limit solutions involved in the statement of \thmref{mainthm2}.  To prepare for the proof of \thmref{mainthm2} we also need to introduce new machinery, the K-martingale approach.

We now detail the structure of the remainder of this article.  The rest of \secref{intro} provides the setting then states the second main result---\thmref{mainthm2}---under the conditions \eqref{cnd}, \secref{pfs} contains the proofs of several lemmas and results leading to the proof of the main results for $\btsie$ and it contains a proof of fourth order SPDEs-BTBM SIEs connection on the lattice, \secref{conc} contains additional concluding remarks.  More specifically, from \secref{btc} to \secref{limitsol}; we give all the ingredients necessary to state our second main existence and regularity result, given in \secref{mres}. In \secref{btc} we introduce a new class of discrete-valued processes that we call Brownian-time chains (BTCs). BTCs are the discretized versions of our BTPs in \cite{Abtp1,Abtp2}.  Of particular interest here is the special case of Brownian-time random walk (BTRW), which we define and link to the lattice version of \eqref{btppdedet} (\lemref{fodde}).  In \secref{limitsol}, we use the density of the BTRW to give a spatially-discretized formulation \eqref{btrwsie} of $\btsie$, which we call BTRW SIE, and we define two notions of solutions to the lattice model: direct solutions and limit solutions (from a finite truncation of the lattice to the whole lattice).  These solutions (both direct and limit) are then used to define two types of lattice BTRW SIEs limit solutions to $\btsie$ (direct limit solutions and double limit solutions), as the size of the lattice mesh shrinks to zero\footnote{We deal with three major types of solutions for $\btsie$ in this article: direct (in which no lattice approximation is used), direct limits solutions (in which solutions are defined as limits of their lattice approximations), and double limits solutions (in which the solutions are defined as limits of their lattice approximations that are, in turns, themselves obtained as limits of their finite truncations).}.  We introduce our K-martingale\footnote{Here, K is for kernel.} approach to kernel SIEs as $\btsie$ in \secref{Kmartsec}.  It is an essential ingredient in the proof of existence for $\btsie$ under the less-than-Lipschitz conditions on $a$ \eqref{cnd}. It is a delicate variant of the well known, and by now classic, martingale problem approach of Stroock and Varadhan to SDEs.  Our K-martingale approach starts by constructing an auxiliary problem to a truncated lattice version of \eqref{btpsol}, for which the existence of solutions implies solutions existence for the truncated lattice model.  We then formulate a martingale problem equivalent to the auxiliary problem (the K-martingale problem).  A key advantage of the K-martingale approach is that it is a unified framework in which the existence and uniqueness of many kernel stochastic integral equations, which are the mild formulation for many SPDEs, may be treated under less than Lipschitz conditions; using only variants of the kernel formulation of the underlying equation.  This includes SPDEs of different orders (including second and fourth), so long as the corresponding spatially-discretized kernel (or density) satisfies Kolmogorov-type bounds on its temporal and spatial differences.  In essence, what the K-martingale approach implies is that if the kernel in the lattice model is nice enough for the lattice model to converge as the lattice mesh shrinks to zero (under appropriate assumptions on $a$), then it is nice enough to guarantee a solution for the lattice model.   We use it here to prove the existence of solutions to \eqref{btpsol} under the conditions \eqref{cnd}, but just as with the Stroock-Varadhan method, it can handle uniqueness as well.
\secref{mres} contains the second main existence and dimension-dependent regularity result for $\btsie$, \thmref{mainthm2}, using lattice-limits solutions.  The upshot is that under both Lipschitz and non-Lipschitz conditions on $a$ (\eqref{lcnd} and \eqref{cnd}) there exist BTRW SIE limit solutions (defined in \secref{limitsol}).  These solutions are extracted as weak limits of a system of BTRW SIEs (the type of limit solution depends on the conditions: direct limit solution for the Lipschitz case and double limits solutions for the less-than-Lipschitz case\footnote{Since we already handle the Lipschitz condition on $a$ in \thmref{lip} directly, without discretization, we relegate the statement and proof sketch of the lattice limits solutions to $\btsie$ under \eqref{lcnd} to Appendix \ref{appB} (\thmref{latlimlip}).}).   This discretization method is similar in spirit to our discretization approach for the simpler second order reaction-diffusion (RD) SPDEs in \cite {Asdde1,Asdde2}.
These lattice-limits solutions have the same spatio-temporal H\"older regularity as their direct counterpart in \thmref{lip}.  We prove here that these limit solutions are temporally locally H\"older continuous with H\"older exponent $\gamma\in(0,\tf{4-d}{8})$ for spatial dimensions $d=1,2,3$.  Spatially, they have an impressive nearly local Lipschitz regularity for $d=1,2$, and nearly local H\"older $1/2$ regularity in $d=3$.    \secref{spdeslinks} gives the BTBM SIEs-SPDEs connections.  It gives a brief look into the indirect KS-type SPDEs connection---via the IBTBAP kernel---to $\btsie$.  It also provides the intriguing connection of $\btsie$ to the parametrized BTBM SPDE on the lattice.  In \secref{pfs} we prove different BTRW and BTBM densities estimates---introducing along the way the notion of $2$-Brownain-times RW and BM; we prove some BTRW SIEs estimates necessary for regularity and tightness; we prove the K-martingale result; and we prove the main existence, uniqueness, and regularity results for $\btsie$ (\thmref{lip} and \thmref{mainthm2}).  Some interesting related results are discussed and/or proved in Appendices \ref{appA}. \ref{appB}, and \ref{appC}.  Appendix \ref{appA} gives the proof of the BTRW fourth order differential-difference equation on lattices in \lemref{fodde}.  Appendix \ref{appB} contains the statement and proof of the existence, uniqueness, and regularity for a direct limit solution to $\btsie$ under Lipschitz conditions (\thmref{latlimlip}) via an iterative-type argument similar to that in the proof of \thmref{lip} together with some estimates obtained in \secref{pfs}.  Appendix \ref{appC} contains more on the BTBM SIEs-SPDEs connections on lattices.  Appendix \ref{glossary} contains frequently used acronyms and notations.

\subsection{Brownian-time random walk and chains on the lattice}\lbl{btc}
In \cite{Asdde1,Asdde2}, standard continuous-time random walks
on a sequence of refining spatial lattices $$\lbr\Xnd:=\prod^{d }_{i=1 }\lbr\ldots, -2\h,-\h,0,\h, 2\h, \ldots\rbr=\h\Zd\rbr_{n\ge1}$$
(with the step size $\h\searrow0$ as $n\to\infty$) played a crucial role---through their densities---in
obtaining our results for second order
RD SPDEs.  Here, in the fourth order Brownian-time setting, that role is played by
Brownian-time random walks on $\Xnd$:
\beq\lbl{btrwdef}
\S^x_{B,\h}(t)\eqdef S^x_{\h}\lpa\lab B_t\rab\rpa;\quad 0\le t<\infty, x\in\Xnd
\eeq
where $S_{\h}^x(t)$ is a standard $d$-dimensional continuous-time symmetric RW
starting from $x\in\Xnd$ and $B$ is an independent one-dimensional
BM starting at $0$.  The subscript $\delta_n$ in \eqref{btrwdef}
is to remind us that the lattice step size is $\delta_n$ in each of the $d$ directions.
It is then clear that the transition probability (density) $\Ktxyn$ of the BTRW $\S_{B,\h}^x(t)$ on
$\Xnd$ is given by
\beq\lbl{btrwdsty}
\Ktxyn=2\int_0^\infty\qsxyn\ptzs ds; \quad 0< t<\infty, \ x,y\in\Xnd
\eeq
with $\ptzs
=\lpa 1/\sqrt{2\pi t}\rpa\exp\lbr- s^2/2t\rbr$ and
$\qtxyn$ is the continuous-time random walk transition density
starting at $x\in\Xnd$ and going to $y\in\Xnd$ in time $t$,
in which the times between transitions are exponentially distributed
with mean $\h^{2d}$.  Throughout this article, $\Ktxn\eqdef\Ktxzn$ (with a similar convention for all transition densities).  I.e., $\qtxn$ is the fundamental solution to the deterministic
heat equation on the lattice $\Xnd:$
\beq\lbl{latheat}
\displaystyle \frac{du_{n}^x(t)}{dt}
 = \frac{1}{2}\Delta_nu_{n}^x(t);\hspace{2mm}(t,x)\in(0,\infty)\times\Xnd
\eeq
where $\mathscr{A}_n\eqdef \Delta_n/2$ is the generator of the RW $S^x_{\h}(t)$ on $\Xnd$.

By mimicking our proof of Theorem 0.1 in \cite{Abtp1} (see the proof in
Appendix \ref{appA}), we easily get the following fourth order differential-difference equation connection to BTRW:
\blm[BTRW's DDE]\lbl{fodde}
Let $u^x_n(t)=\E \lbk u_0\lpa \S^x_{B,\h}(t)\rpa\rbk$ with $u_0$ as in \eqref{cnd}.
Then $u_n$ solves the following fourth order differential-difference equation (DDE) on
$\Rp\times\Xnd$:
\beq\lbl{btrwdde}\bc
\df{d u^x_n(t)}{d t}=\df{\Delta_n u_0(x)}{\sqrt{8\pi t}}+\df{1}{8}\Delta^2_n u_n^x(t),
& (t,x)\in(0,\infty)\times\Xnd\\
u_n^x(0)=u_0(x),&x\in\Xnd
\ec
\eeq
and $\Ktxn$ solves \eqref{btrwdde} on $[0,\infty)\times\Xnd$, with
\beq
u_0(x)=\Kzxn=\qzxn=
\bc
1,&x=0\cr
0,&x\neq0.
\ec
\eeq
\elm
BTRWs are the discretized version of our BTBM in \cite{Abtp1,Abtp2}.  They belong
to a large and new class of discrete-valued processes which we now introduce.  Suppressing the $n$
in the lattices $\Xnd$, let
$B$ be a one-dimensional Brownian motion starting at $0$ and let $D^x$ be
an independent $d$-dimensional $\Xd$-valued continuous-time Markov chain starting at
$x$, both defined on a probability space $\OFFtP$.  We call the process
$\DxBt\eqdef D^x(|B_t|)$  a Brownian-time chain (BTC).  A BTRW
is a special case of BTCs in which $D^x$ is a continuous-time random walk.
Excursions-based Brownian-time chains (EBTCs) are obtained from BTCs by breaking up the
path of $|B_t|$ into excursion intervals---maximal intervals $(r, s)$ of time on which
$|B_t| > 0$---and, on each such interval, we pick an independent copy of the Markov
chain $D^x$ from a finite or an infinite collection.  BTCs and EBTCs may be regarded as
canonical constructions to some quite interesting new processes:
\begin{enumerate}
\item  Markov snake chain: when $|B_t|$ increases we pick a new chain $D^x$, we denote this process
by $\DxBsnkchnt$.
\item  $k$-EBTCs: let $D^{x,1},\ldots,D^{x,k}$ be independent copies of $D^x$ starting from point $x\in\Xd$.
On each $|B_t|$ excursion interval, use one of the copies chosen at random.   we denote such a process by
$\DxkBet$.  
When $k=1$ we obtain a BTC.
\item  $\infty$-EBTCs: we use an independent copy of $D^x$ on each $|B_t|$ excursion interval.
This is the $k\to\infty$ of (2),  It is intermediate between (1) and (2).
Here, we go forward to a new independent chain only after $|B_t|$ reaches $0$.
This process is denoted by $\DxiBet$.
\end{enumerate}
\subsection{Lattice BTRW SIEs and their limits solutions to BTBM SIEs}\lbl{limitsol}
The crucial role of the BTRW density in our approach to the BTBM SIE \eqref{btpsol}
becomes even clearer from the following definition of our approximating spatially-discretized
equations:
\begin{defn}[Lattice BTRW SIEs]\label{btrwsiedef}
By the BTRW SIEs associated with the BTBM SIE $\btsie$ we mean
the system $\lbr\btrwsien\rbr_{n=1}^\infty$ of spatially-discretized stochastic integral equations
on $\Rp\times\Xnd$ given by
\beq\lbl{btrwsie}
\bsp
\uttxn=\sum_{y\in\Xnd}\Ktxyn\uny+\sum_{y\in\Xnd}\int_0^t\Ktsxyn a(\utsyn)\df{d\wsyn}{\sdh},
\end{split}
\eeq
where the BTRW density is given by \eqref{btrwdsty}.  For each $n\in\N$, we think of $\{\wtxn; t\ge0\}$ as
a sequence of independent standard
Brownian motions indexed by the set $\Xnd$ (independence within the same lattice).  We also
assume that if $m\neq n$ and $x\in\Xmd\cap\Xnd$ then $\wtxm=\wtxn$, and if $n>m$ and
$x\in\Xnd\sm\Xmd$ then $\wtxm=0$.
\end{defn}
\bnt\lbl{notdetsto}
We will denote the deterministic and the random parts of \eqref{btrwsie} by $\uttxnD$ and $\uttxnR$
(or $\uttxD$ and $\uttxR$ when we suppress the dependence on $n$), respectively,
whenever convenient.
\ent
We define two types of solutions to BTRW SIEs: direct solutions and limit solutions.
\bdf[Direct BTRW SIE Solutions]\lbl{btrwdsol}
A direct solution to the BTRW SIE system $\lbr\btrwsien\rbr_{n=1}^\infty$ on $\Rp\times\Xnd$
with respect to the Brownian $($in $t$\/$)$ system $\lbr\wtxn;t\ge0\rbr_{(n,x)\in\N\times\Xnd}$ on
the filtered probability space $\OFFtP$ is a sequence of real-valued processes
$\lbr\tilde{U}_{n}\rbr_{n=1}^\infty$  with continuous sample paths in $t$ for each fixed $x\in\Xnd$ and $n\in\N$
such that, for every $(n,x)\in\N\times\Xnd$, $\uttxn$ is $\sFt$-adapted, and equation \eqref{btrwsie} holds
$\P$-a.s.  A solution is said to be strong if $\{\wtxn;t\ge0\}_{(n,x)\in\N\times\Xnd}$ and $\OFFtP$ are fixed a priori$;$
and with
\begin{equation}
\sFt=\sigma\left\{\sigma\left(\wsxn;0\le s\le t,x\in\Xnd,n\in\N\right)\cup
\mathscr N\right\};\quad t\in\Rp,
\label{filt2}
\end{equation}
where $\mathscr N$ is the collection of null sets
$$\left\{O:\exists\,G\in\mathscr{G},O\subseteq G\ \mbox{and}\ \P(G)=0\right\}$$
and where
$$\mathscr{G}=\sigma\left(\bigcup_{t\ge0}\sigma\left(\wsxn;0\le s\le t,x\in\Xnd,n\in\N\right)\right).$$
A solution is termed weak if we are free to
choose $\OFFtP$ and the Brownian system on it and without requiring $\sFt$ to
satisfy $\eqnref{filt2}$.  Replacing $\Rp$ with $\T:=[0,T]$---for some $T>0$ in the
above, we get the definition of  a solution to the BTRW SIE system $\lbr\btrwsien\rbr_{n=1}^\infty$ on
$\T\times\Rd$.
\edf
The next type of BTRW SIE solutions we define is the first step in our K-martingale approach  of \secref{Kmartsec}.   By first reducing $\btrwsien$ to the simpler finite dimensional noise setting, it takes full advantage of the notion of BTRW SIEs limit solutions to BTBM SIEs.
\bdf[Limit BTRW SIE Solutions]\lbl{btrwlsol}
Let $l\in\N$.  By the $l$-truncated BTRW SIE on $\Rp\times\Xnd$ we mean the BTRW SIE obtained from \eqref{btrwsie}
by restricting the sum in the stochastic term to the finite $d$-dimensional lattice $\Xnldf:=\Xnd\cap\lbr[-l,l]^d;l\in\N\rbr$
and leaving unchanged the deterministic term $\uttxnD$:
\beq\lbl{trnctdsie}
\uttxnlf=\bc
\ds\uttxnD
&+\ \ds\sum_{y\in\Xnldf}\int_0^t\kappa^{x,y}_{\h,s,t}\lpa\utsynlf\rpa{d\wsyn}; x\in\Xnldf,\\
\ds\uttxnD;&x\in\Xnd\setminus\Xnldf
\ec
\eeq
where $$\kappa^{x,y}_{\h,s,t}\lpa\utrynlf\rpa:=\df{\Ktsxyn}{\sdh} a(\utrynlf),\quad \forall r,s<t.$$  We denote \eqref{trnctdsie} by
$\btrwsienl$.  Fix $n\in\N$, a solution to the system of truncated BTRW SIEs $\lbr\btrwsienl\rbr_{l=1}^\infty$
on $\Rp\times\Xnd$ with respect to the Brownian $($in $t$\/$)$ system $\lbr\wtxn;t\ge0\rbr_{x\in\Xnd}$ on
the filtered probability space $\OFFtP$ is a sequence of real-valued processes $\lbr\tilde{U}_{n,l}\rbr_{l\in\N}$
with continuous sample paths in $t$ for each fixed $x\in\Xnd$ and $l\in\N$,
such that, for every $(l,x)\in\N\times\Xnd$, $\uttxnlf$ is $\sFt$-adapted, and equation \eqref{trnctdsie} holds
$\P$-a.s.  We call $\utn$ a limit solution to the BTRW SIE \eqref{btrwsie} if $\utn$ is a limit of the truncated
solutions $\utnlf$ (as $l\to\infty$).    If desired, we may indicate the limit type (a.s., in $L^p$, weak, \ldots, etc).
\edf
\brm\lbl{smuthdet}
In both \eqref{trnctdsie} and \eqref{btrwsie}, $\uttxnD=\E \lbk u_0\lpa \S^x_{B,\h}(t)\rpa\rbk$.  So, by \lemref{fodde},
$\uttxnD$ is differentiable in time $t$ and satisfies \eqref{btrwdde}.
Also, using linear interpolation, we can extend the definition of an already continuous-in-time process $\uttxn$
on $\Rp\times\Xnd$, so as to obtain a continuous process on $\Rp\times\Rd$, for each
$n\in\N$, which we will also denote by $\uttxn$.    Henceforth,
any such sequence $\{\utn\}$ of interpolated $\utn$'s will be called a continuous
or an interpolated solution to the system $\lbr\btrwsien\rbr_{n=1}^\infty$.  Similar comments
apply to solutions of the truncated $\btrwsienl$.
\erm
We now define solutions to $\btsie$ based entirely on their approximating $\lbr\btrwsien\rbr$, through their limit.
Since we defined direct and limit solutions to $\btrwsien$, for each fixed $n$, we get two types of BTRW SIEs limit solutions to $\btsie$:
direct BTRW SIEs limit solutions and BTRW SIE double limit solutions.  The ``double'' in the second type of solutions
reminds us that we are taking two limits, one from truncated to nontruncated fixed lattice (as $l\to\infty$)
and the other limit is taken as the lattice mesh size shrinks to zero (as $\h\searrow0$ or equivalently as $n\nearrow\infty$).

\begin{defn}[BTRW SIEs limits solutions to $\btsie$]\label{limitsolns}
We say that the random field $U$ is a BTRW SIE limit solution
to $\btsie$ on $\Rp\times\Rd$ iff there is a solution $\{\tilde{U}_{n}\}_{n\in\N}$
to the lattice SIE system $\lbr\btrwsien\rbr_{n\in\N}$ on a probability space
$\OFFtP$ and with respect to a Brownian system $\{\wtxn;t\ge0\}_{(n,x)\in\N\times\Xnd}$
such that $U$ is the limit or a modification of the limit of $\lbr\utn\rbr_{n\in\N}$ $($or a subsequence thereof).   A BTRW SIE limit solution $U$ is called a direct BTRW SIEs limit solution or a BTRW SIEs double limit solution
according as $\{\utn\}_{n\in\N}$ is a sequence of direct or limit solutions to $\lbr\btrwsien\rbr_{n\in\N}$.   The limits may be taken
in the a.s., probability, $L^p$, or weak sense\footnote{When desired, the types of the solution and the limit are explicitly stated
$($e.g., we say strong $($weak\/$)$ BTRW SIEs weak,
in probability, $L^p(\Omega)$, or a.s.~limit solution to indicate that the solution to the approximating
SIEs system is strong $($weak\/$)$
and that the limit of the SIEs is in the weak, in the probability, in the $L^p(\Omega)$, or in the
 a.s.~sense, respectively\/$)$.  Of course, we may also take limits in any other suitable sense.}.  We say that uniqueness in law holds if whenever $U^{(1)}$ and $U^{(2)}$ are BTRW SIEs limit solutions they have the same law.   We say that pathwise uniqueness holds for BTRW SIEs limit solutions if whenever $\lbr\uotn\rbr$ and $\lbr\utttn\rbr$ are lattice SIEs solutions on the same probability space
and with respect to the same Brownian system, their limits $U^{(1)}$ and $U^{(2)}$ are indistinguishable.
\end{defn}
\subsection{Second main theorem: the lattice-limits solutions case} \lbl{mres}
We can now state our second main result of the paper.  The following theorem gives our lattice-limits solutions result for $\btsie$ under the non-Lipschitz conditions \eqref{cnd} on $a$.  Our limits solutions result under Lipschitz conditions is stated in \thmref{latlimlip}\footnote{The type of limit solutions in the Lipschitz case  is direct weak-limit solutions as opposed to the double weak-limit solution in \thmref{mainthm2}.  Of course, double weak-limit means the two limits are in the weak sense.}, and its proof is outlined in Appendix \ref{appB}.
\bfr
\bthm[Lattice-limits solutions to $\btsie$ for dimensions $1\le d\le3$] \lbl{mainthm2}
 Assume the conditions \eqref{cnd} hold.  Then,
there exists a BTRW SIE double weak-limit solution to $\btsie$, $U$, such that
$U(t,x)$ is $L^p(\Omega,\P)$-bounded on $\T\times\Rd$ for every $p\ge2$ and $U\in\H^{{\tf{4-d}{8}}^{-},{\lpa\tf{4-d}{2}\wedge 1\rpa}^{-}}(\T\times\Rd;\R)$ for every $1\le d\le3$. \ethm
\efr
\brm\lbl{scnv}
There is a subtle distinction between the spatial regularity in the one and two dimensional cases.   We explore this further in \cite{AX} (see also \remref{spatlkerrem}).
Of course, we can use change of measure---as we did in our earlier work on Allen-Cahn SPDEs and other second order SPDEs (see e.g.~\cite{Acom,Acom1,Acom2} and all our change of measure references in \cite{Asdde2} for results and conditions)---to transfer existence, uniqueness, and law equivalence results between $\btsie$ and the BTBM SIE with measurable drift $\btsieb$:
\beq\lbl{btpsiedrft}
\bsp
U(t,x)=&\intrd\KBtxy \uny dy+ \intrdzt\KBtsxy b(\usy)ds dy
\\&+ \intrdzt\KBtsxy a(\usy)\sW(ds\times dy),
\end{split}
\eeq
under the same conditions on the drift/diffusion ratio.  If it is desired to investigate $\btsieb$ on a bounded domain in $\Rd$ with a regular boundary, we simply replace the BTBM density $\KBtxy$ in \eqref{btpsiedrft} with its boundary-reflected or boundary-absorbed version (the BTBM density in which the outside $d$-dimensional BM is either reflected or absorbed at the boundary).
\erm
Even the proof of existence in \thmref{mainthm2} under the conditions \eqref{cnd} is not straightforward---even after obtaining the new non-trivial spatio-temporal regularity estimates (in  \lemref{4thQinequality} and  \lemref{3rdQinequality}) on the unconventional kernel $\KBtxy$.   This is because standard techniques, like the classical martingale problem approach, do not apply directly to kernel equations like the BTBM SIE $\btsie$ or its discretized version (BTRW SIE $\btrwsien$ under \eqrf{cnd}.  This leads us to devise our aforementioned K-martingale approach, which we now introduce.
\subsection{The K-martingale approach}\lbl{Kmartsec}
We now describe our K-martingale approach, which is tailor-made for kernel SIEs like $\btsie$ and other mild formulations for many SPDEs on the lattice.  The first step is to truncate to a finite lattice model as in \eqref{trnctdsie}.  Of course, even after we truncate the lattice, a remaining hurdle to applying
a martingale problem approach is that the finite sum of stochastic integrals in \eqref{trnctdsie} is not a local martingale.
So, we introduce a key ingredient in this K-martingale method: the auxiliary problem associated with the truncated BTRW SIE in \eqref{trnctdsie}, which
we now give.  Fix $(l,n)\in\N^2$ and $\tau\in\Rp$.  We define the $\tau$-auxiliary BTRW SIE associated with \eqref{trnctdsie} on $[0,\tau]\times\Xnd$ by
\beq\lbl{aux}\tag{Aux}
\xttauxnlf=\bc
\ds\uttxnD+\ds\sum_{y\in\Xnldf}\int_0^t\kappa^{x,y}_{\h,s,\tau}\lpa\xstauynlf\rpa{d\wsyn}; &x\in\Xnldf,\\
\ds\uttxnD;&x\in\Xnd\setminus\Xnldf
\ec
\eeq
where the independent BMs sequence $\lbr W_n^y\rbr_{y\in\Xnldf}$ in \eqref{aux} is the same for all $\tau>0$, as well as $x\in\Xnldf$.
We denote \eqref{aux} by $\btrwsienlaux$.  We say that the pair of families $\lpa\lbr\xtaunlf\rbr_{\tau\ge0},\lbr W_n^y\rbr_{y\in\Xnldf}\rpa$ is a solution to $\lbr\btrwsienlaux\rbr_{\tau\ge0}$ on a filtered probability space $\OFFtP$ if
there is one family of independent BMs (up to indistinguishability) $\lbr W_n^y(t);0\le t<\infty\rbr_{y\in\Xnldf}$ on $\OFFtP$ such that, for every fixed $\tau\in\Rp$
\begin{enumerate}\renewcommand{\labelenumi}{(\alph{enumi})}
\item
the process $\lbr\xttauxnlf,\sFt;0\le t\le\tau, x\in\Xnd \rbr$
has continuous sample paths in $t$ for each fixed $x\in\Xnd$ and $\xttauxnlf\in\sFt$ for all $x\in\Xnd$ for every $0\le t\le\tau$; and
\item equation \eqref{aux} holds on $[0,\tau]\times \Xnd$, $\P$-almost surely.
\end{enumerate}
Naturally, implicit in our definition above the assumption that, for each fixed $\tau\in\Rp$, we have
$$\P\lbk\int_0^t\lpa\kappa^{x,y}_{\h,s,\tau}\lpa\xstauynlf\rpa\rpa^{2}{ds}<\infty\rbk=1;\ \forall x,y\in\Xnldf,0\le t\le\tau.$$
For simplicity, we will sometimes say that $\xtaunlf=\lbr\xttauxnlf,\sFt;0\le t\le\tau, x\in\Xnd \rbr$ is a solution to \eqref{aux} to mean the above.  Clearly, if $\xttauxnlf$ satisfies \eqref{aux} then $\uttauxnlf:=\xtautauxnlf$
satisfies \eqref{trnctdsie} at $t=\tau$ for all $x\in\Xnd$.  Also, for each $n$ and each $1\le d\le3$
$$\lab\kappa^{x,y}_{\h,s,\tau}\lpa\xstauynlf\rpa\rab=\lab\df{\Ktausxyn}{\sdh} a(\xstauynlf)\rab\le\df{\lab a\lpa\xstauynlf\rpa\rab}{\sdh}.$$
In addition, for each fixed $\tau\in\Rp$ and each fixed $x,y\in\Xnldf$ we have for a solution $\xtaunlf$ to \eqref{aux} that
$$\kappa^{x,y}_{\h,s,\tau}\lpa\xstauynlf\rpa\in\sFs;\quad \forall s\le\tau,$$
since, of course the deterministic ${\Ktausxyn}/{\sdh}\in\sFs$ and $a(\xstauynlf)\in\sFs$.
Thus, if $\xtaunlf$ solves \eqref{aux}; then, for each fixed $\tau>0$ and $x,y\in\Xnldf$, each stochastic integral in \eqref{aux}
$$I_{n,l}^{\tau,x,y}=\lbr I_{n,l}^{\tau,x,y}(t):=\int_0^t\kappa^{x,y}_{\h,s,\tau}\lpa\xstauynlf\rpa{d\wsyn},\sFt;\ 0\le t\le\tau\rbr$$
is a continuous local martingale in $t$ on $[0,\tau]$.  This is clear since by a standard localization argument we may assume the boundedness of $a$ ($|a(u)|\le C$); in this case we have for each fixed $x,y\in\Xnldf$ and $\tau\in\Rp$ that
\beqs
\bsp
\E\lbk I_{n,l}^{\tau,x,y}(t)\bigg|\sFr\rbk=\int_0^r\kappa^{x,y}_{\h,s,\tau}\lpa\xstauynlf\rpa{d\wsyn}=I_{n,l}^{\tau,x,y}(r),\ r\le t\le\tau.
\end{split}
\eeqs
So, the finite sum over $\Xnldf$ in \eqref{aux} is also  a continuous local martingale in $t$ on $[0,\tau]$.  I.e., for each $\tau>0$ and $x\in\Xnldf$
$$M_{n,l}^{\tau,x}=\lbr M_{n,l}^{\tau,x}(t):=\sum_{y\in\Xnldf}\int_0^t\kappa^{x,y}_{\h,s,\tau}\lpa\xstauynlf\rpa{d\wsyn},\sFt;\ 0\le t\le\tau\rbr\in\locm$$
with quadratic variation
\beq\lbl{qv}
\lqv M_{n,l}^{\tau,x}(\cdot)\rqv_t=\sum_{y\in\Xnldf}\int_0^t\lbk\kappa^{x,y}_{\h,s,\tau}\lpa\xstauynlf\rpa\rbk^2ds
\eeq
where we have used the independence of the BMs $\lbr W_n^y\rbr_{y\in\Xnldf}$ within the lattice $\Xnldf$.
For each $\tau>0$, we call $M_{n,l}^{x,\tau}$ a kernel local martingale (or K-local martingale).

There is another complicating factor in formulating our K-martingale problem approach that is not present in the standard SDEs setting.
To easily extract solutions to the truncated BTRW SIEs in \eqref{trnctdsie} from the family of auxiliary problems
$\lbr\btrwsienlaux\rbr_{\tau>0}$ in \eqref{aux}, we want the independent BMs sequence
$\lbr W_n^y\rbr_{y\in\Xnldf}$ to not depend on the choices of $\tau$ and $x$.   I.e., we want all the K-local martingales in \eqref{aux}
to be stochastic integrals with respect to the same sequence $\lbr W_n^y\rbr_{y\in\Xnldf}$, regardless of $\tau$ and $x$.
With this in mind, we now formulate the K-martingale problem associated with the auxiliary BTRW SIEs in \eqref{aux}.  Let
\beq\lbl{Zspace}
\C_{n,l}:=\lbr u:\Rp\times\lpa\Xnldf\rpa^2\to\Rs; t\mapsto u^{x_1,x_2}(t) \mbox{ is continuous }\forall x_1,x_2\rbr.
\eeq
For $u\in\C_{n,l}$, let $u^{x_1,x_2}(t)=\lpa u_1^{x_1}(t),u_2^{x_2}(t)\rpa$ with $u^x(t)=u^{x,x}(t)$;
and for any $\tau_1,\tau_2>0$ and any $x_1,x_2,y\in\Xnldf$ let
\beq\lbl{ups}
\Upsilon^{x_{i,j},y}_{\h,t,\tau_{i,j}}\lpa u^y(t)\rpa:=\df{\Ktauitxiyn}{\sdh} a(\uity) \df{\Ktaujtxjyn}{\sdh}a(\ujty);\quad 1\le i,j\le2,
\eeq
(we are allowing the cases $\tau_1=\tau_2$ and/or $x_1=x_2$)
where for typesetting convenience we denoted the points $(\tau_i,\tau_j)$ and $(x_i,x_j)$ by $\tau_{i,j}$ and $x_{i,j}$, respectively.
We denote by $\pa_i$ and $\pa^2_{ij}$ the first order partial derivative with respect to the $i$-th argument and the second
order partials with respect to the $i$ and $j$ arguments, respectively.
Let $\C^2=\C^2(\Rs;\R)$ be the class of twice continuously differentiable real-valued functions on $\Rs$  and let
\beq\lbl{twicwcontbdd}
\C_b^2=\lbr f\in\C^2; \mbox{$f$ and its derivatives up to second order are bounded}\rbr.
\eeq
Now, for $\tau_1,\tau_2>0$, for $f\in\C^2_b$, and for $(t,x_1,x_2,u)\in[0,\tau_1\wedge\tau_2]\times\lpa\Xnldf\rpa^2\times \C_{n,l}$ let
\beq\lbl{siegen}
\bsp
\lpa\sAKtauot f\rpa(t,x_1,x_2,u)&:=\sum_{1\le i\le2}\pa_if\lpa u^{x_1,x_2}(t)\rpa \frac{\pa}{\pa t}\uttxinD
\\&+ \frac12\sum_{1\le i,j\le2} \pa^2_{ij}f\lpa u^{x_i,x_j}(t)\rpa\sum_{y\in\Xnldf}\Upsilon^{x_{i,j},y}_{\h,t,\tau_{i,j}}\lpa u^y(t)\rpa
\end{split}
\eeq
Let $X_{n,l}^\tau=\lbr\xttauxnlf;0\le t\le\tau,x\in\Xnd\rbr$ be a continuous in $t$ adapted real-valued process on a filtered probability space $\OFFtP$.
For every $\tau_1,\tau_2>0$ define the two-dimensional stochastic process $Z_{n,l}^{\tau_{1,2}}$:
\beq\lbl{Z}
\lbr Z_{n,l}^{x_{1,2},\tau_{1,2}}(t)=\lpa\xttauoxnlf,\xttautxnlf\rpa; (t,x_1,x_2)\in[0,\tau_1\wedge\tau_2]\times\lpa\Xnldf\rpa^2\rbr
\eeq
with $Z_{n,l}^{y,\tau_{1,2}}(t)=\lpa\xttauoynlf,\xttautynlf\rpa$ and let $U_0^{x_1,x_2}=(\unxo,\unxt)$.  We say that the family $\lbr X_{n,l}^\tau\rbr_{\tau\ge0}$ satisfies the K-martingale problem associated with the auxiliary
BTRW SIEs in \eqref{aux} on $\Rp\times\Xnd$ if for every $ f\in\C_b^2$, $0<\tau_1,\tau_2<\infty$, $\tau=\tau_1\wedge\tau_2$, $t\in[0,\tau]$, $x_1,x_2\in\Xnldf$, and
$x\in\Xnd\setminus\Xnldf$ we have
\beq\lbl{kmart}\tag{KM}
\bc
f(Z_{n,l}^{x_{1,2},\tau_{1,2}}(t))-f(U_0^{x_1,x_2})-\ds\int_0^t \lpa\sAKtauot f\rpa(s,x_1,x_2,Z_{n,l}^{\tau_{1,2}})ds\in\locm;\\
\ds\xttauxnlf=\uttxnD.
\ec
\eeq
We are now ready to state the equivalence of the K-martingale
problem in \eqref{kmart} to the auxiliary SIEs in \eqref{aux} and its implication for the BTRW SIE in \eqref{trnctdsie}.  This result is of independent interest and is stated as the following theorem\footnote{This is because it is easily adaptable to many mild formulations of SPDEs, of different orders, not just for the BTBM SIEs (see also the remarks following \thmref{kmarteqaux}).  Since we don't prove uniqueness under less than Lipschitz conditions for our BTBM SIE, we have not explicitly mentioned the uniqueness implications of our K-martingale approach.  More on that in future articles.}.
\bthm\lbl{kmarteqaux}
The existence of a solution pair $\lpa\lbr\xtaunlf\rbr_{\tau\ge0},\lbr W_n^y\rbr_{y\in\Xnldf}\rpa$ to the $\tau$-auxiliary BTRW SIEs $\lbr\btrwsienlaux\rbr_{\tau\ge0}$ in \eqref{aux} on a filtered probability space $\OFFtP$
is equivalent to the existence of a family of processes $\lbr\xtaunlf\rbr_{\tau\ge0}$ satisfying \eqref{kmart}.  Furthermore, if there is $\lbr\xtaunlf\rbr_{\tau\ge 0}$ satisfying
\eqref{kmart} then there is a solution to \eqref{trnctdsie} on $\Rp\times\Xnd$.
\ethm
The versatility of the K-martingale approach and the fact that it represents a unified way of dealing with many SPDEs that are of different orders
is now clear. The kernel of our BTRW in \eqref{aux} and \eqref{kmart} may be replaced by the discretized version
of the linearized Kuramoto-Sivashinsky kernel (the spatially discretized version of the IBTBAP kernel in \eqref{ibtbapkernelc}) or by the density of a random walk to handle in a unified method
the fourth order linearized KS \eqref{lks} and other related SPDEs or the second order RD SPDEs.  Only minor and obvious modifications are needed to apply this approach to Burgers-type SPDEs.  It can also be adapted to treat Navier-Stokes SPDEs driven by space-time white noise and many hyperbolic SPDEs.
\subsection{BTBM SIEs and their fourth order SPDEs links}\lbl{spdeslinks}
In this section, we briefly discuss the quite interesting connections---both direct and indirect---of the BTBM SIEs $\btsie$ to
fourth order SPDEs driven by space-time white noise.  We start by giving a quick glimpse into the indirect link of
$\btsie$ to Kuramoto-Sivashinsky type SPDEs via our imaginary-Brownian-time-Brownian-angle process (IBTBAP) representation
of the linearized KS PDE (see \cite{Aks,Abtpspde}).  A more extensive treatment of such SPDEs using this IBTBAP approach
is presented in \cite{Aksspde,AL,AD}.   We then end this subsection by giving a connection of
the BTBM SIEs $\btsie$ to an unconventional new parametrized fourth order BTP SPDEs using the spatially discretized versions of these equations.
\subsubsection{BTBM SIEs are cousins of the Kuramoto-Sivashinsky and related SPDEs}\lbl{kslink}
By replacing $\KBtxy$ and $\KBtsxy$ in \eqref{btpsol} with the intimately connected kernels $\KKStxy$ and $\KKStsxy$,
defined by
\beq\lbl{ibtbapkernelc}\tag{KSK}
\bc
\Mpisxy\eqdef\ds\df{\exp\lpa is\rpa}{{\lpa2\pi is \rpa}^{d/2}}e^{-|x-y|^2/2is},\cr
\KKStxy\eqdef\ds\int_{-\infty}^0\Mpisxy\ptzs ds+ \int_{0}^\infty\Mpisxy\ptzs ds
\ec
\eeq
we obtain our definition of IBTBAP solutions
\beq\lbl{ibtbapsol}
\bsp
U(t,x)=\intrd\KKStxy \uny dy+ \intrdzt \KKStsxy a(\usy)\sW(ds\times dy)
\end{split}
\eeq
to the canonical Kuramoto-Sivashinsky SPDE with linearized PDE part:
\begin{equation} \label{lks}
 \begin{cases} \displaystyle\frac{\partial U}{\partial t}=
-\frac18\lap^2U-\frac12\lap U-\frac12U+a(U)\frac{\partial^{d+1} W}{\partial t\partial x}, & (t,x)\in(0,+\infty )\times\Rd;
\cr U(0,x)=\unx, & x\in\Rd.
\end{cases}
\end{equation}
This IBTBAP representation approach for the SPDE \eqref{lks}
is inspired by our earlier work \cite{Aks}, in which we used the deterministic version of \eqref{ibtbapsol} ($a\equiv0$)
to solve the linearized KS PDE obtained from \eqref{lks} by setting $a\equiv0$.  It easily allows for the addition of a nonlinear term to \eqref{lks}.  The nonlinearity  could be an Allen-Cahn one (to get Swift-Hohenberg SPDE), a KPP one (a new KS-type SPDE), a Burgers  one (versions and variants of KS SPDE), and many more interesting nonlinearities.  Quantum mechanics experts will immediately note that, except for the $\exp{(is)}$ term, $\Mpisxy$ in the definition of the IBTBAP kernel in \eqref{ibtbapkernelc} is a $d$-dimensional version of the free propagator associated with Schr\"odinger equation.

One reason why the BTBM SIEs are cousins of the Kuramoto-Sivashinsky and related SPDEs is that the
kernel $\KKStxy$ above may be regarded as the ``density'' of the IBTBAP,  as in \cite{Aks}.  By construction, the
IBTBAP---which we also called the linearized Kuramoto-Sivashinsky
process or LKSP---is intimately connected to BTPs (e.g. \cite{Abtp1,Abtp2,Aks}).  Moreover, as we showed in \cite{Abtp1,Abtp2,Aks} for the PDEs case,
the kernels  $\KBtxy$ and $\KKStxy$ have similar regularizing effects on their corresponding equations.  Analogously, our solutions to the BTBM SIEs \eqref{btpsol} have similar $d$-dependent regularity properties
to the IBTBAP mild solutions  \eqref{ibtbapsol} to KS-type SPDEs, including \eqref{lks} and nonlinear versions of it,
for $d=1,2,3$  (see \cite{Aksspde} and followup papers).
We are also currently using this approach to investigate asymptotic and other qualitative behaviors of several related
nonlinear SPDEs in applied mathematics (e.g., \cite{AL,AD}).
\subsubsection{BTBM SIEs and the parametrized BTBM SPDEs}\lbl{btspdelink}
We now link  $\btsie$ to a new fourth order parametrized SPDE on discrete spatial lattices.  First,
we note that this direct link of BTBM SIEs to SPDEs is not as straightforward as one might be tempted to believe;
and it cannot be resolved by using a conventional SPDE, as we explain below.  Even though
it is certainly true that the deterministic term on the right hand side of \eqref{btpsol}
\beq\lbl{Asol}
v(t,x):=\intrd\KBtxy \uny dy
\eeq
solves, and is the kernel integral form of, the deterministic fourth order PDE \eqref{btppdedet}
(see Allouba et al.~\cite{Abtp2,Abtp1}), it turns out that the BTBM SIE $\btsie$
is {\it different} from---i.e., not the mild stochastic integral form of---the naturally guessed SPDE \eqref{btpspde} (see \cite{Abtpspde}).
This is due precisely to the Laplacian acting on the initial data $\un$.   Instead, \eqref{btpspde} may be viewed as a degenerate version of the parametrized BTBM SPDEs linked here to our BTBM SIEs (see Appendix \ref{appC} for a brief discussion involving the kernel formulation of \eqref{btpspde} in terms of its spatially-discretized version)\footnote{So, calling solutions to $\btsie$ the BTP solutions to \eqref{btpspde}---as we did in \cite{Abtpspde}---is not precise, and we call $\btsie$ the
BTBM SIE instead. We study BTBM SPDEs and related stochastic fractional PDEs in upcoming articles.}.     Our main interest here is in the BTBM SIEs, and here also we use the spatially-discretized version of our BTBM SIE $\btsie$ (the BTRW SIE $\btrwsien$ in \eqref{btrwsie}) to connect it to fourth order SPDEs.  Figuring out this correct and subtle SPDE-link to $\btsie$ via our easier-to-see discretized versions of the equations is another advantage of this multiscale approach over the direct one.

To intuitively see the correct SPDE link to $\btsie$, we go back to the deterministic BTP PDE case \eqref{btppdedet},
and observe again that the solution in \eqref{Asol} $v\in \mathrm{C}^{1,4}(\Rp\times\Rd;\R)$
is indeed very smooth for \textit{all} times and \textit{all} spatial dimensions $d\ge1$ (see \cite{Abtp2,Abtp1})
despite the presence of the positive biLaplacian.
In order to heuristically pinpoint the cause of the smoothing effect of the BTBM kernel $\KBtxy$ in \eqref{Asol}
and connect it to the corresponding PDE \eqref{btppdedet}, we must look at both
terms $\Delta^2u$ and $\Delta u_0/\sqrt{8\pi t}$ {\it together}.  We observe that the smoothing effect of the Laplacian term in
\eqref{btppdedet} gets arbitrarily large as the time $t$ in $\KBtxy$ goes to zero ($t\searrow0$) and fades as $t\nearrow\infty$ at the rate of $1/\sqrt{8\pi t}$;
and the Laplacian is acting on the smooth initial solution $u_0$ (the solution
at time $t=0$).  Heuristically, this suggests a recipe for obtaining eternally-smooth solutions involving the positive biLaplacian
coupled with the smoothing Laplacian.

Turning our attention now towards the BTBM SIE $\btsie$.  We first give a heuristic prelude to the new notion of parametrized BTBM SPDEs, and then we precisely state the result.  We know from the results of Sections \ref{intro} and \ref{pfs}, (as well as those in Appendix \ref{appB} and the comments in \secref{conc}) that the BTBM density $\KBtxy$ has a significant smoothing effect on solutions to $\btsie$ as compared to the Brownian motion density in the mild formulation of standard second order RD SPDEs.
We also see that the stochastic white noise term in $\btsie$
involves the BTBM density at $t-s$, viz.~$\KBtsxy$.  So, based on the heuristic above and generalizing it, any SPDE that
captures the $\KBtsxy$ smoothing effect will include a positive bi-Laplacian term along with a Laplacian term
whose coefficient grows arbitrarily large as $s\nearrow t$ at the rate of $1/\sqrt{8\pi(t-s)}$.  On one hand, this Laplacian will have to act
on the solution of the SPDE for all times $r\le s$; on the other hand it has to also act on the solution at spatial points $x=y$, since $\KBtsxy=0$ except at $x=y$ when $t-s=0$.  That is, we need to keep track of four parameters in both time and space (not just $(t,x)$) to encode the smoothing effect of the BTBM kernel in \eqref{btpsol} into an SPDE; and we are led to the notion of
parametrized SPDEs associated with our $\btsie$.  We give this link on the lattice in \lemref{tsddesie}, and we prove it in \secref{linkpf} for spatially discretized versions of $\btsie$.  More precisely, we discretize $\Rd$ into the lattices $\Xnd$.  We then fix an arbitrary $n\in\N$ and consider the following parametrized stochastic differential-difference equation (PSDDE), written in integral form as:
\beq
\lbl{psdde}
\bsp
&\utstxyn-\uny=\int_0^s\lbk\df{\Dn\utrtxyn\big|_{y=x}}{\sqrt{8\pi(t-r)}}+\df18\Dns\utrtxyn\rbk dr\\
&+\int_0^sa(\utryn)\df{d\wryn}{\sdh};\ 0\le s\le t<\infty,\ x,y\in\Xnd,
\end{split}
\eeq
which we denote by $\pbtsdden$, where $\Delta_n$ and $\Dns$ are the $d$-dimensional discrete Laplacian and bi-Laplacian on $\Xnd$, acting on the second spatial argument respectively:
\beq
\bsp
&\Dn\utrtxyn\big|_{y=x}=\sum_{i=1}^{d}\Delta_{n,i}\utrtxyn\big|_{y=x}
\\&:=\sum_{i=1}^{d} \df{\tilde{U}_{n}^{x;x_{1},\ldots, x_{i}+\h,\ldots, x_{d}}(r,t)-2\tilde{U}_{n}^{x,x}(r,t)+ \tilde{U}_{n}^{x; x_{1},\ldots, x_{i}-\h,\ldots, x_{d}}(r,t)}{\h^{2}};
\\&\Dns\utrtxyn=\sum_{i,j=1}^{d}\Delta_{n,j}\Delta_{n,i}\utrtxyn
\end{split}
\eeq
and where $\utryn:=\tilde{U}^{y,y}_n(r,r)$, and $\tilde{U}_n^{x,y}(0,t)=\uny$ for all $(t,x,y)\in\Rp\times\Xnds$.
By a solution to the PSDDE system $\lbr\pbtsdden\rbr_{n=1}^\infty$ on $\Rp\times\Xnd$ with respect to the Brownian $($in $t$\/$)$ system $\lbr\wtyn\rbr_{(n,y)\in\N\times\Xnd}$ on the filtered probability space $\OFFtP$ we mean a sequence of real-valued processes $\lbr\tilde{U}_{n}\rbr_{n=1}^\infty$  with continuous sample paths in $s$ for each fixed $(t,x,y)\in\Rp\times\Xnds$ and $n\in\N$
such that, for every $(n,t,x,y)\in\N\times\Rp\times\Xnds$, $\utstxyn$ is $\sFs$-adapted, and equation \eqref{psdde} holds $\P$-a.s.  In particular,
$$\P\lbk \int_0^s\lpa\lab\df{\Dn\utrtxyn\big|_{y=x}}{\sqrt{8\pi(t-r)}}+\df18\Dns\utrtxyn\rab +\frac{a^{2}(\utryn)}{\hd}\rpa dr<\infty\rbk=1$$
holds for every $0\le s\le t<\infty,\ x,y\in\Xnd$.

We then prove (see \secref{linkpf}) that if $\utstxyn$ solves \eqref{psdde} then $\uttxn$
solves the spatially-discretized version of $\btsie$---the BTRW SIE $\btrwsien$ given by \eqref{btrwsie}.
We stress here that it is not enough for the ``diagonal terms''  $\utttxxn=\uttxn$ to satisfy \eqref{psdde} (the special case
of \eqref{psdde} $s=t$ and $y=x$) to conclude that $\uttxn$ satisfies the BTRW SIE in \eqref{btrwsie}:
{\it all} of the $\utstxyn$ must satisfy the PSDDE \eqref{psdde} for us to have this implication
(see \lemref{tsddesie} and its proof in \secref{linkpf}).

It is in the above spatially-discretized sense that we say that the BTBM SIE $\btsie$
is associated with the parametrized BTBM SPDE
\beq\lbl{pspde}
\bsp
\df{\pa U_t^x(s,y)}{\pa s}&=\lbk\df{\D U_t^x(s,y)\big|_{y=x}} {\sqrt{8\pi(t-s)}}+\df18\Ds U_t^x(s,y)\rbk ds +a(U_s^y(s,y))\df{\pa^{d+1}W}{\pa s\pa y}
\\
U_t^x(0,y) &= \uny
\end{split}
\eeq

The system of parametrized  stochastic differential-difference equations (PSDDEs) \eqref{psdde},
may also be written in differential form as:
\beq\lbl{dbtrwsdde}
\bc d\utstxyn=\lbk\df{\Dn\utstxyn\big|_{y=x}}{\sqrt{8\pi(t-s)}}+\df18\Dns\utstxyn\rbk ds +a(\utsyn)\df{d\wsyn}{\sdh};& \cr
\utztxyn = \uny; &
\ec
\eeq

The proof of the following lemma follows from an application of It\^o's rule.  The interesting point here is the unorthodox nature of the parametrized SDDE connected to the lattice version of our BTBM SIE $\btsie$.
\blm[Relation between $\btrwsien$ and $\pbtsdden$]\lbl{tsddesie}
Fix an arbitrary $n\in\N$ and assume the conditions in \eqref{cnd} and let
$$\utn:=\lbr\utstxyn;(s,t,x,y)\in[0,t]\times\Rp\times\Xnds\rbr$$
be a continuous-in-$s$ solution to the parametrized SDDEs \eqref{psdde} such that,
for any fixed pair $(t,x)$, $\E\lab\utstxyn\rab^2\le C$ for all $(s,y)\in[0,t]\times\Xnd$ for some constant $C>0$.
Then $\uttxn:=\utttxxn$ solves the BTRW SIE \eqref{btrwsie}.
\elm
\brm\lbl{momcndonpsddes}
The moment boundedness condition in \lemref{tsddesie} is technical and may be relaxed.  Also, in light of the fact that \lemref{tsddesie}  holds for all $n$, we may call the BTRW SIE $\btrwsien$ the lattice-kernel-diagonal form of the PSDDE \eqref{psdde}; and we may then call the BTBM SIE $\btsie$ the limiting-lattice-kernel-diagonal form of the parametrized BTBM SPDE \eqref{pspde}.  A converse of \lemref{tsddesie} is given in Appendix \ref{appC} in \lemref{converse}.
\erm

\section{Proof of results}\lbl{pfs}
\subsection{Density regularity estimates}\lbl{destm}
The first set of estimates\footnote{As is customary, $C$ will denote a constant that may change its value from one line to the next.  We will denote the Euclidean norm on $d$-dimensional spaces by $\lab\cdot\rab$.} we need are bounds on the square of the Brownian-time Brownian motion density $\KBtxy$ and its associated lattice Brownian-time random walk density $\Ktxn$ and their temporal and spatial differences.  We obtain these estimates for both kernels simultaneously. The method of proof is to reduce, via an asymptotic argument, these estimates for the BTRW to the corresponding ones for the BTBM density $\KBtx$ and perform the computations in the setting of the BTBM.  Since all the results in this part hold for all $n\ge N^*$ (equivalently for all $\h\le\delta_{N^*}$) for some positive integer $N^*$, we will suppress the dependence on $n$, except when it is needed or helpful, to simplify the notation.   Also, whenever we need these estimates, we assume that  $n\ge N^*$ without explicitly stating it every time; and when we do, we let
\beq\lbl{star}
\Ns:=\lbr n\in\N;n\ge N^*\rbr
\eeq

We start by observing that in the classical setting of Brownian motion and its discretized version continuous-time random walk on $\Xnd=\h\Zd$, we have  the following well known asymptotic result relating their densities (see e.g., \cite{SZ})
\beq\lbl{as1}
\qtxynint\sim\ptxy\hd\mbox{ as } n\to\infty\mbox{ (as $\h\to0$); }\forall t>0,\ x,y\in\Rd,
\eeq
where for each $x\in\Rd$ we use $[x]_{\h}$ to denote the element of $\Xnd$ obtained by replacing each coordinate $x_{i}$ with $\h$ times the integer part of $\h^{-1}x_{i}$, and $a_n\sim b_n$ as $n\to\infty$ means $a_n/b_n\to1$ as $n\to\infty$.
Now, for every continuous and bounded $\un:\Rd\to\R$, we have
\beq\lbl{minusx}
\lim_{\h\searrow0}\sum_{y\in\Xnd\setminus\lbr x\rbr} \KBtxy \uny\hd=\intrd\KBtxy\uny dy;\  t>0,\ x\in\Rd,\ d\ge1,
\eeq
and by the dominated convergence theorem we obtain
\beq\lbl{BTRWsol2BTBMsol}
\bsp
&\lim_{\h\searrow0}\lab\sum_{y\in\Xnd}\Ktxynint\uny-\sum_{y\in\Xnd\setminus\lbr x\rbr} \KBtxy \uny\hd\rab\\
&=\lab\int_{0}^{\infty}\lbr\lim_{\h\searrow0}\sum_{y\in\Xnd\setminus\lbr x\rbr}\lbk\qsxynint-\psxy\hd\rbk\uny\rbr2\ptzs ds\rab=0
\end{split}
\eeq
for $t>0$, $x\in\Rd$, and $d\ge1$; since, by \eqref{as1},
$$\lim_{\h\searrow0}\sum_{y\in\Xnd}\qsxynint\uny=\lim_{\h\searrow0}\sum_{y\in\Xnd}\psxy\uny\hd=\intrd\psxy\uny dy$$
for every $(s,x)\in(0,\infty)\times\Rd$.  We then straightforwardly get the following result.
\blm\lbl{btrwasbtp}  For every continuous and bounded $\un:\Rd\to\R$ and for every $d\ge1$
\beq\lbl{ddegreen2pdegreen}
\lim_{\h\searrow0}\sum_{y\in\Xnd}\Ktxynint\uny=\intrd\KBtxy\uny dy;\forall (t,x)\in(0,\infty)\times\Rd,
\eeq
and the following asymptotic relation holds between the BTP and BTRW densities:
\beq\lbl{btrwbtp}
\Ktxynint\sim\KBtxy\hd \mbox{ as } n\to\infty\ (\mbox{as }\h\to0);\ t>0,\ x,y\in\Rd,\ x\neq y.
\eeq
\elm
\brm\lbl{ddetopde}
Equation \eqref{ddegreen2pdegreen} confirms the intuitively clear fact that the kernel form of the BTRW DDE \eqref{btrwdde}
converges pointwise---as $\h\searrow0$---to the kernel form of its continuous version, the BTBM PDE \eqref{btppdedet}.  We also remind the reader that the right hand side of \eqref{ddegreen2pdegreen} is in $\mathrm{C}^{1,4}$ for all $(t,x)\in(0,\infty)\times\Rd$ under the conditions \eqref{cnd} on $\un$.
 \erm
By reducing the computations to the setting of the Brownian-time Brownian motion, using \lemref{btrwasbtp} together with scaling, our method of proof for the next three lemmas shows that the estimates in these lemmas all hold for the BTBM density as well as for the BTRW one with obvious changes from the discrete to the continuous settings (see \lemref{2ndQinequality}, \lemref{4thQinequality}, and \lemref{3rdQinequality} below).  Thus, these lemmas are stated for both densities $\Ktx$ and $\KBtx$.   This, in turn, allows us to easily indicate how to prove the same H\"older regularity of solutions to the BTBM SIE $\btsie$ directly (without discretization) in the case of Lipschitz conditions (see \secref{pf1stmain} for the detailed existence, uniqueness, $L^{p}(\Omega)$-boundedness and  the H\"older regularity proof of \thmref{lip}).  We start with
\begin{lem}\label{2ndQinequality}
There are constants $C$ and $\tilde{C}$, depending only on $d$, and a $\d^*>0$ such that for all $\d\le\d^*$
\[\intrd\lbk\KBtx\rbk^2dx=\df{C}{t^{d/4}}\mbox{ and }\sum_{x\in\Xd} \lbk\Ktx\rbk^2 \leq \tilde{C}\df{\dd}{t^{d/4}};\quad t>0,\ 1\le d\le3\]
{ and hence }
\[\int_0^t\intrd\lbk\KBsx\rbk^2dx ds=Ct^{\tf{4-d}{4}}\mbox{ and } \int_0^t\sum_{x\in\Xd} \lbk\Ksx\rbk^2  ds \leq \tilde{C}\dd t^{\tf{4-d}{4}};\]
for all $t>0,\ 1\le d\le3$.
\end{lem}
\begin{proof}
Using \lemref{btrwasbtp} we obtain
\beq\lbl{allcomp}
\bsp
&\lim_{\d\searrow0}\sum_{x\in\Xd} \df{\lbk\Ktx\rbk^2}{\dd}=\intrd\lbk\KBtx\rbk^{2} dx\\
&=4\lbk\int_0^\infty\int_0^\infty\lbk\int_{\Rd} \psoxz \pstxz dx\rbk \ptzso \ptzst ds_1ds_2\rbk\\
&=4\lbk\int_0^\infty\int_0^\infty\lbk\df{1}{\lbk2\pi (s_1+s_2)\rbk^{d/2}}\rbk \ptzso \ptzst ds_1ds_2\rbk\\
&=4\lbk\int_0^{\pi/2}\int_0^\infty\df{e^{-\rho^2/2t}}{(2\pi)^{\tf{d+2}{2}}
\lbk\rho (\sin(\theta)+\cos(\theta))\rbk^{d/2}t}\rho d\rho d\theta\rbk\\
&=\bc\df{4C_{d}}{t^{d/4}};& 1\le d\le 3,\cr\infty;&d\ge4,\ec
\end{split}
\eeq
where $C_{d}$ is dimension-dependent\footnote{$C_{1}\approx 0.0914$, $C_{2}\approx0.0396$, and $C_{3}\approx 0.0243$.}.  Then there is a $\d^{*}>0$ such that, whenever $\d\le\d^{*}$, we obtain
$$\df1\dd\sum_{x\in\Xd} \lbk\Ktx\rbk^2\le\df{\tilde{C}_{d}}{t^{d/4}};\ 1\le d\le 3,$$
with a constant $\tilde{C}_{d}>4C_{d}$.  The last assertion of the lemma trivially follows upon integration over the time interval $(0,t]$.
\end{proof}

The following lemma is key to our H\"older regularity result in time for $1\le d\le3$.  We give a probabilistically-flavored proof using the notion of $2$-Brownian-times random walk and $2$-Brownian-times Brownian motion given below.
\begin{lem}\label{4thQinequality}
There is a constant $C$, depending only on $d$, and a $\d^*>0$ such that for $\d\le\d^*$
\beq\lbl{tmpkernel}
\bc
\ds\int_0^t\intrd{\lbk\KBtsx - \KBrsx\rbk}^2 dx ds \leq C(t-r)^{\tf{4-d}{4}};&\hspace{-2mm}0<r<t,1\le d\le3,\\
\ds\int_0^t\sum_{x\in\Xd} {\lbk\Ktsx - \Krsx\rbk}^2 ds \leq {C}\dd(t-r)^{\tf{4-d}{4}};&\hspace{-2mm}0<r<t, 1\le d\le 3,
\ec
\eeq
with the convention that $\Ktx=0=\KBtx$ if $t<0$.
\end{lem}
\bpf
We will prove that
\beq\lbl{tmphl}
\int_0^t\sum_{x\in\Xd} {\lbk\Ksptmrx - \Ksx\rbk}^2 ds \leq C\dd(t-r)^{\tf{4-d}{4}}; \quad 1\le d\le 3.
\eeq
for all $\delta\le\delta^{*}$, for some $\d^{*}>0$, simultaneously with its corresponding BTBM density statement.  The first step is to show the identity
\beq\lbl{ido}
\sum_{x\in\Xd} {\lbk\Ksptmrx - \Ksx\rbk}^2=\KBBsptmrz+\KBBsz-2\KBBsptmrsz
\eeq
where
\beq
\KBBuvz=4\int_{0}^{\infty}\int_{0}^{\infty}\qroprtz\puzro\pvzrt dr_{1} dr_{2}
\eeq
is the density of the $2$-Brownian-times random walk
\beq\lbl{twobtrwdef}
\S^0_{B^{(1)},B^{(2)},\h}(u,v)\eqdef S^0_{\h}\lpa\lab B^{(1)}_u\rab+\lab B^{(2)}_v \rab\rpa;\quad 0\le u,v<\infty,
\eeq
in which the $d$-dimensional random walk $S^0_{\h}$ (on $\Xnd$) and the two one-dimensional BMs $B^{1}$ and $B^{2}$ are all independent.
But,
\beq\lbl{idz}
\bsp
&\sum_{x\in\Xd} \Kux\Kvx\\&=4\int_{0}^{\infty}\int_{0}^{\infty}\lbk\sum_{x\in\Xd}\qrox\qrtx\rbk\puzro\pvzrt dr_{1} dr_{2}\\
\\&=4\int_{0}^{\infty}\int_{0}^{\infty}\qroprtz\puzro\pvzrt dr_{1} dr_{2}=\KBBuvz.
\end{split}
\eeq
The identity \eqref{ido} immediately follows from \eqref{idz}.  Similarly,  we get the corresponding identity for the BTBM setting
\beq\lbl{idtwo}
\int_{\Rd} {\lbk\KBsptmrx - \KBsx\rbk}^2 dx=\KBBBsptmrz+\KBBBsz-2\KBBBsptmrsz
\eeq
where
\beq
\KBBBuvz=4\int_{0}^{\infty}\int_{0}^{\infty}\proprtz\puzro\pvzrt dr_{1} dr_{2}
\eeq
is the density of the $2$-Brownian-times Brownian motion
\beq\lbl{2btbmdef}
\X^0_{B^{(1)},B^{(2)}}(u,v)\eqdef X^0\lpa\lab B^{(1)}_u\rab+\lab B^{(2)}_v \rab\rpa;\quad 0\le u,v<\infty,
\eeq
in which the $d$-dimensional BM $X^{0}$ and the two one dimensional BMs $B^{(1)}$ and $B^{(2)}$ are all independent.
Using the identities \eqref{ido} and \eqref{idtwo}, along with a similar asymptotic argument to the one we used in the proof of
\lemref{2ndQinequality} together with the dominated convergence theorem, yield
\beq\lbl{hlbd}
\bsp
&\lim_{\d\searrow0}\frac{1}{\dd}\lbk\int_{0}^{t}\KBBsptmrz ds + \int_{0}^{t} \KBBsz ds-2\int_{0}^{t} \KBBsptmrsz ds\rbk
\\&=\lim_{\d\searrow0}\int_0^t\sum_{x\in\Xd} \frac{{\lbk\Ksptmrx - \Ksx\rbk}^2}{\dd} ds\\&=\int_{0}^{t}\int_{\Rd} {\lbk\KBsptmrx - \KBsx\rbk}^2 dxds
\\&=\lbk\int_{0}^{t}\KBBBsptmrz ds + \int_{0}^{t} \KBBBsz ds-2\int_{0}^{t} \KBBBsptmrsz ds\rbk
\\&=\lbk\int_{0}^{t}\tilde{\K}_{2s+2(t-r)} ds + \int_{0}^{t} \tilde{\K}_{2s} ds-2\int_{0}^{t} \tilde{\K}_{2s+(t-r)}ds\rbk
\\&=\lbk\int_{0}^{\tf{t-r}{2}}\tilde{\K}_{2s} ds-\int_{\tf{t-r}{2}}^{t-r}\tilde{\K}_{2s} ds-\int_{t}^{t+\tf{t-r}{2}}\tilde{\K}_{2s} ds
+\int_{t+\tf{t-r}{2}}^{2t-r}\tilde{\K}_{2s} ds\rbk
\end{split}
\eeq
for $1\le d\le 3$, where $\tilde{\K}_{w}$ is defined in terms of $\KBBBuvz$ by the relation
\beq\lbl{trans}
\bsp
&\tilde{\K}_{w}=\KBBBuvz\iff w=u+v\mbox{ and }(u,v)\mbox{ has one of the forms} \\&(u,v)=(a,a)\mbox{ or }(u,v)=(a+b,a)\mbox{ or }(u,v)=(a,a+b)\mbox{ for some}\ a,b\ge0.
\end{split}
\eeq

We observe that
\beq\lbl{Ktildebtrw}
\bsp
\tilde{\K}_{2u}=\KBBBuz&=4\int_{0}^{\infty}\int_{0}^{\infty}\proprtz\puzro\puzrt dr_{1} dr_{2}\\
&=4\int_{0}^{\infty}\int_{0}^{\infty}\lbk\int_{\Rd}\prox\prtx dx\rbk\puzro\puzrt dr_{1} dr_{2}\\
&=\int_{\Rd} \lbk\KBux\rbk^{2}dx=\df{C_{d}}{u^{d/4}};\quad1\le d\le3.
\end{split}
\eeq
The last assertion follows from the computation in \eqref{allcomp} (or see p.~531 in \cite{Abtpspde}).  It is clear then that $\tilde{\K}_{2u}$ is decreasing in $u$.  Thus, the sum of the last three terms of the \eqref{hlbd} is $\le0$.   This and \eqref{Ktildebtrw} give us \eqref{tmphl} for all $\delta\le\delta^{*}$, for some $\d^{*}>0$ and for some constant $C>0$, together with its corresponding BTBM density statement; and \lemref{4thQinequality} follows at once.\epf
The following spatial difference second moment inequality for the BTRW and BTBM densities captures their impressive spatial-regularizing effect on our solutions.
\begin{lem}\label{3rdQinequality}
For $1\le d\le3$, define the intervals
$$I_{d}=\bc (0,1];&d=1,\\
(0,1);&d=2,\\
(0,\frac12);&d=3.
\ec$$
For any given $d=1,2,3$ there exists a constant $C_{d}$ depending only on $d$ and $\alpha_{d}\in I_{d}$, and  a $\d^*>0$ such that for $\d\le\d^*$
\beq\lbl{sptlkernel}
\bc
\ds\int_0^t\intrd {\lbk\KBsx - \KBsxpz\rbk}^2 dx ds \leq C_{d} |z|^{2\alpha_{d}}t^{p_{d}(\alpha_{d})};&\alpha_{d}\in I_{d}, t>0,\\
\ds\int_0^t\sum_{x\in\Xd} {\lbk\Ksx - \Ksxpz\rbk}^2 ds \leq C_{d}\delta^{d} |z|^{2\alpha_{d}}t^{p_{d}(\alpha_{d})};& \alpha_{d}\in I_{d}, t>0,
\ec
\eeq
where $0<C_{d}<\infty$ and $0\le p_{d}(\alpha_{d})<1$ for every $\alpha_{d}\in I_{d}$ for $d=1,2,3$.
\end{lem}
\brm\lbl{spatlkerrem}
In the case $d=\alpha_{d}=1$ the power $p_{1}(\alpha_{1})=1/4$.  Also, the constants $C_{d}$'s are increasing in $\alpha_{d}$, with $0<C_{1}\le c<\infty$ for some absolute constant $c$ for all $0<\alpha_{1}\le1$; whereas $\lim_{\alpha_{2}\to1}C_{2}=+\infty=\lim_{\alpha_{3}\to1/2}C_{3}$.  Moreover, while their exact values are not needed, the following limits hold for the powers $p_{d}(\alpha_{d})$
\beq\lbl{powers}
\bc\ds
\lim_{\alpha_{d}\to1}p_{d}(\alpha_{d})=\frac14,\ \lim_{\alpha_{d}\to0}p_{d}(\alpha_{d})=\frac34; &d=1,\\ \\
\ds\lim_{\alpha_{d}\to1}p_{d}(\alpha_{d})=0,\ \lim_{\alpha_{d}\to0}p_{d}(\alpha_{d})=\frac12; &d=2,\\ \\
\ds\lim_{\alpha_{d}\to\tf12}p_{d}(\alpha_{d})=0,\ \lim_{\alpha_{d}\to0}p_{d}(\alpha_{d})=\frac14; &d=3.
\ec
\eeq
On any compact time interval $\T=[0,T]$, the inequality \eqref{sptlkernel} may---for any given value $\alpha_{d}\in I_{d}$---thus be rewritten as
\beq\lbl{sptlkernel2}
\bc
\ds\int_0^t\sum_{x\in\Xd} {\lbk\Ksx - \Ksxpz\rbk}^2 ds \leq
C_{d}\lpa\sup_{\alpha_{d}\in I_{d}}T^{p_{d}(\alpha_{d})}\rpa\delta^{d} |z|^{2\alpha_{d}}
\\ \qquad\qquad\qquad\qquad\qquad\qquad\qquad=\tilde{C}_{d}\delta^{d} |z|^{2\alpha_{d}};&\\
\ds\int_0^t\intrd {\lbk\KBsx - \KBsxpz\rbk}^2 dx ds \leq \tilde{C}_{d} |z|^{2\alpha_{d}};
\ec
\eeq
with $\tilde{C}_{d}<\infty$ in \eqref{sptlkernel2}.  We emphasize that, while we may take $\alpha_{1}=1$ in \eqref{sptlkernel} or \eqref{sptlkernel2}; we can neither take the limiting values $\alpha_{2}=1$ nor $\alpha_{3}=\frac12$ since $\lim_{\alpha_{2}\to1}C_{2}=+\infty=\lim_{\alpha_{3}\to1/2}C_{3}$.
\erm
\bpf
Starting with the $L^{2}$ estimate involving the spatial difference of the BTBM density in \eqref{sptlkernel}, we have
\beq\lbl{sptlebaktar}
\bsp
&\int_0^t\intrd {\lbk\KBsx - \KBsxpz\rbk}^2dx ds
\\&=4\int_0^t\lbk\int_0^\infty\int_0^\infty\intrd\prod_{i=1}^{2}{\lpa\prix - \prixpz\rpa}\pszri dxdr_{1}dr_{2}\rbk ds
\\&=4\int_{0}^{t}\int_0^\infty\int_{0}^{\infty}{\lpa2\proprtz - 2\proprtzz\rpa} \pszro\pszrt  dr_{1}dr_{2}ds
\\&=8\int_{0}^{t}\int_0^\infty\int_{0}^{\infty}\frac{1-e^{-\frac{|z|^{2}}{2(r_{1}+r_{2})}}}{\lbk2\pi (r_1+r_2)\rbk^{d/2}}\pszro\pszrt  dr_{1}dr_{2}ds
\\&=8\int_{0}^{t}\int_0^{\pi/2}\int_0^\infty\df{\lpa1-e^{-\frac{|z|^{2}}{2\rho (\sin(\theta)+\cos(\theta))}}\rpa e^{-\rho^2/2s}}{(2\pi)^{\tf{d+2}{2}}
\lbk\rho (\sin(\theta)+\cos(\theta))\rbk^{d/2}s}\rho d\rho d\theta ds
\\&\le C \int_{0}^{t}\int_0^\infty\df{\lpa1-e^{-\frac{|z|^{2}}{2\rho}}\rpa e^{-\rho^2/2s}}{\rho^{d/2}s}\rho d\rho  ds
\\&\le C \int_{0}^{t}\int_0^\infty\df{{{|z|^{2\alpha}}} e^{-\rho^2/2s}}{\rho^{\alpha+d/2}s}\rho d\rho ds
\\&\le\bc
 C_{1} |z|^{2\alpha}t^{p_{1}(\alpha)}; &d=1,\alpha\in(0,1], \\
 C_{2} |z|^{2\alpha}t^{p_{2}(\alpha)}; &d=2,\alpha\in(0,1), \\
 C_{3}|z|^{2\alpha}t^{p_{3}(\alpha)}; &d=3,\alpha\in(0,\tf{1}{2}),
\ec
\end{split}
\eeq
for some finite constants $C_{i}$, $i=1,2,3$, where $C_{2}$ and $C_{3}$ depend on $\alpha$\footnote{See \remref{spatlkerrem}.}, and where we have used the simple facts that $\min_{0\le\theta\le\pi/2}\lbk\sin(\theta)+\cos(\theta)\rbk=1$ and that $1-e^{-u}\le u^{\alpha}$ for $u\ge0$ and $0<\alpha\le1$.  This proves the $L^{2}$ estimate for the BTBM density in \eqref{sptlkernel}.  Then, an asymptotic argument similar to the one in the proofs of \lemref{2ndQinequality} and \lemref{4thQinequality} yields
\beq\lbl{ha}
\lim_{\d\searrow0}\int_0^t\sum_{x\in\Xd}\frac{{\lbk\Ksx - \Ksxpz\rbk}^2}{\dd} ds=\int_0^t\intrd {\lbk\KBsx - \KBsxpz\rbk}^2 dx ds,
\eeq
together with the desired BTRW density $L^{2}$ estimate in \eqref{sptlkernel} for all $\d\le\d^{*}$, for some $\d^{*}>0$, with possibly different constants.
\epf

\subsection{Spatio-temporal estimates for BTRW SIEs}\lbl{regandtight}
 In this subsection, and assuming only the less-than-Lipschitz conditions \eqref{cnd} on $a$---together with a temporary moment condition---we obtain spatial and temporal differences moments estimates that are crucial in obtaining the regularity
 of the BTRW SIE $\btrwsien$ for each fixed $n\in\Ns$ (see \eqref{star}), the tightness of the BTRW SIEs sequence $\lbr\btrwsien\rbr_{n\in\Ns}$, as well as the H\"older regularity for the limit.

 Fix $n\in\Ns$, and assume $\uttxn$ solves $\btrwsien$ in \eqref{btrwsie}.   Suppressing the dependence on $n$,
 let $M_q(t) = \sup_x \mathbb{E}|\tilde{U}^x(t)|^{2q}$, $q\ge1$.  Writing $\uttx$ in terms of
 its deterministic and random parts  $\uttx=\uttxD+\uttxR$, we observe that the deterministic part
 $\uttxD$ is smooth in time by \lemref{fodde}.  The next two lemmas give us estimates on the random part.
 \begin{lem}[Spatial differences] \lbl{sptdiff}
Assume that \eqref{cnd} holds and that $M_q(t)$ is bounded on any time interval\footnote{This is the aforementioned temporary  moment condition.  It is assumed here (in \lemref{sptdiff} and \lemref{tmpdiff} below) only to simplify the presentation and to get to the proof of \thmref{lip} as quickly as possible in \secref{pf1stmain}.  In \secref{regandtight2}, this moment condition is shown to automatically hold under \eqref{cnd}.} $\T=[0,T]$.
There exists a constant $C_{d}$ depending only on $q$, $\max_x |u_0(x)|$, the spatial dimension
$1\le d\le3$, $\alpha_{d}$, and $T$ such that
\[
\mathbb{E}\left| \tilde{U}_R^x(t) - \tilde{U}_R^y(t)\right|^{2q}
\le C_{d} |x-y|^{2q\alpha_{d}};\ \alpha_{d}\in I_{d},
\]
for all $x,y \in \Xd$, $t\in\T$, and $1\le d\le3;$ where $\alpha_{d}$ and $I_{d}$ are as in \lemref{3rdQinequality}.  I.e., in $d=1$, we may take $\alpha_{1}=1;$ in $d=2$ we may take any fixed $\alpha_{2}\in(0,1);$ and in $d=3$, $\alpha_{3}$ may be taken to be any fixed value in $(0,\frac12)$.
\end{lem}
\begin{proof}
Using Burkholder inequality, we have for any $(t,x,y)\in\T\times\X^{2d}$
\begin{equation}
\begin{split}
\E\lab \tilde{U}_R^x(t) - \tilde{U}_R^y(t)\rab^{2q}
\leq C\E{\lab\sum_{z\in\Xd}\int_0^t
{\lbk\Ktsxz- \Ktsyz\rbk}^2a^2
(\tilde{U}^{z}(s)) \df{ds}{\dd}\rab}^q
\end{split}
\label{spdiff1}
\end{equation}
For any fixed but
arbitrary point $(t,x,y)\in\T\times\X^{2d}$ let $\mu_t^{x,y}$ be the measure defined on
$[0,t]\times\Xd$ by
\beqs
\bsp
d\mu_t^{x,y}(s,z)&=\lbk\Ktsxz-\Ktsyz\rbk^2\df{ds}{\dd},
\end{split}
\eeqs
and let
$|\mu_t^{x,y}| = \mu_t^{x,y}([0,t] \times \Xd)$.    We see from \eqref{spdiff1}, Jensen's inequality
applied to the probability measure $\mu_t^{x,y}/\lab\mu_t^{x,y}\rab$,
the growth condition on $a$, the definition of $M_q(t)$, and elementary inequalities, that we have
\begin{equation}
\begin{split}
\mathbb{E}\left| \tilde{U}_R^x(t) - \tilde{U}_R^y(t)\right|^{2q}
&\le C\mathbb{E}\Big[\int_{[0,t] \times \Xd}
\left|a(\tilde{U}^{z}(s))\right|^{2q}\frac{d\mu_t^{x,y}(s,z)}{|\mu_t^{x,y}|}\Big]
{|\mu_t^{x,y}|}^q \\
&\le  C\Big[\int_{[0,t] \times \Xd}\left (1+M_q(s)\right)
\frac{d\mu_t^{x,y}(s,z)}{|\mu_t^{x,y}|}\Big] {|\mu_t^{x,y}|}^q
\end{split}
\label{spdiff2}
\end{equation}
Now, using the boundedness assumption on $M_q$ on $\mathbb{T}$ for $1\le d\le 3$,
we get
 \begin{equation*}
\begin{split}
\mathbb{E}\left| \tilde{U}_R^x(t) - \tilde{U}_R^y(t)\right|^{2q}\le C\lab\mu_t^{x,y}\rab^q
\le \lbk C_{d}t^{p_{d}(\alpha_{d})}\rbk^{q} |x-y|^{2q\alpha_{d}}\le \tilde{C}_{d}|x-y|^{2q\alpha_{d}};\hspace{-.1mm} \alpha_{d}\hspace{-.1mm}\in\hspace{-.1mm} I_{d},
\end{split}
\label{spdiff3}
\end{equation*}
where the last inequality follows from \lemref{3rdQinequality} and \eqref{sptlkernel2} in \remref{spatlkerrem}, and where the constant $\tilde{C}_{d}=\lbk C_{d}\sup_{\alpha_{d}\in I_{d}}T^{p_{d}(\alpha_{d})}\rbk^{q}<\infty$.
\end{proof}

\begin{lem}[Temporal differences]\label{tmpdiff}
Assume that \eqref{cnd} holds and that $M_q(t)$ is bounded on any time interval $\T=[0,T]$.
 There exists a constant
$C$ depending only on $q$, $\max_x|u_0(x)|$, the spatial dimension
$1\le d\le3$, and $T$ such that
\[\mathbb{E}\left| \tilde{U}_R^x(t) - \tilde{U}_R^x(r) \right|^{2q}
\leq C\lab t-r\rab^{\tf{(4-d)q}{4}},\] for all $x \in \Xd$, for all $t,r \in\mathbb{T}$, and for $1\le d\le 3$.
\end{lem}

\begin{proof}
Assume without loss of generality  that $r<t$.
Using Burkholder inequality, and using the change of variable $\rho=t-s$,
we have for $(r,t,x)\in\T^2\times\X^{d}$
\begin{equation}
\begin{split}
\E\lab \tilde{U}_R^x(t) - \tilde{U}_R^x(r)\rab^{2q}
&\leq C\E{\lab\sum_{z\in\Xd}\int_0^r{\lbk\Ktsxz- \Krsxz\rbk}^2a^2(\tilde{U}^{z}(s)) \df{ds}{\dd}\rab}^q
\\&+C\E{\lab\sum_{z\in\Xd}\int_0^{t-r}
{\lbk\Kroxz\rbk}^2a^2
(\tilde{U}^{z}(t-\rho)) \df{d\rho}{\dd}\rab}^q
\end{split}
\label{tmpdiff1}
\end{equation}
For a fixed point $(r,t,x)$,
let $\mu_{t,r}^x$ be the measure defined on $[0,r]\times\Xd$ by
\begin{equation*}
\bsp
d\mu_{t,r}^x(s,z)&=\lbk\Ktsxz-\Krsxz\rbk^2\df{ds}{\dd}
\end{split}
\end{equation*}
and let $|\mu_{t,r}^x| = \mu_{t,r}^x([0,r] \times \Xd)$.
Also, for a fixed $x\in\X^d$, let $\kappa^x$ be the measure defined
on $[0,t-r]\times\Xd$ by
\begin{equation*}
\bsp
d\kappa^x(\rho)&=\lbk\Kroxz\rbk^2\df{d\rho}{\dd}
\end{split}
\end{equation*}
and let $|\kappa^x| = \kappa^x([0,t-r] \times \Xd)$.
Then, arguing as in Lemma \ref{sptdiff}
above we get that
 \begin{equation*}
 \begin{split}
\mathbb{E}\left| \tilde{U}_R^x(t) - \tilde{U}_R^x(r)\right|^{2q}&\le
C\lpa\lab\mu_{t,r}^x\rab^q+\lab\kappa^x\rab^q\rpa\le C(t-r)^{\tf{(4-d)q}{4}},
\end{split}
\end{equation*}
for $1\le d\le3$, where the last inequality follows from \lemref{2ndQinequality} and
\lemref{4thQinequality}, completing the proof.
\end{proof}
\subsection{Proof of the first main Theorem}\lbl{pf1stmain}
Here, we prove \thmref{lip}.  We start first by recalling a useful elementary Gronwall-type lemma whose proof can be found in Walsh \cite{W}.
\blm\lbl{grnwlrate}
Let $\lbr g_n(t)\rbr_{n=0}^\infty$ be a sequence of positive functions such that $g_0$
is bounded on $\T=[0,T]$ and
$$g_n(t)\le C \intzt g_{n-1}(s) (t-s)^\alpha ds, \quad n=1,2,\ldots$$
for some constants $C>0$ and $\alpha>-1$.  Then, there exists a (possibly different)
constant $C>0$ and an integer $k>1$ such that for each $n\ge1$ and $t\in\T$
$$g_{n+mk}(t)\le C^m\intzt g_n(s)\df{t-s}{(m-1)!} ds;\quad m=1,2,\ldots.$$
\elm
We are now ready for our proof.
\begin{proof}[Proof of \thmref{lip}]
For the existence proof, we construct a solution iteratively.  So, given a space-time white noise $\sW$, on some $\OFFtP$,
define
\beq\lbl{itdef}
\bc
\zivtx=\ds\intrd\KBtxy \uny dy&\cr
\ds\npoivtx=\zivtx\ds+\intrdzt\KBtsxy a(\nivsy)\sW(ds\times dy)
\ec
\eeq
We will show that, for any $p\ge2$ and all $1\le d\le3$, the sequence $\lbr \nivtx\rbr_{n\ge1}$
converges in $L^p(\Omega)$ to a solution.  Let
$$D_{n,p}(t,x):=\E\lab\npoivtx-\nivtx\rab^p$$

$$D^*_{n,p}(t):=\sup_{x\in\Rd}\dnptx.$$
Starting with the case  $p>2$, we bound $D_{n,p}$ using Burkholder inequality,
the Lipschitz condition $(a)$ in \eqref{lcnd}, and then H\"older inequality with $0\le\epsilon\le1$
and $q=p/(p-2)$ to get
\beqs
\bsp
&\dnptx=\E\lab\intrdzt\KBtsxy\lbk a(\nivsy)-a(\nmoivsy)\rbk\sW(ds\times dy)\rab^{p}
\\&\le C\E\lab\intrdzt\lpa\KBtsxy\rpa^2\lbk \nivsy-\nmoivsy\rbk^2ds dy\rab^{p/2}
\\&\le C\lpa\intrdzt\lbk\KBtsxy\rbk^{2\epsilon q}dsdy\rpa^{p/2q}
\\&\times\intrdzt\lpa\KBtsxy\rpa^{(1-\epsilon)p}\dnmopsy ds dy
\end{split}
\eeqs
Take $\epsilon=(p-2)/p$ in the above ($2\epsilon q=(1-\epsilon)p=2$), take the supremum over the space variables,
and use the computation on p.~531 of \cite{Abtpspde}\footnote{See also \lemref{2ndQinequality} and its proof.}
to see that, for $1\le d\le3$ the above reduces to
\beq\lbl{mdiff}
\bsp
\mnpt\le C \lpa t^{\tf{4-d}{4}}\rpa^{\tf{p-2}{2}}\intzt\mnmops\lbk t-s\rbk^{-\tf{d}{4}} ds
\end{split}
\eeq
The case $p=2$ is simpler.  We apply Burkholder's inequality to $D_{n,2}$ and then take the space supremum
to get
\beq\lbl{mdiff1}
\bsp
D^*_{n,2}(t)\le C \intzt D^*_{n-1,2}(s)\lbk t-s\rbk^{-\tf{d}{4}} ds
\end{split}
\eeq
I.e., on any time interval $\T=[0,T]$, the integral multiplier on the r.h.s.~of \eqref{mdiff} is bounded; and
if $D^*_{n-1,p}$ is bounded on $\T$ then so is $D^*_{n,p}$, for every $p\ge2$.  Now,
$$D^*_{0,p}(t)\le C\sup_{x\in\Rd}\E\lab\intrdzt\lbk\KBtsxy\rbk^{2}a^2\lpa\zivsy\rpa dsdy\rab^{\tf p 2}$$
Since $\un$ is bounded and deterministic, then so are  $U^{(0)}$ and $a(U^{(0)})$.  The latter assertion follows
from the growth condition on $a$ in \eqref{lcnd}.  Thus, by the computation on p.~531 of \cite{Abtpspde},
$D^*_{0,p}$ is bounded on $\T$  for $1\le d\le3$ and so are all the
$D^*_{n,p}$.  \lemref{grnwlrate} now implies that for each $1\le d\le 3$, the series
$\sum_{m=0}^\infty \lbk D^*_{n+mk,p}(t)\rbk^{1/p}$ converges uniformly on compacts for
each $n$, which in turn implies that the series $\sum_{n=0}^\infty\lbk D^*_{n,p}(t)\rbk^{1/p}$ converges uniformly
on compacts.  Thus $U^{(n)}$ converges in $L^p(\Omega)$ for $p\ge2$, uniformly on $\T\times\Rd$ for $1\le d\le3$.  Let
$U(t,x):=\lim_{n\to\infty}\nivtx$.
It is easy to see that $U$  satisfies \eqref{btpsol},
and hence solves the BTBM SIE $\btsie$.  This follows from \eqref{itdef} since
the Lipschitz condition in \eqref{lcnd} gives
\beqs
\bsp
\E\lab a(U(t,x))-a(\nivtx)\rab^2\le C\E\lab U(t,x)-\nivtx\rab^2\to0\quad\mbox{as }n\to\infty
\end{split}
\eeqs
uniformly on $\T\times\Rd$.  Therefore, the stochastic integral term in \eqref{itdef} converges to
the same term with $U^{(n)}$ replaced with the limiting $U$---i.e., it converges to
the corresponding term in $\btsie$---as $n\to\infty$, for
\beqs
\bsp
&\E\lbk\intrdzt\KBtsxy\lpa a(U(s,y))- a(\nivsy)\rpa\sW(ds\times dy)\rbk^2
\\&\le  C \intrdzt\lbk\KBtsxy\rbk^2 \E\lbk U(s,y)-\nivsy\rbk^2 ds dy
\longrightarrow 0
\end{split}
\eeqs
as $n\to\infty$.  It follows that $U$ satisfies the BTBM SIE $\btsie$. Also, the solution
is strong since the $U^{(n)}$ are constructed for a given white noise $\sW$, and the
limit $U$ satisfies \eqref{btpsol} with respect to that same $\sW$.
Clearly $U$ is $L^p(\Omega)$ bounded on $\T\times\Rd$, $1\le d\le3$, for any $p\ge2$ and for any $T>0$.

To show uniqueness let $1\le d\le 3$, let $T>0$ be fixed but arbitrary, and let $U_1$ and $U_2$ be two solutions to the BTP
SIE \eqref{btpsol} that are $L^2(\Omega)$-bounded on $\T\times\Rd$.
Fix an arbitrary $(t,x)\in\Rp\times\Rd$.
Let $D(t,x)=U_2(t,x)-U_1(t,x)$, $L_2(t,x)=\E D^2(t,x)$,
and $L^*_2(t)=\sup_{x\in\Rd}L_2(t,x)$ (which is bounded on $\T$ by hypothesis).
Then, using \eqref{btpsol}, the Lipschitz condition in \eqref{lcnd},
and taking the supremum over the space variable and using the computation on p.~531 of \cite{Abtpspde}
we have
\beq\lbl{un1}
\bsp
L_2(t,x)&=\intrdzt \E\lbk a(U_2(s,y))- a(U_1(s,y))\rbk^2\lbk\KBtsxy\rbk^2dsdy
\\&\le C\intrdzt L_2(s,y)\lbk\KBtsxy\rbk^2dsdy
\\&\le C\intzt L^*_2(s)\intrd\lbk\KBtsxy\rbk^2dyds
\le C\intzt\df{L^*_2(s)}{(t-s)^{\tf d4}}ds
\end{split}
\eeq
Iterating and interchanging the order of integration we get
 \beq\lbl{un2}
\bsp
L_2(t,x)&\leq  C\lbr \int_0^t{L^*_2(r)}\lpa\int_r^t\frac{ds}{{(t - s)^{d/4}}{(s-r)^{d/4}}}
 \rpa dr\rbr\\
 &\le C\lpa\ds\int_0^t L^*_2(s)ds \rpa
\end{split}
\eeq
for any $1\le d\le 3$.  Hence,
\beq\lbl{un3}
\bsp
L^*_2(t)\le C\lpa\ds\int_0^t L^*_2(s)ds \rpa
\end{split}
\eeq
for every $t\ge0$.  An easy application of Gronwall's lemma gives that $L^*_2\equiv0$.  So for every $(t,x)\in\Rp\times\Rd$ and $1\le d\le3$ we have $U_1(t,x)=U_2(t,x)$ with probability one.   The indistinguishability of $U_{1}$ from $U_{2}$, and hence pathwise uniqueness, follows immediately from their H\"older regularity, which we now turn to.

We have just shown that, under the Lipschitz conditions \eqref{lcnd}, our BTBM SIE in \eqref{btpsol} has an $L^{p}(\Omega)$-bounded solution $U(t,x)$ on $\T\times\Rd$ for any $T>0$ and any $p\ge2$.  Equivalently, $M_q(t)= \sup_x \mathbb{E}|U(t,x)|^{2q}$, $q\ge1$, is bounded on any time interval $\T$.   Recalling that the deterministic part\footnote{Of course, the deterministic part of $\btsie$ is, as discussed before, the integral $\intrd\KBtxy \uny dy$; and the random part is  $\intrdzt\KBtsxy a(\usy)\sW(ds\times dy)$.} of $U$ is a $\mathrm{C}^{1,4}(\Rp,\Rd)$ function, we can then repeat the same arguments in the proofs of \lemref{sptdiff} and \lemref{tmpdiff} above, on the random part of $U$, with obvious modifications---replace $\tilde{U}$ with $U$, sums over $\Xd$ with integrals over $\Rd$ (with  $ds dz$ instead of $ds/\dd$ or $d\rho/\dd$), and $\BTRW$ with $\BTBM$---and use the BTBM statements in \lemref{2ndQinequality}, \lemref{4thQinequality}, and \lemref{3rdQinequality}) to get the same estimates on the spatial and temporal differences of $U$, with possibly different constants:
\beq\lbl{spttempdirect}
\bc
\E\lab\utx-\uty\rab^{2q}\le C_{d} |x-y|^{2q\alpha_{d}};\ \alpha_{d}\in I_{d},\\
\E\lab\utx-\urx\rab^{2q}\le C_{d}\lab t-r\rab^{\tf{(4-d)q}{4}},
\ec
\eeq
for $d=1,2,3$.  This in turn straightforwardly leads to the desired local H\"older regularity for the direct solution of $\btsie$, $U$, as follows: we let $q_n=n+d$ for $n\in\{0,1,\ldots\}$ and
let $n=m+d$ for $m=\{0,1,\ldots\}$, we then have from \eqref{spttempdirect}
that
\beq\lbl{directbtpsiehldr}
\bc
\E\lab\utx-\uty\rab^{2n+2d}\le C_{d}\lab x-y\rab^{(2n+2d)\alpha_{d}},\\
\E\lab\utx-\urx\rab^{2m+4d}\le C\lab t-r\rab^{\tf{(4-d)(m+2d)}{4}}.
\ec
\eeq
for $1\le d\le 3$.  Thus as in Theorem 2.8 p.~53 and Problem 2.9 p.~55 in \cite{KS} we get that the spatial H\"older exponent is
$\gamma_s\in\lpa 0,\tf{2(n+d)\alpha_{d}-d}{2n+2d}\rpa$ and the temporal exponent is $\gamma_t\in\lpa0,\tf{m\lpa1-d/4\rpa+d(1-d/2)}{2m+4d}\rpa$
$\forall m,n$.  Taking the limits as $m,n\to\infty$, we get
$\gamma_t\in\lpa0,\tf{4-d}{8}\rpa$ and $\gamma_s\in\lpa0,\alpha_{d}\rpa$, for $1\le d\le 3$.
The proof is complete.
\end{proof}
\subsection{Regularity and tightness without the Lipschitz condition}\lbl{regandtight2}
As we mentioned in \secref{regandtight}, the finiteness assumption of $M_{q}(t)$ on $\T$ in \lemref{sptdiff} and \lemref{tmpdiff} is for convenience only.  We now proceed to show how to remove that assumption by showing it automatically holds under the weaker conditions \eqref{cnd}. It is easily seen that if $a$ is bounded then,  for all spatial dimensions $1\le d\le 3$, $M_q$ is bounded on any compact time interval $\T=[0,T]$ (see \remref{bddrem} below).  The following Proposition gives an exponential upper bound on the growth of $M_q$ in time in all $1\le d\le3$ under the conditions in \eqref{cnd}.
\begin{prop}[Exponential bound for $M_q$] \label{Expbd}
Assume that $\uttx$ is a solution of the BTRW SIE $\btrwsien$ on $\T\times\Xd$, and assume that conditions \eqref{cnd} are satisfied.  There exists a constant $C$ depending only on $q$, $\max_x|u_0(x)|$, the dimension $d$, and $T$ such that
\begin{equation*}
M_q(t) \leq  C\lpa1+ \ds\int_0^t M_q(s)ds \rpa; \quad\forall\,0\le t\le T, \ q\ge1,\mbox{ and }1\le d\le 3,
\end{equation*}
and hence
\begin{equation*}
M_q(t) \leq C\exp{\{Ct\}}; \quad \forall\,0\le t\le T, \ q\ge1,\mbox{ and }1\le d\le 3.
\end{equation*}
   In particular, $M_q$ is bounded on $\mathbb{T}$ for all $q\ge1$
and $1\le d\le 3$.
 \end{prop}

 The proof of Proposition \ref{Expbd} proceeds via the following lemma and its
 corollary.

\begin{lem} Under the same assumptions as in \propref{Expbd}
there exists a constant $C$ depending only on $q\ge1$, $\max_x |u_0(x)|$, the dimension $d$,
and $T$ such that
\begin{equation*}
M_q(t) \leq\bc C \lpa1 + \ds{\int_0^t\df{M_q(s)}{{(t - s)}^{d/4}} ds}\rpa;
   &0< t\le T,\ q\ge1,\mbox{ and } 1\le d\le 3,\\
  C;&t=0, \ q\ge1,\mbox{ and } 1\le d\le 3.
  \ec
\end{equation*}
\label{1stboundonUtilde}
\end{lem}

\begin{proof}
Fix $q\ge1$, let $\ds\tilde{U}_D^{x}(t)\overset{\triangle}{=} \sum_{y\in\Xd} \Ktxy u_0(y)$
(the deterministic part of $\tilde{U}$).  Then,  for any $(t,x)\in\T\times\Xd$,
we apply Burkholder inequality to the random term  $\tilde{U}_R^x(t)$ to get
\begin{equation}\label{afterBurkholder}
\begin{split}
\E\lab\tilde{U}^x(t)\rab^{2q}&= \E\lab \sum_{y\in {\Xd}}\int_0^t \Ktsxy
\frac{a(\tilde{U}^y(s))}{\sdd} dW^{y}(s)+ \tilde{U}_D^{x}(t) \rab^{2q} \\
\\&\leq C\lpa\E\lab\sum_{y\in {\Xd}} \int_0^t {\lpa\Ktsxy\rpa}^2\frac{a^2(\tilde{U}^y(s))}{\dd}ds\rab^q + \lab\tilde{U}_D^{x}(t)\rab^{2q}\rpa.
\end{split}
\end{equation}
Now, for a fixed point $(t,x)\in\T\times\Xd$ let $\mu_t^x$ be the measure on $[0,t] \times\Xd$ defined by
$d\mu_t^x(s,y)=\lbk{\lpa\Ktsxy\rpa}^2/\dd\rbk ds$, and let $|\mu_t^x| = \mu_t^x([0,t] \times \Xd)$.  Then, we can rewrite \eqref{afterBurkholder} as
 \begin{equation}\label{afterBurkholder2}
\begin{split}
 \E\lab\tilde{U}^x(t)\rab^{2q}
\leq C\lpa\E\left|\int_{[0,t] \times \Xd}a^2(\tilde{U}^y(s))
\frac{d\mu_t^x(s,y)}{|\mu_t^x|}\right|^q {|\mu_t^x|}^q +|\tilde{U}_D^{x}(t)|^{2q}
\rpa.
\end{split}
\end{equation}
Observing that $\mu_t^x/|\mu_t^x|$ is a probability measure, we apply Jensen's inequality,
the growth condition on $a$ in \eqref{cnd}, and other elementary inequalities
to \eqref{afterBurkholder2} to obtain
\begin{align*}
&\E\lab\tilde{U}^x(t)\rab^{2q}
\leq C\lpa\E\lbk\int_{[0,t] \times \Xd}\lab a(\tilde{U}^y(s))
\rab^{2q}\frac{d\mu_t^x(s,y)}{|\mu_t^x|}\rbk
{|\mu_t^x|}^q + \lab\tilde{U}_D^{x}(t)\rab^{2q}\rpa\\
& \leq C\lbk\int_{[0,t] \times \Xd}
\lpa1 + \E\lab\tilde{U}^y(s)\rab^{2q}\rpa{d\mu_t^x(s,y)}\rbk
{|\mu_t^x|}^{q-1} + C \lab\tilde{U}_D^{x}(t)\rab^{2q} \\
&= C\lpa\lbk\sum_{y\in {\Xd}}\int_0^t\frac{\lpa\Ktsxy\rpa^2}
{\dd}\lpa1 + \E\lab\tilde{U}^y(s)\rab^{2q}\rpa ds\rbk
{|\mu_t^x|}^{q-1}+\lab\tilde{U}_D^{x}(t)\rab^{2q} \rpa
\end{align*}
Using \lemref{2ndQinequality} we see that  $|\mu_t^x|$ is
uniformly bounded for $t\leq T$ and $1\le d\le 3$.
So, using the boundedness of $u_0$, and hence of $\tilde{U}_D^{x}(t)$
by the simple fact that $\sum_{y\in\Xd}\Ktxy=1$,
Lemma \ref{2ndQinequality} and the definition of $M_q(s)$, we get
\begin{align*}
\E\lab\tilde{U}^x(t)\rab^{2q}
&\leq C\lpa1+\sum_{y\in {\Xd}}
\int_0^t\frac{\lpa\Ktsxy\rpa^2}{\dd}M_q(s) ds\rpa\\
&\overset{R_1}{\leq} C \lpa1 + \int_0^t\frac{M_q(s)}{{(t - s)}^{d/4}} ds
\rpa.
\end{align*}
Here, $R_1$ holds for $1\le d\le 3$.   This implies that
\begin{equation*}
M_q(t) \leq C \lpa1 + \int_0^t\frac{M_q(s)}{{(t - s)}^{d/4}} ds\rpa.
\end{equation*}
Of course, $M_q(0)=\sup_{x}\lab\unx\rab^{2q}\le C$, by the boundedness and nonrandomness
assumptions on $\unx$ in \eqref{cnd}. The proof is complete.
\end{proof}
\brm\lbl{bddrem}
It is clear that for a bounded $a$, $M_q$ is locally bounded in time.  This follows immediately from \lemref{2ndQinequality}
along with \eqref{afterBurkholder2} above.
\erm
\begin{cor}
Under the same assumptions as those in \propref{Expbd}
there exists a constant $C$ depending only on $q$, $\max_x |u_0(x)|$, the dimension $d$, and $T$
such that
\begin{equation*}
M_q(t) \leq C\lpa1+ \ds\int_0^t M_q(s)ds \rpa, 0\le t\le T, \ q\ge1,\mbox{ and }1\le d\le 3;
\end{equation*}
and hence
\begin{equation*}
M_q(t) \leq C\exp{\{Ct\}}; \quad \forall\,0\le t\le T, \ q\ge1,\mbox{ and }1\le d\le 3.
\end{equation*}
\label{befG}
\end{cor}

\begin{proof}
Iterating the bound in Lemma \ref{1stboundonUtilde} once, and changing the
order of integration, we obtain
\begin{equation}
\begin{split}
 M_q(t) &\leq  C\lbr 1 + C\lbk\int_0^t\frac{ds}{(t-s)^{d/4}}+
 \int_0^t{M_q(r)}\lpa\int_r^t\frac{ds}{{(t - s)^{d/4}}{(s-r)^{d/4}}}
 \rpa dr\rbk\rbr\\
 &\le C\lpa1+ \ds\int_0^t M_q(s)ds \rpa
 \end{split} \label{mom1}
\end{equation}
for $1\le d\le3$.  The proof of the last statement is a straightforward application of Gronwall's
lemma to \eqref{mom1}.  This finishes the proof of \coref{befG} and thus of \propref{Expbd}.
\end{proof}

The regularity, tightness, and weak limit conclusions for the BTRW SIEs now follow.

\blm[Regularity and tightness]\lbl{regtight}
Assume that the conditions \eqref{cnd} hold, and that
$\lbr\uttxn\rbr_{n\in\Ns}$ is a sequence of spatially-linearly-interpolated solutions to
the BTRW SIEs $\lbr\btrwsien\rbr_{n\in\Ns}$ in \eqref{btrwsie}.   Then
\begin{enumerate}\renewcommand{\labelenumi}{$($\alph{enumi}$)$}
\item For every $n$, $\uttxn$ is continuous on $\Rp\times\Rd$.  Moreover, with probability one,
the continuous map $(t,x)\mapsto\uttxn$ is locally
$\gamma_t$-H\"older continuous in time with $\gamma_t\in\lpa0,\tf{4-d}{8}\rpa$ for $1\le d\le3$.
\item  There is a BTRW SIE weak limit solution to $\btsie$,  call it $U$, such that
$U(t,x)$ is $L^p(\Omega,\P)$-bounded on $\T\times\Rd$ for every $p\ge2$ and $U\in\mathrm{H}^{\tf{4-d}{8}^{-},\alpha^{-}_{d}}(\T\times\Rd;\R)$ for every $1\le d\le3$ and $\alpha_{d}\in I_{d}$, where $\alpha_{d}$ and $I_{d}$ are as in \lemref{3rdQinequality}.
\end{enumerate}
\elm
\brm
Of course in part (a) above, even without linear interpolation in space, $\uttx$ is locally H\"older continuous in time with
H\"older exponent $\gamma_t\in\lpa0,\tf{4-d}{8}\rpa$ for $1\le d\le 3$.
\erm
\bpf  For each $n$, let $\uttxn=\uttxnD+\uttxnR$ be the decomposition of $\uttxn$ in \eqref{btrwsie}
into its deterministic and random parts, respectively.
\begin{enumerate}\renewcommand{\labelenumi}{$($\alph{enumi}$)$}
\item By \lemref{fodde}, $\uttxnD$ is clearly smooth in time;
so it is enough to consider the random term $\uttxnR$.
We let $q_m=m+2$ for $m\in\{0,1,\ldots\}$,
we then have from \lemref{tmpdiff} that
\beq\lbl{siedishldrt}
\E\lab\tilde{U}_R^x(t) - \tilde{U}_R^x(r) \rab^{4+2m}\le C\lab t-r\rab^{\tf{(4-d)(m+2)}{4}}.
\eeq
for $1\le d\le 3$.  Thus as in Theorem 2.8 p.~53 \cite{KS} we get that the temporal H\"older exponent
$\gamma_t\in\lpa0,\tf{(8+(4-d)m-6d)/4}{2m+4}\rpa$
for every $m$.  Taking the limit as $m\to\infty$, we get
$\gamma_t\in\lpa0,\tf{4-d}{8}\rpa$ for $1\le d\le 3$.
\item   By \lemref{btrwasbtp} it follows that $\uttxnD$ converges pointwise
to the deterministic part of $\btsie$ in \eqref{btpsol}; i.e.,
\beq\lbl{det2det}
\lim_{n\to\infty}\uttxnD=\intrd\KBtxy \uny dy.
\eeq
We also conclude from \lemref{sptdiff} and \lemref{tmpdiff} that the random part sequence $\lbr\uttxnR\rbr_{n\in\Ns}$ is tight on $\mathrm{C}(\T\times\Rd)$  for $1\le d\le3$.  Thus, there exists a weakly convergent
subsequence $\lbr\utnk\rbr_{k\in\N}$ and hence a BTRW SIE weak limit solution $U$ to $\btsie$.
Then, following Skorokhod, we can construct processes\footnote{As usual, $\eqd$ means equal in law.} $\yk\eqd\utnk$ on some probability space $\Skspace$ such
that with probability $1$, as $k\to\infty$, $\ytxk$ converges to a random field $\ytx$ uniformly on compact subsets of
$\T\times\Rd$ for $1\le d\le3$.  Now, for the BTRW SIEs limit regularity assertions, clearly the deterministic term on the right hand side of
\eqref{det2det} is $\mathrm{C}^{1,4}$ and bounded as in \cite{Abtp1,Abtp2}, so we use \propref{Expbd}, \lemref{sptdiff}, and \lemref{tmpdiff}
to obtain the regularity results for the random part.  We provide the steps here for completeness.
First, $Y_k\overset{d}{=} \Tilde{U}_{n_k}$ and so \propref{Expbd} gives us, for each $p\ge2$:
\beq\lbl{reg111}
\E\lab\ytxk\rab^{p}=\E\lab\uttxnk\rab^{p}\le C<\infty;\forall(t,x,k)\in\T\times\Rd\times\N,
1\le d\le 3,
\eeq
for some constant $C$ that is independent of $k, t, x$ but that depends on the dimension $d$.
It follows that, for each $(t, x)\in\T\times\Rd$ the sequence $\lbr|Y_k(t, x)|^p\rbr_k$ is uniformly
integrable for each $p\ge2$ and each $1\le d\le3$. Thus,
\beq\lbl{spdemmntsbd}
\E\lab\utx\rab^{p}=\E\lab\ytx\rab^{p}=\lim_{k\to\infty}\E\lab\ytxk\rab^{p}\le C<\infty;\ \forall(t,x)\in\T\times\Rd,
\eeq
for all $1\le d\le 3$ and $p\ge2$.   Equation \eqref{spdemmntsbd} establishes the $L^p$ boundedness assertion.  In addition, for $q\ge1$ and $1\le d\le 3$ we have by \propref{Expbd}
\beq\lbl{ui2}
\bsp
&\E \lab Y_k(t, x) - Y_k(t, y)\rab^{2q} + \E\lab Y_k(t, x) - Y_k(r, x)\rab^{2q}
\\&\le C \lbk\E\lab Y_k(t, x)\rab^{2q} + \E\lab Y_k(t, y)\rab^{2q} + \E\lab Y_k(r, x)\rab^{2q}\rbk
\\&\le C; \ \forall(k, r, t, x, y)\in\N\times\T^2\times \R^2.
\end{split}
\eeq
So, for each $(r, t, x, y)\in\T^2\times\R^2$, the sequences $\lbr\lab Y_k(t, x) - Y_k(t, y)\rab^{2q}\rbr_k$ and
$\lbr\lab Y_k(t, x) - Y_k(r, x)\rab^{2q}\rbr_k$
are uniformly integrable, for each $q\ge1$. Therefore, using
\lemref{sptdiff} and \lemref{tmpdiff}, we obtain
\beq\lbl{spdehldr1}
\bc
\E\lab\utx-\uty\rab^{2q}=\E\lab\ytx-\yty\rab^{2q}&\cr
=\ds\lim_{k\to\infty}\E\lab\ytxk-\ytyk\rab^{2q}\le C_{d} |x-y|^{2q\alpha_{d}};\ \alpha_{d}\in I_{d},\\
\E\lab\utx-\urx\rab^{2q}=\E\lab\ytx-\yrx\rab^{2q}&\cr
=\ds\lim_{k\to\infty}\E\lab\ytxk-\yrxk\rab^{2q}\le C\lab t-r\rab^{\tf{(4-d)q}{4}},
\ec
\eeq
for $d=1,2,3$.  The local H\"older regularity is then obtained using exactly the same steps following \eqref{spttempdirect}.
\end{enumerate}
The proof is complete
\epf
\subsection{From K-martingale problems to truncated BTRW SIEs}\lbl{kmart2trnctd}
We now establish \thmref{kmarteqaux}.
\begin{proof}[Proof of \thmref{kmarteqaux}]
Assume that $\lpa\lbr\xtaunlf \rbr_{\tau\ge0},\lbr\wtyn\rbr_{y\in\Xnldf}\rpa$ is a solution to
$\lbr\btrwsienlaux\rbr_{\tau\ge0}$ in \eqref{aux} on a filtered probability space $\OFFtP$; then clearly $\xztauxnlf=\unx$ for $x\in\Xnd$
and $\xttauxnlf=\uttxnD$ for $x\in\Xnd\sm\Xnldf$, for every $\tau\ge0$ and $t\in[0,\tau]$.  Fix any arbitrary $\tau_1,\tau_2>0$  and $x_1,x_2\in\Xnldf$.
By It\^o's formula we have
\beq\lbl{itokmart}
\bsp
&f(Z_{n,l}^{x_{1,2},\tau_{1,2}}(t))-f(u_0^{x_1,x_2})-\ds\int_0^t \lpa\sAKtauot f\rpa(s,x_1,x_2,Z_{n,l}^{\tau_{1,2}})ds
\\&=\sum_{i=1}^2\sum_{y\in\Xnldf}\int_0^t\partial_if(Z_{n,l}^{x_{1,2},\tau_{1,2}}(s))\kappa^{x_i,y}_{\h,s,\tau_i}\lpa\xstauiynlf\rpa{d\wsyn}\in\locm
\end{split}
\eeq
for $t\in[0,\tau_1\wedge\tau_2]$ and \eqref{kmart} is satisfied.

Conversely, if a family of adapted $$\lbr\xtaunlf\rbr_{\tau\ge0}=\lbr\lbr\xttauxnlf;0\le t\le\tau, x\in\Xnd\rbr\rbr_{\tau\ge0}$$ defined on a
probability space $\OFFtP$ satisfies \eqref{kmart}; then fixing any two $\tau_{1},\tau_{2}>0,$ letting $\B_{0,R}=\lbr u=(u_1,u_2)\in\Rs; \lab u\rab\le R\rbr$
and choosing $f_1,f_2,f\in\C_b^2(\R;\R)$ such that $f_i(u)=u_i$ for $i=1,2$ and $f(u)=u_1u_2$ whenever $u\in\B_{0,R}$ we see that
\beq\lbl{2mart}
\bc
\mtautxi&:=\xttauxinlf-\uttxinD\in\locm; i=1,2,\\
\ntautxot&:=\ds\prod_{i=1}^2\xttauxinlf-\prod_{i=1}^2\unxi
\\&\ds-\sum_{\substack{1\le i,j\le2\\i\neq j}}\int_0^t\xstauxinlf d\utsxjnD
\\&\ds-\int_0^t\sum_{y\in\Xnldf}\Upsilon^{x_{1,2},y}_{\h,t,\tau_{1,2}}\lpa Z_{n,l}^{y,\tau_{1,2}}(s)\rpa ds\in\locm
\ec
\eeq
for all $t\in[0,\tau_{1}\wedge\tau_{2}]$ and $x_1,x_2\in\Xnldf$, where we have used the notation $Z_{n,l}^{\tau_{1,2}}$ for the two-dimensional process defined in \eqref{Z}.  We then have that
\beq\lbl{qmart}
\bsp
&\prod_{i=1}^2\mtautxi-\sum_{y\in\Xnldf}\int_0^t\Upsilon^{x_{1,2},y}_{\h,t,\tau_{1,2}}\lpa Z_{n,l}^{y,\tau_{1,2}}(s)\rpa ds
\\&=\ T_1^{x_{1,2},\tau_{1,2}}(t)+T_2^{x_{1,2},\tau_{1,2}}(t)
\end{split}
 \eeq
where
\beq\lbl{m1}
T_1^{x_{1,2},\tau_{1,2}}(t):=\ntautxot-\sum_{\substack{1\le i,j\le2\\i\neq j}}\unxi\mtautxj\in\locm
\eeq
and
\beq\lbl{m2}
\bsp
T_2^{x_{1,2},\tau_{1,2}}(t):&=\sum_{\substack{1\le i,j\le2\\i\neq j}}\int_0^t\lbk\xstauxinlf-\xttauxinlf\rbk d\utsxjnD
+\prod_{i=1}^2\lbk\uttxinD-\unxi\rbk
\\&=\sum_{\substack{1\le i,j\le2\\i\neq j}}\int_0^t\lbk\unxi-\utuxinD\rbk d\mtauuxj\in\locm.
\end{split}
\eeq
Thus,
\beq\lbl{qvkmart}
\lqv\mtauxo,\mtauxt\rqv_t= \sum_{y\in\Xnldf}\int_0^t\Upsilon^{x_{1,2},y}_{\h,t,\tau_{1,2}}\lpa Z_{n,l}^{y,\tau_{1,2}}(s)\rpa  ds.
\eeq
Equations \eqref{2mart} and \eqref{qvkmart} imply that there exists a set of $r=\#\lbr y;y\in\Xnldf\rbr$ independent
Brownian motions $\lbr\wtyn;t\in\Rp\rbr_{y\in\Xnldf}$ on an extension $\OFFtPt$ such that
\beq\lbl{repres}
\mtautx=\sum_{y\in\Xnldf}\int_0^t\kappa^{x,y}_{\h,s,\tau}\lpa\xstauynlf\rpa{d\wsyn},\ \forall(x,\tau,t)\in\Xnldf\times\Rp\times[0,\tau].
\eeq
In fact, fixing any $\tau>0$ and labeling the $\lbr x;x\in\Xnldf\rbr$ as $\lbr x_1,\ldots,x_r\rbr$,
the restriction of the desired family of BMs $\lbr\wtyn;t\in\Rp\rbr_{y\in\Xnldf}$ to the time interval $[0,\tau]$ is obtained from the matrix equation (written in differential form)
\beq\lbl{MatrixBMs}
\bsp
\begin{bmatrix}d{W}_n^{x_1}(t)\\\vdots\\d{W}_n^{x_r}(t)\end{bmatrix}&=\begin{bmatrix}\kappa^{x_1,x_1}_{\h,t,\tau}\lpa\xttauxonlf\rpa&\cdots&\kappa^{x_1,x_r}_{\h,t,\tau}\lpa\xttauxrnlf\rpa\\
\vdots&\vdots&\vdots\\
\kappa^{x_r,x_1}_{\h,t,\tau}\lpa\xttauxonlf\rpa&\cdots&\kappa^{x_r,x_r}_{\h,t,\tau}\lpa\xttauxrnlf\rpa\end{bmatrix}^{-1}
\begin{bmatrix}dM^{\tau,x_{1}}(t)\\\vdots\\dM^{\tau,x_{r}}(t)\end{bmatrix}
\\&:=\begin{bmatrix}d\tilde{W}_n^{\tau,x_1}(t)\\\vdots\\d\tilde{W}_n^{\tau,x_r}(t)\end{bmatrix}
\end{split}
\eeq
whenever the middle inverse kernel-diffusion coefficient $r\times r$ matrix, denoted by $A^{-1}$, exists (the determinant  ${\mathrm{det}}(A)\neq0$) almost surely. If this fails we can proceed similar to the standard finite dimensional SDE case cf. Ikeda and Watanabe \cite{IW} or Doob \cite{Do}.  It is now straightforward to verify that, for any $\tau_{1},\tau_{2}>0$ and any $t\in[0,\tau_{1}\wedge\tau_{2}]$, the two families of BMs $\lbr \lbr\tilde{W}_n^{\tau_{k},y}\rbr_{y\in\Xnldf};k=1,2\rbr$ satisfy
\beq\lbl{covbms}
\langle \tilde{W}_{n}^{\tau_{1},x_{i}}(\cdot), \tilde{W}_{n}^{\tau_{2},x_{j}}(\cdot)\rangle_{t}
=\bc
t,&1\le i=j\le r\\
0,&1\le i\neq j\le r
\ec
\eeq
almost surely, whether or not $\tau_{1}=\tau_{2}$.  I.e., we get one family of independent BMs $\lbr\wtyn\rbr_{y\in\Xnldf}$, such that $\lpa\lbr\xtaunlf\rbr_{\tau\ge0},\lbr W_n^y\rbr_{y\in\Xnldf}\rpa$ is a solution pair to  $\lbr\btrwsienlaux\rbr_{\tau\ge0}$ in \eqref{aux} on $\OFFtPt$.
Hence,  on the probability space $\OFFtPt$, the pair $\lpa\lbr\xttxnlf\rbr_{t\ge0,x\in\Xnd},\lbr\wtyn\rbr_{y\in\Xnldf}\rpa$ solves the $l$-truncated BTRW SIE $\btrwsienl$ in \eqref{trnctdsie} on $\Rp\times\Xnd$.
\end{proof}
\subsection{Completing the proof of the second main result}
\lbl{lipvnlip}
We now complete the proof of \thmref{mainthm2}.  In \secref{regandtight} and \secref{regandtight2} we assumed the existence of a BTRW SIE solution and
we obtained regularity and tightness for the sequence of lattice SIEs $\lbr\btrwsien\rbr_{n\in\Ns}$.  This, in turn,
implied the existence and regularity for a BTRW SIE limit solution to our $\btsie$ in \eqref{btpsol}.   To complete the existence of the desired double limit solution\footnote{The type of our lattice limit  solution to $\btsie$ in \eqref{btpsol} depends on the conditions: under the Lipschitz conditions \eqref{lcnd} we get a direct solution to the lattice SIE $\btrwsien$ for every $n$ and a direct BTRW SIE limit solution to $\btsie$ (see \thmref{latlimlip} in Appendix \ref{appB}); whereas under the non-Lipschitz conditions in \eqref{cnd} we obtain a limit BTRW SIE solution, thanks to our K-martingale approach, and a BTRW SIEs double limit solution to $\btsie$.} $\btsie$ it suffices then to prove the existence of a solution to $\btrwsien$ for each fixed $n\in\Ns$, under the condition \eqref{cnd}, that is uniformly $L^{p}(\Omega,\P)$ bounded on $[0,T]\times \X^d$ for every $T>0$ and every $p\ge2$.  We establish this existence via our K-martingale approach, using \thmref{kmarteqaux} which we just proved in \secref{kmart2trnctd}.

First, the following proposition summarizes the results in this case for the BTRW SIEs spatial lattice scale.

\bpr[Existence for BTRW SIEs with non-Lipschitz $a$]\lbl{ewlip}
Assume the conditions \eqref{cnd} hold.  Then,
\begin{enumerate}\renewcommand{\labelenumi}{$($\alph{enumi}$)$}
\item For every $(n,l)\in\Ns\times\N$ and for every $p\ge2$, there exists an $L^p$-bounded solution $\uttxnlf$ to the truncated BTRW SIE
\eqref{trnctdsie} on $\T\times\Xnd$.  Moreover, if we linearly interpolate $\uttxnlf$ in space; then, with probability one,
the continuous map $(t,x)\mapsto\uttxnlf$ is locally $\gamma_t$-H\"older continuous in time with
$\gamma_t\in\lpa0,\tf{4-d}{8}\rpa$ for $1\le d\le3$.
\item  For any fixed $n\in\Ns$, the sequence $\lbr\uttxnlf\rbr_{l\in\N}$ of linearly-interpolated solutions in $(a)$
has a subsequential weak limit $\utn$ in $\C(\T\times\Rd;\R)$.  We thus have a limit solution $\utn$
to $\btrwsien$, and $\utn$ is locally $\gamma_t$-H\"older continuous in time with
$\gamma_t\in\lpa0,\tf{4-d}{8}\rpa$ for $1\le d\le3$.
\end{enumerate}
\epr
\bpf
\begin{enumerate}\renewcommand{\labelenumi}{$($\alph{enumi}$)$}
\item  First, recall that the deterministic term $\uttxD$ in \eqref{trnctdsie} is completely determined by $\un$.
Moreover, under the conditions in \eqref{cnd} on $\un$, $\uttxD$ is clearly bounded and it is smooth in time as
in \remref{smuthdet}.  Fix an arbitrary $T>0$, and let $\T=[0,T]$.  We now prove the existence of a family of adapted processes $\lbr\tilde{X}^\tau_{n,l}\rbr_{\tau\in\T}$ satisfying our K-martingale problem
\eqref{kmart}, which by \thmref{kmarteqaux} implies the existence of a solution to the $l$-truncated BTRW SIE
\eqref{trnctdsie} on $\T\times\Xnd$.  On a probability space $\OFFtP$ we prepare  a family of $r$-independent BMs $\lbr\wtyn\rbr_{y\in\Xnldf}$.  For each $\tau\in\T$ and each $i=1,2,\ldots$ define a continuous process $X_{n,l,i}^\tau$ on $[0,\tau]\times\Xnd$
inductively for $k/2^i\le t\le((k+1)/2^i)\wedge\tau$
$(k=0,1,2,\ldots)$ as follows: $X_{n,l,i}^{x,\tau}(0)=\unx$ $(x\in\Xnd)$ and if $X_{n,l,i}^{x,\tau}(t)$ is defined for $t\le k/2^i$, then we define
$X_{n,l,i}^{x,\tau}(t)$ for $k/2^i\le t\le((k+1)/2^i)\wedge\tau$, by
\beq\lbl{indproc}\bsp
&X_{n,l,i}^{x,\tau}(t)\\&=\bc
\ds X_{n,l,i}^{x,\tau}\lpa\tf{k}{2^{i}}\rpa+\ds\sum_{y\in\Xnldf}\kappa^{x,y}_{\h,\frac{k}{2^i},\tau}\lpa X_{n,l,i}^{y,\tau}(\tfrac{k}{2^i})\rpa\lpa\Delta_{t,\tf{k}{2^{i}}}W_{n}^{y}\rpa\\ +\lbk\uttxnD-\tilde{U}_{n,D}^{x}\lpa\frac{k}{2^{i}}\rpa \rbk; &x\in\Xnldf, \\
\ds\uttxnD;&x\in\Xnd\setminus\Xnldf,
\ec
\end{split}
\eeq
where $\Delta_{t,\tf{k}{2^{i}}}W_{n}^{y}=\wtyn-W_n^y(\tfrac{k}{2^i})$.   Clearly, $X_{n,l,i}^\tau$ is the solution to the equation
\beq\lbl{indsie}\bsp
&X^{x,\tau}(t)\\&=
\bc
\ds\sum_{y\in\Xnldf}\int_0^t\kappa^{x,y}_{\h,\phi_i(s),\tau}\lpa X^{y,\tau}(\phi_i(s))\rpa{d\wsyn}
+\uttxnD;&x\in\Xnldf,\\
\ds\uttxnD;&x\in\Xnd\setminus\Xnldf
\ec
\end{split}
\eeq
with $X^{x,\tau}(0)=\unx$, where $\phi_i(t)=k/2^i$ for $k/2^i\le t<(k+1)/2^i\wedge\tau$ $(k=0,1,2,\ldots)$.

Now, for $q\ge1$, let $M_{q,l,i}^\tau(t)=\sup_{x\in\Xnd}\E\lab X_{n,l,i}^{x,\tau}(t)\rab^{2q}$.  By the boundedness of
$\uttxnD$ over the whole infinite lattice $\Xnd$, we have
\beq\lbl{trnctdlpbdd}
\bsp
M_{q,l,i}^\tau(t)\le C+\sup_{x\in\Xnldf}\E\lab X_{n,l,i}^{x,\tau}(t)\rab^{2q}
\end{split}
\eeq
Then, replacing $\Xnd$ by $\Xnldf$ and following the same steps as in the proof of \propref{Expbd}, we get that
\beq\lbl{trnctdulpbdd}
\sup_{\tau\in\T}\sup_{t\in[0,\tau]}M_{q,l,i}^\tau(t)\le C,\quad1\le d\le3,
\eeq
where, here and in the remainder of the proof, the constant $C$ depends only on $q$, $\max_x|u_0(x)|$, the spatial dimension
$1\le d\le3$, and $T$ but may change its value from one line to the next.  Remembering that $\h\searrow0$ as $n\nearrow\infty$ and $n\in\Ns$,
the independence in $l$ is trivially seen since \lemref{2ndQinequality} implies
$$\sum_{y\in\Xnldf} \lbk\Ktxyn\rbk^2\le\sum_{y\in\Xnd} \lbk\Ktxyn\rbk^2\le \df{C}{t^{d/4}};\forall1\le d\le3, l\in\N$$
Similarly, letting $X_{n,l,i,R}^{x,\tau}$ denote the random part of $X_{n,l,i}^{x,\tau}$ on the truncated lattice $\Xnldf$, using
\eqref{trnctdulpbdd}, and repeating the arguments in \lemref{sptdiff} and \lemref{tmpdiff}---replacing $\Xnd$ by $\Xnldf$ and
noting that \lemref{4thQinequality} and \lemref{3rdQinequality} hold on $\Xnldf$---we obtain
\beq\lbl{trnctdspttmpdiff}
\bsp
&\E\lab X_{n,l,i,R}^{x,\tau_{1}}(t) - X_{n,l,i,R}^{y,\tau_{1}}(t)\rab^{2q}+\E\lab X_{n,l,i,R}^{x,\tau_{2}}(t) - X_{n,l,i,R}^{y,\tau_{2}}(t)\rab^{2q}\\&\le C_{d} |x-y|^{2q\alpha_{d}};\ \alpha_{d}\in I_{d},\\
&\E\lab X_{n,l,i,R}^{x,\tau_{1}}(t) - X_{n,l,i,R}^{x,\tau_{1}}(r)\rab^{2q}+\E\lab X_{n,l,i,R}^{x,\tau_{2}}(t) - X_{n,l,i,R}^{x,\tau_{2}}(r)\rab^{2q}\\&\le C\lab t-r\rab^{\tf{(4-d)q}{4}},
\end{split}
\eeq
for all $x,y \in \Xnldf$,  $r,t\in[0,\tau_{1}\wedge\tau_{2}]$, $\tau_{1},\tau_{2}\in\T$, and $1\le d\le3$.   It follows that, for every point $\tau_{1,2}=\lpa\tau_{1},\tau_{2}\rpa\in\T^{2}$, there is a subsequence $\lbr\lpa\tilde{X}_{n,l,i_{m}}^{\tau_{1}},\tilde{X}_{n,l,i_{m}}^{\tau_{2}}\rpa\rbr_{m=1}^\infty$ on a probability space $\OtauotFtauotPtauott$ such that $\lpa\tilde{X}_{n,l,i_{m}}^{\tau_{1}},\tilde{X}_{n,l,i_{m}}^{\tau_{2}}\rpa\overset{\mathscr L}{=} \lpa{X}_{n,l,i_{m}}^{\tau_{1}},{X}_{n,l,i_{m}}^{\tau_{2}}\rpa $ and $$\lpa\tilde{X}_{n,l,i_{m}}^{x,\tau_{1}}(t),\tilde{X}_{n,l,i_{m}}^{x,\tau_{2}}(t)\rpa\longrightarrow \lpa\tilde{X}_{n,l}^{x,\tau_{1}}(t),\tilde{X}_{n,l}^{x,\tau_{2}}(t)\rpa$$  uniformly on compact subsets
of $[0,\tau_{1}\wedge\tau_{2}]\times\Xnd$, as $m\to\infty$ a.s.  Let $\TQ=\T\cap\Q$, where $\Q$ is the set of rationals, and define the product probability space
$$\OFPt:=\lpa{\bigotimes_{\tau_{1,2}\in\TQ^{2}}}\tilde{\Omega}_{\tau_{1,2}},\bigotimes_{\tau_{1,2}\in\TQ^{2}}\tilde{\sF}_{\tau_{1,2}},\bigotimes_{\tau_{1,2}\in\TQ^{2}}\tilde{\P}_{\tau_{1,2}}\rpa.$$
If $s<t$, then for every $ f\in\C_b^2(\Rs;\R)$, $\tau_1,\tau_2\in\TQ\setminus\lbr0\rbr$,
$t\in[0,\tau_1\wedge\tau_2]$, $x_1,x_2\in\Xnldf$, and for every bounded continuous $F:\C\lpa\Rp;\Rs\rpa\to\R$
that is ${\mathscr{B}}_s\lpa\C\lpa\Rp;\Rs\rpa\rpa:=\sigma\lpa z(r); 0\le r\le s\rpa$-measurable function, we have
\beq\lbl{mmart}
\bsp
&\EPt\lbk\lbr f(\tilde{Z}_{n,l}^{x_{1,2},\tau_{1,2}}(t))-f(\tilde{Z}_{n,l}^{x_{1,2},\tau_{1,2}}(s))\right.\right.
\\&-\left.\left.\int_s^t \lpa\sAKtauot f\rpa(r,x_1,x_2,\tilde{Z}_{n,l}^{\tau_{1,2}})dr\rbr F\lpa\tilde{Z}_{n,l}^{x_{1,2},\tau_{1,2}}(\cdot)\rpa\rbk \\
&=\lim_{m\to\infty}\EPt\lbk\lbr f(\tilde{Z}_{n,l,i_m}^{x_{1,2},\tau_{1,2}}(t))-f(\tilde{Z}_{n,l,i_m}^{x_{1,2},\tau_{1,2}}(s))\right.\right.
\\&-\left.\left.\int_s^t \lpa\sAKtauotim f\rpa(r,x_1,x_2,\tilde{Z}_{n,l,i_m}^{\tau_{1,2}})dr\rbr F\lpa\tilde{Z}_{n,l,i_m}^{x_{1,2},\tau_{1,2}}(\cdot)\rpa\rbk=0,
\end{split}
\eeq
where, by a standard localization argument, we have assumed that $a$ is also bounded; and where $\tilde{Z}_{n,l}^{\tau_{1,2}}$ and $\tilde{Z}_{n,l,i_m}^{\tau_{1,2}}$ are obtained from the definition of
${Z}_{n,l}^{\tau_{1,2}}$ in \eqref{Z} by replacing $X_{n,l}^{\tau_j}$ by $\tilde{X}_{n,l}^{\tau_j}$ and $\tilde{X}_{n,l,i_m}^{\tau_j}$, $j=1,2$,
respectively.  The operator $\sAKtauotim$ is obtained from $\sAKtauot$ by replacing $\Upsilon^{x_{i,j},y}_{\h,t,\tau_{i,j}}\lpa u^y(t)\rpa$ in \eqref{siegen} by $\Upsilon^{x_{i,j},y}_{\h,\phi_{i_m}(t),\tau_{i,j}}\lpa u^y(\phi_{i_m}(t))\rpa$.  Also, obviously, for any $\tau\in\TQ$ and $t\in[0,\tau]$
\beq\lbl{mmartd}
\tilde{X}_{n,l}^{x,\tau}(t)=\lim_{m\to\infty}\tilde{X}_{n,l,i_m}^{x,\tau}(t)=\uttxnD; \qquad x\in\Xnd\setminus\Xnldf, \mbox{ a.s.~}\Pt.
\eeq
It follows from \eqref{mmart} and \eqref{mmartd} that  $\lbr\tilde{X}_{n,l}^{\tau}\rbr_{\tau\in\TQ}$ satisfies the K-martingale problem \eqref{kmart} with respect to the filtration $\{\tilde{\sFt}\}$, with
$$\tilde{\sFt}=\bigcap_{\epsilon>0}\sigma\lbr \tilde{X}_{n,l}^{x,\tau}(u);u\le (t+\epsilon)\wedge\tau, \tau\in\TQ\cap(t,T]\rbr.$$
Thus, by \thmref{kmarteqaux}, with $\tau\in\Rp$ replaced by $\tau\in\TQ$, there is a solution $\uttxnlf$ to the $l$-truncated BTRW SIE \eqref{trnctdsie} on $\TQ\times\Xnd$.  Use continuous extension in time of $\uttxnlf$ to extend its definition to $\T\times\Xnd$, and denote the extension also by $\uttxnlf$.  Clearly $\uttxnlf$ solves the $l$-truncated BTRW SIE \eqref{trnctdsie} on $\T\times\Xnd$.

Now, for $q\ge1$, let $M_{q,l}(t)=\sup_{x\in\Xnd}\E\lab\uttxnlf\rab^{2q}$.  As above, the boundedness of $\uttxnD$, implies
\beq\lbl{trnctdlpbdd2}
\bsp
M_{q,l}(t)\le C+\sup_{x\in\Xnldf}\E\lab\uttxnlf\rab^{2q}.
\end{split}
\eeq
Then, replacing $\Xnd$ by $\Xnldf$ and following the same steps as in the proof of \propref{Expbd}, we get that
\beq\lbl{trnctdulpbdd2}
M_{q,l}(t)\le C,\quad\forall t\in\T\mbox{ and }1\le d\le3.
\eeq

Similarly, letting $\uttxnlfR$ denote the random part of $\uttxnlf$ on the truncated lattice $\Xnldf$, using
\eqref{trnctdulpbdd2}, and repeating the arguments in \lemref{sptdiff} and \lemref{tmpdiff}---replacing $\Xnd$ by $\Xnldf$ and
noting that the inequalities in \lemref{4thQinequality} and \lemref{3rdQinequality} trivially hold if we replace $\Xnd$ by $\Xnldf$---we
obtain
\beq\lbl{trnctdspttmpdiff2}
\bsp
&\E\lab\uttxnlfR - \uttynlfR\rab^{2q}\le C_{d} |x-y|^{2q\alpha_{d}};\ \alpha_{d}\in I_{d},\\
&\E\lab\uttxnlfR - \utrxnlfR\rab^{2q}\le C\lab t-r\rab^{\tf{(4-d)q}{4}},
\end{split}
\eeq
for all $x,y \in \Xnldf$, $r,t\in\T$, and $1\le d\le3$.   By \remref{smuthdet},
$\uttxnD$ is differentiable in $t$.  So, linearly interpolating $\uttxnlf$ in space and using  \eqref{trnctdspttmpdiff2} and
arguing as in the proof of part (a) of \lemref{regtight}, we get that the continuous map $(t,x)\mapsto\uttxnlf$
is locally $\gamma_t$-H\"older continuous in time with $\gamma_t\in\lpa0,\tf{4-d}{8}\rpa$ for $1\le d\le3$.
\item  Clearly, $\uttxnD$ in \eqref{trnctdsie} is the same for every $l$, so it is enough to show convergence of the random part $\uttxnlfR$.
Using \eqref{trnctdspttmpdiff2} we get tightness for $\lbr\uttxnlfR\rbr_l$ and consequently a subsequential weak limit $\utn$,
which  is our limit solution for $\btrwsien$.  For the regularity assertion, $\uttxnD$ is smooth and bounded as noted above.  So,
using \eqref{trnctdulpbdd2} and \eqref{trnctdspttmpdiff2}, and imitating the argument in the proof of part (b) of \lemref{regtight}
(remembering that here we are taking the limit as $l\to\infty$); we
get the desired $L^p$ boundedness for $\utn$ as in \propref{Expbd} and the spatial and temporal moments bounds in
\lemref{sptdiff} and \lemref{tmpdiff}\beq\lbl{est}
\bc
\E\lab\uttxn\rab^{2q}\le C&\\
\E\lab\uttxnR - \uttynR\rab^{2q}\le C_{d} |x-y|^{2q\alpha_{d}};\ \alpha_{d}\in I_{d},&\\
\E\lab\uttxnR - \utrxnR\rab^{2q}\le C\lab t-r\rab^{\tf{(4-d)q}{4}},
\ec
\eeq
for $(t,x,n)\in\T\times\Xnd\times\Ns$ and for $d=1,2,3$ and $q\ge1$ and the desired H\"older regularity follows.
\end{enumerate}
The proof is complete.
\epf
We now get \thmref{mainthm2} for $\btsie$ as the following corollary.
\bcr \thmref{mainthm2} holds.
\ecr
\bpf
The desired conclusion follows upon using the argument in the proof of part (b) of \lemref{regtight} along with \defnref{limitsolns} and the $L^p $-boundedness and the spatial and temporal moments bounds for $\lbr \utn\rbr_n$ that we got in \eqref{est} above.
\epf
\subsection{BTBM SIEs and fourth order parametrized BTBM SPDEs link}\lbl{linkpf}
The proof of \lemref{tsddesie} can be easily handled using an application of It\^o's formula.  We give a simple derivation below.
\begin{proof}[Proof of \lemref{tsddesie}]
  Fix an arbitrary $(n,t,x)\in\N\times\Rp\times\Xnd$.   For any $l\in\N$, let $\Xnlxdf$ be the finite sublattice of $\Xnd$ centered around $x$ and of radius $l$; i.e., $\Xnlxdf={\ds\Xnd\cap{\otimes_{i=1}^d}}[x_i-l,x_i+l]$
 ($x_i$ is the $i$-th coordinate of $x$) and let $r=\#\lbr y;y\in\Xnlxdf\rbr$ and $\lbr y;y\in\Xnlxdf\rbr=\lbr y^{(1)},\ldots,y^{{(r)}}\rbr$.  Let $\utn$ be a continuous (in $s$) solution to $\pbtsdden$; i.e., the semimartingale (in $s$) satisfying \eqref{psdde}.  Let
\beqs
\bsp
 F_l(t-s,x,\utstxdn)&:=\sum_{y\in\Xnlxdf}\Ktsxyn\utstxyn=\sum_{k=1}^{r}\Ktsxykn\tilde{U}_{n}^{x,y^{(k)}}(s,t)\\
\end{split}
\eeqs
We denote by $\pa_iF_l$ the first derivative of $F_l$ in the $i$-th variable, $i=1,2$, with $\nabla_3F_l$ being the $r$-dimensional gradient vector of $F_l$ in the third argument $\utn$ at all $y^{(1)},\ldots, y^{(r)}\in\Xnlxdf$; i.e., formally,
$$\nabla_3F_l(t-s,x,\utstxdn):=\lpa\df{\pa F_l}{\pa \tilde{U}^{x,y^{(1)}}_n(s,t)},\ldots,
\df{\pa F_l}{\pa \tilde{U}^{x,y^{(r)}}_n(s,t)}\rpa.$$
Applying It\^o's formula to $F_{l}$, and remembering that $X_s^{(1)}=t-s$ is of bounded variation with $dX_s^{(1)}=-ds$,
$X_s^{(2)}=x$ and $dX_s^{(2)}=0$, $X_{s}^{(k)}=\tilde{U}_{n}^{x,y^{(k-2)}}(s,t)$ for $k=3,\ldots, 2+r$, $\lbr\wsyn\rbr_{y\in\Xnd}$ is a collection of independent BMs, all second partials for $k=3,\ldots, 2+r$ are zero at every $y\in\Xnlxdf$---since $F_l$ is linear in $\lbr \tilde{U}_{n}^{x,y^{(k-2)}}(s,t)\rbr_{k=3}^{2+r}$, and all the differentials below are in $s$, we get
\beqs
\bsp
&\uttxn-\sum_{y\in\Xnlxdf}\Ktxyn\uny= F_l(0,x,\utttxdn)-F_l(t,x,\un)
\\&=-\int_0^t\pa_1F_l(t-s,x,\utstxdn)ds
+\int_0^t\lpa\nabla_3F_l(t-s,x,\utstxdn),d\utstxdn\rpa
\\&=-\sum_{y\in\Xnlxdf}\int_0^t\lbk\df{\partial}{\partial(t-s)}\Ktsxyn\rbk\utstxyn ds
+\sum_{y\in\Xnlxdf}\int_0^t\Ktsxyn d\utstxyn
\end{split}
\eeqs
Using \lemref{fodde} and \eqref{dbtrwsdde}, we see the above is
\beqs
\bsp
&=-\sum_{y\in\Xnlxdf}\int_0^t\lbk\df{\Dn\Kzxyn}{\sqrt{8\pi (t-s)}}+\df18\Dns\Ktsxyn\rbk\utstxyn ds
\\&+\sum_{y\in\Xnlxdf}\int_0^t\Ktsxyn\lbk\df{\Dn\utstxyn\big|_{y=x}}{\sqrt{8\pi(t- s)}}+\df18\Dns\utstxyn\rbk ds
\\&+\sum_{y\in\Xnlxdf}\int_0^t\Ktsxyn\lbk a(\utsyn)\df{d\wsyn}{\sdh}\rbk
\end{split}
\eeqs
Taking the limit as $l\to\infty$, using the conditions on $\tilde{U}$, we get
\beqs\lbl{ito1}
\bsp
&\uttxn-\sum_{y\in\Xnd}\Ktxyn\uny
=\sum_{y\in\Xnd}\int_0^t\Ktsxyn\lbk a(\utsyn)\df{d\wsyn}{\sdh}\rbk
\\&-\sum_{y\in\Xnd}\int_0^t\lbk\df{\Dn\Kzxyn}{\sqrt{8\pi (t-s)}}+\df18\Dns\Ktsxyn\rbk\utstxyn ds
\\&+\sum_{y\in\Xnd}\int_0^t\Ktsxyn\lbk\df{\Dn\utstxyn\big|_{y=x}}{\sqrt{8\pi(t- s)}}+\df18\Dns\utstxyn\rbk ds
\end{split}
\eeqs
Now, remember that $\sum_{y\in\Xnd}\Ktsxyn\equiv1$ for every $(s,t,x)\in[0,t]\times\Rp\times\Xnd$,
that the term $\Dn\utstxyn\big|_{y=x}/\sqrt{8\pi(t- s)}$ does not depend on $y$, and use summation by parts four
times to get
\beqs\lbl{ito2}
\bsp
&\uttxn-\sum_{y\in\Xnd}\Ktxyn\uny=\sum_{y\in\Xnd}\int_0^t\Ktsxyn a(\utsyn)\df{d\wsyn}{\sdh}
\\&-\int_0^t\df{\Dn\utstxyn\big|_{y=x}}{\sqrt{8\pi (t-s)}}ds
-\sum_{y\in\Xnd}\int_0^t\df18\Dns\Ktsxyn\utstxyn ds
\\&+\int_0^t\df{\Dn\utstxyn\big|_{y=x}}{\sqrt{8\pi(t- s)}}ds
+\sum_{y\in\Xnd}\int_0^t\df18\Dns\Ktsxyn\utstxyn ds
\\&=\sum_{y\in\Xnd}\int_0^t\Ktsxyn a(\utsyn)\df{d\wsyn}{\sdh},
\end{split}
\eeqs
which is what we wanted to show.
\end{proof}
\section{Conclusions}\lbl{conc}
We considered the multiplicative noise case for our recently introduced (see \cite{Abtpspde})
 fourth order BTBM stochastic integral equation $\btsie$ in \eqref{btpsol} on $\Rp\times\Rd$,
 under both Lipschitz and less than Lipschitz conditions on the diffusion coefficient $a$.
 We showed striking spatio-temporal dimension-dependent H\"older regularity for such SIEs that are not only real-valued up to spatial dimension $d=3$, but---even more impressively---are spatially nearly locally Lipschitz, and hence roughly twice as smooth in space as the Brownian sheet of the driving noise in $d=1,2$.    This gives, for the first time, an example of a kernel---that is also the density of an interesting stochastic process---that is able to regularize a space-time white noise driven equation so that its solutions are pushed beyond the traditional nearly H\"older $1/2$ spatial regularity.   Of course this contrasts sharply with  second order reaction-diffusion (RD) SPDEs driven by space-time white noise whose counterpart real-valued solutions are confined to the case $d=1$.  This contrast is summarized in Table \ref{tabl} below.
\begin{table}[h]
\bt[t]{|c|c|c|c|c|}\hline
Spatial Dimension & \multicolumn{2}{c|}{Real-Valued Solution}
&\multicolumn{2}{c|}{H\"older Exponent (time, space)}
\\
 \cline{2-5}
 & RD SPDEs & BTBM SIEs & RD SPDEs & BTBM SIEs\\
 \hline
$d=1$ & Yes & Yes & $\lpa\tf14^{-}, \tf12^{-}\rpa$ & $\lpa\tf38^{-}, 1^{-}\rpa$ \\
$d=2$ & No & Yes & N/A & $\lpa\tf14^{-}, 1^{-}\rpa$ \\
$d=3$ & No & Yes & N/A & $\lpa\tf18^{-},\tf12^{-}\rpa$  \\
\hline
\et
\smallskip\\
\caption{BTBM SIEs vs RD SPDEs}
\lbl{tabl}
\end{table}
For the precise regularity results please consult \thmref{lip} and \thmref{mainthm2} and \remref{scnv}.

 We analyzed our BTBM SIE using both a direct and a numerically-flavored lattice approach similar to  our discretized method for the simpler second order RD SPDEs driven by space-time white noise (see \cite{Asdde1,Asdde2}).  We used the second approach in the case of non-Lipschitz conditions.  In it, we discretize space and formulate solutions to the resulting spatially-discrete stochastic integral equations
 in terms of the density of Brownian-time random walk or BTRW---which we introduce in this article along with the
 general class of Brownian-time chains (BTCs), of which BTRW is a special case.
 As with their continuous counterpart BTPs, which we introduced in \cite{Abtp1,Abtp2},
 BTCs are interesting new processes outside of the current well established
 theory; and we believe they merit further study.  In the course of proving our results, we prove several interesting facts about the BTRW.   These include a connection to fourth order differential-difference equation that is proved in Appendix \ref{appA} and different estimates which lead to the definition of $2$-Brownian-times Brownian motion and $2$-Brownian-times random walk.    We define two notions of solutions to the lattice model: direct solutions and limit solutions (from a finite truncation of the lattice to the whole lattice).  These solutions (both direct and limit) are then used to define two types of BTRW SIEs limit solutions to $\btsie$ (direct limit solutions and double limit solutions), as the size of the lattice mesh shrinks to zero.

To deal with existence under the non-Lipschitz condition on $a$, we introduce our K-martingale approach, which is tailor-made for kernel SIEs as $\btsie$ as well as for other mild formulations of many other different SPDEs.  It is a delicate variant of the well known, and by now classic, martingale problem approach of Stroock and Varadhan for SDEs.  This K-martingale approach starts by constructing an auxiliary problem to a truncated lattice version of \eqref{btpsol}, for which the existence of solutions implies solutions existence for the truncated lattice model.  We then formulate a martingale problem equivalent to the auxiliary problem (the K-martingale problem).
  A key advantage of the K-martingale approach is that it is a unified framework in which the existence and uniqueness\footnote{As we noted earlier, we don't prove uniqueness under less-than-Lipschitz conditions, and we therefore don't discuss further the implications of our K-martingale approach to uniqueness in this article.} of many kernel stochastic integral equations, which are the mild formulation for many SPDEs, may be treated using only variants of the kernel SIE.  This includes SPDEs of different orders (including second and fourth), so long as the corresponding spatially-discretized kernel (or density) satisfies Kolmogorov-type bounds on its temporal and spatial differences.  In essence, what the K-martingale approach implies is that if the kernel in the lattice model is nice enough for the lattice model to converge as the lattice mesh shrinks to zero (under appropriate assumptions on $a$), then it is nice enough to guarantee a solution for the lattice model.  We use it here to prove the existence of BTRW SIEs double limit solutions to \eqref{btpsol} under the conditions \eqref{cnd}, but just as with the Stroock-Varadhan method, it can handle uniqueness as well.
The densities of BTRW and BTBM have a considerable regularizing effect on stochastic kernel equations driven by space-time white noise as compared to the standard Brownian motion or continuous time random walk densities (the usual green functions for second order RD equations and their spatially-discretized versions).  The unconventional memory-preserving fourth order PDEs associated with the BTP density are
highly regular: their solutions are $\mathrm{C}^{1,4}(\Rp\times\Rd,\R)$ for all times and all $d\ge1$, despite the fact that
a positive bi-Laplacian term is part of the equation \cite{Abtp1,Abtp2}.  However, this bi-Laplacian is coupled in
a very specific way---dictated by the BTBM probability density function---with a Laplacian acting on the smooth initial data and whose coefficient grows arbitrarily large as time approaches the initial time (zero) at the rate of $1/\sqrt{8\pi t}$.  One way to understand this smoothing effect is to note that the BTBM density in $\btsie$ is intimately connected---and shares regularity properties---with the kernel associated with our recently introduced imaginary-Brownian-time-Brownian-angle process, which we used to give a solution to a Kuramoto-Sivashinsky-type PDE in \cite{Aks}.   In the stochastic setting of this article, this BTBM density smoothing effect on $\btsie$ is evidenced in a regularity of solutions that is much higher than typical second order space-time white noise driven RD SPDEs.  This regularizing effect is such that we are able to obtain $\gamma$-H\"older continuous solutions to our BTBM SIE $\btsie$ for spatial dimensions $1\le d\le3$.   We show that the H\"older exponent is dimension-dependent with $\gamma\in(0,\tf{4-d}{8})$, $1\le d\le3$.  In addition, we are able to show ultra spatial regularity by showing a nearly local Lipschitz regularity for $d=1,2$, and nearly local H\"older $1/2$ regularity in $d=3$.   To get the smoothing effect of the BTBM density, we prove BTBM and BTRW estimates that enable us to extract a BTRW SIEs weak limit (direct limit in the Lipschitz case and double limits in the non-Lipschitz one) H\"older continuous solution to $\btsie$ in spatial dimensions $d=1,2,3$, in spite of the presence of the driving space-time white noise.  Again, the effective H\"older exponent depends on the space dimension through the expression $(4-d)/8$.  This ultra regularity in multispatial dimensions naturally motivates the study of the variations (temporal and spatial) of BTBM SIEs, which we undertake---among other aspects of BTBM SIEs---in \cite{AX} and a followup article.

Encoding this smoothing effect from our BTBM SIE $\btsie$ into a fourth order SPDE
involving the bi-Laplacian coupled with a Laplacian term requires extra parameters.  We
give what we call the fourth order parametrized SPDE corresponding to $\btsie$, linking
the spatially-discretized BTRW SIE to the diagonals of a parametrized stochastic differential-difference equation
(BTRW PSDDE) on the lattice.

  As we recently started doing for PDEs \cite{Aks}, we adapt the methods presented here and in \cite{Aks} to give an
  entirely new approach---in terms of our Linearized Kuramoto-Sivashinsky process (or imaginary-Brownian-time-Brownian-angle process) and related processes---to study the multi-spatial dimensions SPDEs version of famous fourth order applied mathematics PDEs like the Kuramoto-Sivashinsky (several different versions), the Cahn-Hilliard, and the Swift-Hohenberg PDEs.   We  illustrate this in upcoming papers (\cite{Aksspde,AL,AD}) and planned followup papers.  Traditional semigroup analytical methods alone are not adequate for this since the existence of the KS semigroup in $d>1$ is not settled analytically.   In another direction, we have discovered interesting connections between BTBM and related processes and stochastic fractional PDEs, we address these connections and their consequences in upcoming articles as well.

  Also, SIEs corresponding to other BTP processes we introduced in \cite{Abtp2} may also be studied by adapting
  and generalizing our approach here.    We believe BTPs, their PDEs, their SIEs, and their discretized cousins (the BTCs and their equations) can play a useful role by adding new, currently unavailable, insights and models to the ever growing mathematical finance theory (see \cite{Carr} for an example).    We also hope to explore these aspects in future papers.
\subsection*{Acknowledgement}  I would like to thank two anonymous rederees for their careful reading of the paper and for their positive comments.
\appendix
\section{Proof of BTRW-DDE Connection}\lbl{appA}
In this appendix, we give the proof of \lemref{fodde} linking the density of BTRW to fourth order differential-difference equations.
\begin{proof}[Proof of \lemref{fodde}]
Let $u^x_n(t)=\E \lbk u_0\lpa \S^x_{B,\h}(t)\rpa\rbk$ with $u_0$ as in \eqref{cnd}.
Observe that
\beq\lbl{sg}
\E \lbk u_0\lpa \S^x_{B,\h}(t)\rpa\rbk=2\int_0^\infty\ptzs\E\lbk u_0\lpa S^x_{\h}(s)\rpa\rbk ds
\eeq
where $\ptzs$ the transition density of the one dimensional BM $B(t)$. Differentiating \eqref{sg} with
respect to $t$ and putting the derivative under the integral, which is easily justified
by the dominated convergence theorem, then using the fact that
$\partial\ptzs/\partial t=\tf12\partial^2 \ptzs/\partial s^2$ we have
\beq\lbl{tder}
\bsp
\df{d}{d t}\E \lbk u_0\lpa \S^x_{B,\h}(t)\rpa\rbk&=
2\int_0^\infty\df{\pa}{\pa t}\ptzs\E\lbk u_0\lpa S^x_{\h}(s)\rpa\rbk ds
\\&=\int_0^\infty\df{\pas}{\pa s^2}\ptzs\E\lbk u_0\lpa S^x_{\h}(s)\rpa\rbk ds.
\end{split}
\eeq
Letting $\sTns \unx=\E \un(S^x_{\h}(s))$ be the action of the semigroup $\sTns$
associated with the standard continuous-time symmetric random walk $S^x_{\h}$ on the lattice $\Xnd$.
Then, the generator of $S^x_{\h}$ on $\Xnd$ is given by $\mathscr{A}_n=\Dn/2$.
Alternatively, noting that $\qtxn$ is
the fundamental solution to the deterministic heat equation \eqref{latheat} on the lattice $\Xnd$, we
get
\beqs
\bsp
\df{d}{d s}\sTns \unx&=\df{d}{d s}\E \lbk u_0\lpa S^x_{\h}(s)\rpa\rbk=
\sum_{y\in\Xnd}\uny \df{d}{d s}\qsxyn
\\&=\df12\sum_{y\in\Xnd}\uny\Dn\qsxyn
=\df12\Dn\sum_{y\in\Xnd}\uny\qsxyn
\\&=\df12\Dn\E \lbk u_0\lpa S^x_{\h}(s)\rpa\rbk=\mathscr{A}_n\sTns \unx.
\end{split}
\eeqs
So, we integrate \eqref{tder} by parts twice, and we observe that the boundary terms
always vanish at $\infty$ (as $s\nearrow\infty$) and that we have $(\pa/\pa s)\ptzs=0$ at
$s=0$ but $\ptzz>0$.  This gives us
\beq\lbl{fstp}
\bsp
\df{d}{d t}u^x_n(t)&=\df{d}{d t}\E \lbk u_0\lpa \S^x_{B,\h}(t)\rpa\rbk=
-\int_0^\infty\df{\pa}{\pa s}\ptzs\df{d}{d s}\sTns \unx ds
\\&=\ptzz\sAn\unx+
\int_0^\infty\ptzs\sAns\sTns\unx ds
\\&=\ptzz\sAn\unx+\sAns\int_0^\infty\ptzs\sTns\unx ds
\\&=\df{\Dn\unx}{\sqrt{8\pi t}}+\df18\Dns u^x_n(t).
\end{split}
\eeq
Obviously, $u^x_n(0)=\unx$, and we have proven \eqref{btrwdde}.

Of course, if $u^x_n(t)=\Ktxn$, then $\unx=\Kzxn$ as given in \lemref{fodde},
and we have by \eqref{btrwdsty} and by the steps above
\beqs
\bsp
\df{d}{d t}u^x_n(t)&=\df{d}{d t}\Ktxn=2\int_0^\infty\df{\pa}{\pa t}\ptzs\qsxn ds
=\int_0^\infty\df{\pas}{\pa s^2}\ptzs\qsxn ds
\\&=\df{\Dn\qzxn}{\sqrt{8\pi t}}+\df18\Dns u^x_n(t)
=\df{\Dn\Kzxn}{\sqrt{8\pi t}}+\df18\Dns u^x_n(t).
\end{split}
\eeqs
The proof is complete.
\end{proof}
\section{Limit solutions in the  Lipschitz case}\lbl{appB}
We now state our lattice-limit solution  existence, uniqueness, and regularity  for our BTBM SIE on $\Rp\times\Rd$ under Lipschitz conditions.
\bthm[Lattice-limits solutions: the Lipschitz case]\lbl{latlimlip}
Under the Lipschitz conditions
there exists a unique-in-law direct BTRW SIE weak-limit solution to\break $\btsie$, $U$, such that
$U(t,x)$ is $L^p(\Omega,\P)$-bounded on $\T\times\Rd$ for every $p\ge2$ and $U\in\H^{{\tf{4-d}{8}}^{-},{\lpa\tf{4-d}{2}\wedge 1\rpa}^{-}}(\T\times\Rd;\R)$ for every $1\le d\le3$.
\end{thm}
\thmref{latlimlip} follows as a corollary to the results of \secref{regandtight} combined with the following proposition.
\bpr\lbl{lipprop}
Under the Lipschitz conditions $\eqref{lcnd}$ there exists a unique direct solution to $\btrwsien$, $\utn$, on some filtered probability space $\OFFtP$ that is $L^p(\Omega,\P)$-bounded on $[0,T]\times \X_{n}^d$ for every $T>0$, $p\ge2$, $n\in\Ns, $and $1\le d\le3$.
\epr
The proof of \propref{lipprop} follows the same steps as the non-discretization Picard-type direct proof of the corresponding part in the continuous case in \secref{pf1stmain}, with obvious changes, and we leave the details to the interested reader.

\bcr  \thmref{latlimlip} holds.
\ecr
\bpf  The conclusion follows from \propref{lipprop}, \lemref{sptdiff}, \lemref{tmpdiff}, and \lemref{regtight} (b).
\epf
\brm
With extra work, it is possible to prove the existence of a strong limit solution under Lipschitz conditions.  We plan to address that in a future article.
\erm

\section{BTBM SPDE Kernel Formulation,  brief remarks, and a converse to \lemref{tsddesie}}\lbl{appC}
Our main interest in this paper is in the BTBM SIE $\btsie$, but---for the reader's convenience---we show here that
the spatially-discretized version of the SPDE \eqref{btpspde}   given by
\begin{equation}\label{sdde}
\begin{cases}
\df{d\uttxn}{dt}= \df{\Delta_n\unx }{\sqrt{8\pi t}}+\df{\Delta_n^2\uttxn}{8}
+a(\uttxn) \df{d\wtxn}{\sdh dt};& \cr
\utzeroxn = \unx; &
\end{cases}
\end{equation}
is different from the BTRW SIE $\btrwsien$, and we give its kernel formulation.  The derivation, which is similar to the proof of \lemref{tsddesie}, is now sketched.  The integral form of \eqref{sdde} is of course
\begin{equation}\label{ISDDE}
\begin{split}
\uttxn &=u_0(x)+\int_{0}^t\left[\left(\df{\Dn\unx}{\sqrt{8\pi s}}+\frac{1}{8}\Delta^2_n\utsxn\right)ds+
a(\utsxn)\frac{d\wsxn}{\sdh}\right];\\ &\qquad(t,x)\in\T\times\Xnd,\ n\in\N,\mbox{ a.s. }\P.
\end{split}
\end{equation}
We fix an arbitrary $(n,t,x)\in\N\times\Rp\times\Xnd$. Assume that $\utn$ is a semimartingale (in $t$) satisfying \eqref{ISDDE}.   Applying It\^o's formula and proceeding as in the proof of \lemref{tsddesie} we get
\beqs\lbl{ito2app}
\bsp
&\uttxn-\sum_{y\in\Xnd}\Ktxyn\uny=\int_0^t\sum_{y\in\Xnd}\Ktsxyn\lbk a(\utsyn)\df{d\wsyn}{\sdh}\rbk
\\&-\int_0^t\sum_{y\in\Xnd}\lbk\df{\Dn\Kzxyn}{\sqrt{8\pi (t-s)}}+\df18\Dns\Ktsxyn\rbk\utsyn ds
\\&+\int_0^t\sum_{y\in\Xnd}\Ktsxyn\lbk\df{\Dn\uny}{\sqrt{8\pi s}}+\df18\Dns\utsyn\rbk ds
\\&=\int_0^t\sum_{y\in\Xnd}\Ktsxyn\df{\Dn\uny}{\sqrt{8\pi s}}ds
-\int_0^t\df{\Dn\utsxn}{\sqrt{8\pi (t-s)}}ds
\\&+  \int_0^t\sum_{y\in\Xnd}\Ktsxyn\lbk a(\utsyn)\df{d\wsyn}{\sdh}\rbk
\end{split}
\eeqs
I.e, a solution to the SDDE system in equation $\eqnref{sdde}$
satisfies the Brownian-time random walk (BTRW) kernel (density function) formulation:
\begin{equation}\label{GSDDEeq}
\begin{split}
\uttxn  &=\sum_{y\in\Xnd} \lbk\Ktxyn\uny+\Dn\uny\int_0^t\df{\Ktsxyn}{\sqrt{8\pi s}}ds\rbk
\\&-\int_0^t\df{\Dn\utsxn}{\sqrt{8\pi(t-s)}}ds
+\sum_{y\in\Xnd}\int_0^t \Ktsxyn a(\utsyn)\frac{d\wsyn}{\sdh},
\end{split}
\end{equation}
for $(t,x)\in\Rp\times\Xnd$, where the BTRW density is given by \eqref{btrwdsty}.

Several observations on the special form
of \eqref{GSDDEeq} are in order here.  First, unlike second order SPDEs on the lattice (RD,
Burgers, etc\dots) and unlike other fourth order SPDEs on the lattice (e.g. see \cite{AL,AD} for equations like
KS, CH, SH, etc\ldots) and unlike the BTRW SIEs in \eqref{btrwsie}, the BTRW kernel formulation
\eqref{GSDDEeq}
involves terms (the second and third terms on the r.h.s.~of \eqref{GSDDEeq})
that are simply not there in other standard Green function type formulations.
This is due to the unique form of the BTP PDE and its discretized version, involving the initial data $\un$
in the equation itself.  Second, the third term on the r.h.s.~of \eqref{GSDDEeq} involves the discrete Laplacian of the
solution $\Dn\uttxn$ with no kernel term, which is fine on the lattice $\Xnd$, but the continuous-space Laplacian
$\Delta U(t,x)$ of the solution to the BTBM SPDE $\btspde$ is not defined in the classical sense.
This is why we regard the BTBM SPDE as a degenerate version of the BTBM SIE and its associated parametrized BTBM SPDE.
So, with an eye on the limiting SPDE $\btspde$, one way to start addressing this difficulty is by reformulating
\eqref{GSDDEeq} into a weaker test function formulation which we now give.
Multiplying \eqref{GSDDEeq} by $\h$  and $\xi\in C^2_{c}(\Rd;\R)$, summing over $x$, and
using summation by parts on the third term on the r.h.s.~of \eqref{GSDDEeq}, we obtain
\beq\lbl{tf-btrw-sdde}
\bsp
&\uttxin:=\sum_{x\in\Xnd} \uttxn\xi(x)\hd=
\sum_{x\in\Xnd}\xi(x)\lbr\sum_{y\in\Xnd} \Ktxyn\uny\rbr\hd
\\&+\sum_{x\in\Xnd}\xi(x)\lbr\sum_{y\in\Xnd} \Dn\uny\int_0^t\df{\Ktsxyn}{\sqrt{8\pi s}}ds\rbr\hd
\\&-\sum_{x\in\Xnd}\lbr\int_0^t\df{\utsxn}{\sqrt{8\pi(t-s)}}ds\rbr\Dn\xi(x)\hd
\\&+\sum_{x\in\Xnd}\xi(x)\lbr\sum_{y\in\Xnd}\int_0^t \Ktsxyn a(\utsyn)\frac{d\wsyn}{\sdh}\rbr\hd
\end{split}
\eeq
Here again it is interesting to observe that the form of \eqref{tf-btrw-sdde} is special, for it has
the unusual feature of having a mix of both kernel and test function terms in the same equation,
with the third term on the r.h.s.~of \eqref{tf-btrw-sdde} involving only the test function $\xi$ with no kernel
terms.  It is clear by the derivation of \eqref{tf-btrw-sdde} above that a solution to $\btsdden$
will satisfy \eqref{tf-btrw-sdde} for every $\xi\in C^2_{c}(\Rd;\R)$.

We end with a converse to \lemref{tsddesie}
\blm\lbl{converse}
Assume that for each fixed $(t,x,y)\in\Rp\times\Xnds$ the process $\utn$ is a continuous
semimartingale in $s$ having the form
\beq\lbl{dec}
\utstxyn=\uny+V^{x,y}_n(s,t)+M^y_n(s),
\eeq
where $V_n$ is the process of bounded variation on compacts (in $s$)
and $M_n$ is the local martingale (in $s$) in the decomposition of the semimartingale
$\utn$.  Assume further that $\uttxn:=\utttxxn$ satisfies \eqref{btrwsie}, and that for any fixed pair $(t,x)$, $\E\lab\utstxyn\rab^2\le C$ for all $(s,y)\in[0,t]\times\Xnd$ for some constant $C>0$.  Then $\utn$ satisfies \eqref{psdde}.
\elm
\brm\lbl{momcndonpsddes2}
Again, the moment boundedness condition above is for convenience and may be relaxed.
\erm
\begin{proof}[Sketch of the proof]
Assume that the $(t,x)$-parametrized random field $\utstxyn$ on $[0,t]\times\Rp\times\Xnds$ is a continuous semimartingale in $s$ such that $\uttxn=\utttxxn$ satisfies the BTRW SIE \eqref{btrwsie}, and assume $\utztxyn = \uny$ for all $(t,x,y)\in\Rp\times\Xnds$; then as in the proof of \lemref{tsddesie} we have by It\^o's rule and \lemref{fodde}
\beqs
\bsp
&\int_0^t\sum_{y\in\Xnd}\Ktsxyn\lbk a(\utsyn)\df{d\wsyn}{\sdh}\rbk=\uttxn-\sum_{y\in\Xnd}\Ktxyn\uny
\\&=-\int_0^t\sum_{y\in\Xnd}\lbk\df{\Dn\Kzxyn}{\sqrt{8\pi (t-s)}}+\df18\Dns\Ktsxyn\rbk\utstxyn ds
\\&+\int_0^t\sum_{y\in\Xnd}\Ktsxyn d\utstxyn
 \end{split}
 \eeqs
which is
\beqs
\bsp
\\&=-\int_0^t\sum_{y\in\Xnd}\Ktsxyn\lbk\df{\Dn\utstxyn\big|_{y=x}}{\sqrt{8\pi (t-s)}}
+\df18\Dns\utstxyn \rbk ds
\\&+\int_0^t\sum_{y\in\Xnd}\Ktsxyn dV^{x,y}_n(s,t)+\int_0^t\sum_{y\in\Xnd}\Ktsxyn dM^{y}_n(s)
\end{split}
\eeqs
for every $(t,x)\in\Rp\times\Xnd$.  It is then straightforward to see that
\beq
\bsp
&\int_0^t\sum_{y\in\Xnd}\Ktsxyn\lbk a(\utsyn)\df{d\wsyn}{\sdh}- dM^{y}_n(s) \rbk
\\&=-\int_0^t\sum_{y\in\Xnd}\Ktsxyn\lbk\df{\Dn\utstxyn\big|_{y=x}}{\sqrt{8\pi (t-s)}}
+\df18\Dns\utstxyn \rbk ds
\\&+\int_0^t\sum_{y\in\Xnd}\Ktsxyn dV^{x,y}_n(s,t)
\end{split}
\eeq
and that $dM^{y}_n(s)=a(\utsyn)\tf{d\wsyn}{\sdh}$; i.e.,
\beqs
\bsp
d\utstxyn=\lbk\df{\Dn\utstxyn\big|_{y=x}}{\sqrt{8\pi(t-s)}}+\df18\Dns\utstxyn\rbk ds
+a(\utsyn)\df{d\wsyn}{\sdh}.
\end{split}
\eeqs
We are now done.
\end{proof}
\section{Glossary of frequently used acronyms and notations}\lbl{glossary}
\begin{enumerate}\renewcommand{\labelenumi}{\Roman{enumi}.}
\item {\textbf{Acronyms}}\vspace{2mm}
\bit
\item BM: Brownian motion
\item BTBM: Brownian-time Brownian motion.
\item BTP:  Brownian-time processe.
\item BTC: Brownian-time chain.
\item BTRW: Brownian-time random walk.
\item DDE, SDDE, and PSDDE: Differential difference equation, Stochastic DDE, and Parametrized SDDE.
\item IBTBAP: Imaginary-Brownian-time Brownian-angle process.
\item KS: Kuramoto-Sivashinsky.
\item RW: Random walk.
\item SIE: Stochastic integral equation.
\eit
\vspace{2.5mm}
\item {\textbf{Notations}}\vspace{2mm}
\bit
\item $\qtxyn$:  The $d$-dimensional continuous-time random walk transition density
starting at $x\in\Xnd$ and going to $y\in\Xnd$ in time $t$.
\item $\psxy$: The density of a $d$-dimensional BM.
\item $\ptzs$: The  density of a $1$-dimensional BM, starting at $0$.
\item $\KBtxy$: The kernel or density of a $d$-dimensional Brownian-time Brownian motion.\vspace{0.5 mm}
\item $\Ktxyn$: The kernel or density of a $d$-dimensional Brownian-time random walk on a spatial lattice with step size $\h$ in each of the $d$-dimensions.
\item $\btsie$: The BTBM SIE with diffusion coefficient $a$ and initial function $\un$.
\item $\btrwsien$: The BTRW SIE on the lattice $\Xnd=\h\Zd$ with diffusion coefficient $a$ and initial function $\un$.
\item $\mathrm{C}_{b}^2(\Rd;\R)$ the space of twice continuously differentiable real-valued functions on $\Rd$ with bounded derivatives of all orders $k=0,1,2$.
\eit
\end{enumerate}

\end{document}